\documentclass[journal]{IEEEtran}

\usepackage{cite}     
\usepackage{url}   
\usepackage{amsmath}
\usepackage{amssymb}
\usepackage{graphicx}

\usepackage{parskip}
\usepackage[labelfont=bf,textfont=it]{caption}
\usepackage[marginclue,footnote,silent,draft]{fixme}

\begin{document}
\title{Continuous Fuzzy Transform as Integral Operator}
%
\author{Giuseppe Patan\'e 
\IEEEcompsocitemizethanks{\IEEEcompsocthanksitem G. Patan\'e is with CNR-IMATI, Consiglio Nazionale delle Ricerche, Istituto di Matematica Applicata e Tecnologie Informatiche
Genova, Italy.\protect\\ E-mail: patane@ge.imati.cnr.it}}
%
%
%

\markboth{IEEE TRANSACTIONS ON FUZZY SYSTEMS}%
{}

\maketitle

\begin{abstract}
The Fuzzy transform is ubiquitous in different research fields and applications, such as image and data compression, data mining, knowledge discovery, and the analysis of linguistic expressions. As a generalisation of the Fuzzy transform, we introduce the continuous Fuzzy transform and its inverse, as an integral operator induced by a kernel function. Through the relation between membership functions and integral kernels, we show that the main properties (e.g., continuity, symmetry) of the membership functions are inherited by the continuous Fuzzy transform. Then, the relation between the continuous Fuzzy transform and integral operators is used to introduce a data-driven Fuzzy transform, which encodes intrinsic information (e.g., structure, geometry, sampling density) about the input data. In this way, we avoid coarse fuzzy partitions, which group data into large clusters that do not adapt to their local behaviour, or a too dense fuzzy partition, which generally has cells that are not covered by the data, thus being redundant and resulting in a higher computational cost. To this end, the data-driven membership functions are defined by properly filtering the spectrum of the Laplace-Beltrami operator associated with the input data. Finally, we introduce the space of continuous Fuzzy transforms, which is useful for the comparison of different continuous Fuzzy transforms and for their efficient computation.
\end{abstract}
%
%
%
%
\section{Introduction\label{sec:INTRODUCTION}}
\IEEEPARstart{D}ue to the increasing availability of data, which is supported by ongoing technological advances in acquisition, storage, and processing, several transformations (e.g., the Fourier transform, the Laplace transform, the Fuzzy transform) have been proposed to solve problems that spread from signal analysis to the solution of partial differential equations, from the analysis to the approximation of signals, from fuzzy logic to fuzzy modelling. From a general perspective, a transform is typically defined as a linear operator between functional spaces, and its discretisation reduces to a matrix-vector multiplication. Main examples include the definition of the Fourier and Laplace transforms as integral operators induced by a complex and a real exponential kernel, respectively.

In fuzzy modelling (Sect.~\ref{sec:PREVIOUS-WORK}), the \emph{Fuzzy transform} (F-transform)~\cite{PERFILIEVA2006,PERFILIEVA2016A} maps the space of continuous functions to vectors in~$\mathbb{R}^{n}$, and computations on~ the input functions are then converted into discrete operations on~$\mathbb{R}^{n}$. Viceversa, the \emph{inverse F-transform} converts discrete samples of the input signal to a continuous approximation. The ubiquity of the F-transform is due to its different constructions, which apply linear algebra in vector spaces of finite dimension, residuated lattice, specialised basis functions (e.g., B-splines, Shepard kernels, and Bernstein polynomials), and higher-degree formulations. For instance, the F-transform has been applied to image~\cite{PATERNIAN2017,PERFILIEVA2010,DIMARTINO2008} and data~\cite{SZTYBER2014,ABDELLAL2013,BASHLOVKINA2015,GAETA2014} compression by optimising different aspects, such as bandwidth, allocated memory space, signal curvature, and data reduction~\cite{DIMARTINO2014,PERFILIEVA2017}. The current definition of the F-transform is limited mainly to 1D signals and 2D data organised as a regular grid (e.g., 2D images). This aspect limits the potential application of the F-transform to arbitrary data in terms of dimensionality and structure, and is mainly due to the difficulty to define membership functions on arbitrary data. 

\textbf{Overview and contribution}
In this paper, we introduce the \emph{continuous F-transform} \mbox{$\mathcal{F}:\mathcal{L}^{2}(\Omega)\rightarrow\mathcal{C}^{0}(\Omega)$}, \mbox{$\mathcal{F}f(\mathbf{p})=\int_{\Omega}K(\mathbf{p},\mathbf{q})f(\mathbf{q})d\mathbf{q}$}, as an integral operator induced by a kernel \mbox{$K:\Omega\times\Omega\rightarrow\mathbb{R}$}, and as a generalisation of the F-transform (Sect.~\ref{sec:CONTINUOUS-FT}). Here, \mbox{$\mathcal{L}^{2}(\Omega)$} and \mbox{$\mathcal{C}^{0}(\Omega)$} are the spaces of square-integrable and continuous functions defined on~$\Omega$, respectively. To this end, we exploit the mutual relation between membership functions and kernels, i.e., (i) any membership function induces an integral kernel through a normalisation of its values by the mean of the membership function over~$\Omega$, and (ii) any kernel generates a family of membership functions. For the definition of the continuous F-transform, we consider different classes of membership functions (Sect.~\ref{sec:INT-OPERATOR-SP}), which are defined analytically as polynomials and radial kernels, or generated as a tensor product, a linear combination, a pointwise product, or as the limit of a sequence of kernels.

According to the relation between membership functions and integral kernels, we show that the main properties (e.g., square integrability, continuity, symmetry, positiveness) of the membership functions are inherited by the continuous F-transform. Furthermore, the continuous F-transform \mbox{$\mathcal{F}f$} interpolates the values \mbox{$\mathbf{F}_{n}:=(F_{i})_{i=1}^{n}$} of the F-transform associated with the set \mbox{$(f(\mathbf{q}_{i}))_{i=1}^{s}$} of~$f$-values sampled at a discrete set of points, i.e., \mbox{$\mathcal{F}f(\mathbf{p}_{i})=F_{i}$}. Indeed, the continuous F-transform can be interpreted as a generalisation of the F-transform. We further study the generalisation properties of the continuous F-transform through restriction and out-of-sample operators, based on meshless approximations with radial kernels.

The relation between the continuous F-transform and integral operators allows us to introduce a \emph{data-driven F-transform} (Sect.~\ref{sec:CONT-DATA-FT-SHORT}) through the definition of \emph{data-driven membership functions}, which encode \emph{intrinsic} information (e.g., structure, geometry, sampling density) about the input data. In this way, we avoid coarse fuzzy partitions, which group data into large clusters that do not adapt to their local behaviour, or a too dense fuzzy partition, which generally has cells that are not covered by the data, thus being redundant and resulting in a higher computational cost. To this end, the data-driven membership functions are defined by properly filtering the spectrum of the Laplace-Beltrami operator associated with the input data. The aforementioned properties of the membership functions are then inherited by the continuous F-transform and are important in case of structured (e.g., regular, irregular), sparse, or time-depending data. Indeed, the proposed generalisation allows us to define and efficiently compute the F-transform on arbitrary data by properly encoding their properties in the membership functions.

The representation of the inverse continuous F-transform is derived according  to the structure of the underlying space as a Hilbert or a Reproducing Kernel Hilbert Space. Applying the Mercer theorem for integral operators, we represent the continuous F-transform in terms of its spectrum (i.e., the eigenvalues/eigenfunctions of the integral operator) and express the inverse F-transform in terms of the pseudo-inverse of integral operators (Sect.~\ref{sec:SPECTRAL-REPRESENTATION-FT}).

Through integral operators, we introduce the \emph{space of continuous F-transforms}, which is endowed with a Hilbert Space structure. This space is useful to compare the discrete and continuous F-transforms, according to the underlying scalar product and the corresponding norm, and to approximate a given continuous F-transform in order to make its computation more efficient and numerically stable. In the discrete setting, we obtain analogous relations and reduce the evaluation of the continuous F-transform to numerical linear algebra. Finally, we discuss experimental results (Sect.~\ref{sec:DISCUSSION}) and possible extensions of the proposed approach (Sect.~\ref{sec:CONCLUSION}). 

Main \emph{contributions} of the paper are: (i) the generalisation of the F-transform to continuous signals, which can be applied to arbitrary data, in terms of dimensionality and structure; (ii) a characterisation of its properties through the theory of integral operators; (iii) the definition of the space of continuous F-transforms for the modelling and comparison of F-transforms induced by different kernels; (iv) the definition of data-driven membership functions and continuous F-transforms, which encode and adapt to the local properties of the input data, in terms of geometric features and density.
\begin{figure}[t]
\centering
\includegraphics[height=110pt]{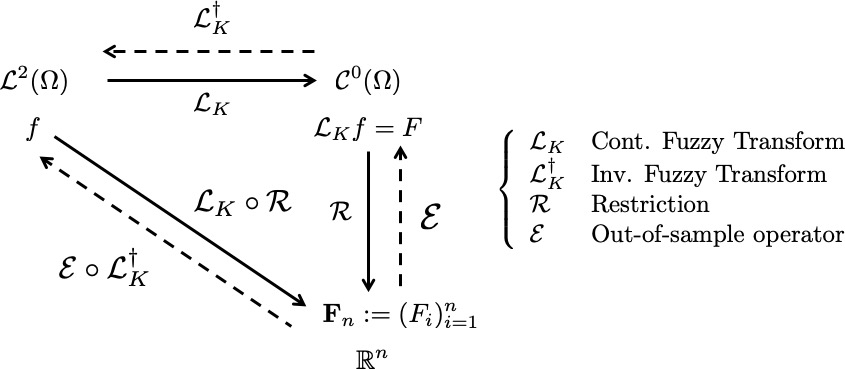}
\caption{Commutative diagram for the continuous and discrete F-transform and their inverse.\label{fig:DIAGRAM}}
\end{figure}
\section{Previous work\label{sec:PREVIOUS-WORK}}
We briefly review previous work on the F-transform~\cite{PERFILIEVA2006,PERFILIEVA2008} (Sect.~\ref{sec:PREVIOUS-WORK-DFT}) and integral operators~\cite{SCHOELKOPF02} (Sect.~\ref{sec:INTEGRAL-OPERATOR}).

\subsection{Discrete F-transform\label{sec:PREVIOUS-WORK-DFT}}
In fuzzy modelling, the \emph{F-transform} \mbox{$\mathcal{F}:\mathcal{C}^{0}(\Omega)\rightarrow\mathbb{R}^{n}$} provides a relation between the space \mbox{$\mathcal{C}^{0}(\Omega)$} of continuous functions defined on a domain~$\Omega$ and~$\mathbb{R}^{n}$. Let \mbox{$\mathcal{I}:=\{\Omega_{i}\}_{i=1}^{n}$} be a partition of~$\Omega$ and \mbox{$\mathcal{P}:=\{\mathbf{p}_{i}\}_{i=1}^{n}$} a set of points such that \mbox{$\mathbf{p}_{i}\in\Omega$}, \mbox{$i=1,\ldots,n$}. A family of functions \mbox{$\mathcal{A}:=\{A_{i}:\Omega\rightarrow [0,1]\}_{i=1}^{n}$} is a \emph{fuzzy partition} of~$\Omega$ if the following properties hold for each~$i$
\begin{itemize}
\item$A_{i}(\mathbf{p})\neq 0$, \mbox{$\mathbf{p}\in\Omega_{i}$}, and \mbox{$A_{i}(\mathbf{p}_{i})=1$}; 
\item$A_{i}$ is continuous and has its unique maximum at~$\mathbf{p}_{i}$;
\item for all \mbox{$\mathbf{p}\in\Omega$}, \mbox{$\sum_{i=1}^{n}A_{i}(\mathbf{p})=1$}.
\end{itemize}
Under these assumptions, the \emph{F-transform}~\cite{PERFILIEVA2006,PERFILIEVA2008,PERFILIEVA2010} of a function \mbox{$f:\Omega\subseteq\mathbb{R}^{d}\rightarrow\mathbb{R}$} is defined as the array \mbox{$\mathbf{F}_{n}:=(F_{i})_{i=1}^{n}\in\mathbb{R}^{n}$} with components
\begin{equation}\label{eq:FUZZY-TRANSFORM}
F_{i}:=\frac{\int_{\Omega} A_{i}(\mathbf{p})f(\mathbf{p})d\mathbf{p}}{\int_{\Omega} A_{i}(\mathbf{p})d\mathbf{p}},\qquad i=1,\ldots,n.
\end{equation}
Among the properties of the F-transform, we mention the linearity with respect to the input function and the least-squares property, which guarantees that the~$i$-th component of~$\mathbf{F}_{n}$ minimises the quadratic least-squares error \mbox{$\Phi_{i}(t):=\int A_{i}(\mathbf{p})\vert f(\mathbf{p})-t\vert^{2}d\mathbf{p}$}, \mbox{$t\in\mathbb{R}$}, associated with~$A_{i}$.

In real cases, where the function~$f$ is known at a set of points \mbox{$\mathcal{Q}:=\{\mathbf{q}_{i}\}_{i=1}^{s}$}, the definition (\ref{eq:FUZZY-TRANSFORM}) is replaced by the \emph{discrete F-transform} \mbox{$\mathbf{F}_{n}:=(F_{i})_{i=1}^{n}\in\mathbb{R}^{n}$}, whose components are
\begin{equation}\label{eq:DISCRETE-FT}
F_{i}:=\frac{\sum_{j=1}^{s}A_{i}(\mathbf{q}_{j})f(\mathbf{q}_{j})}{\sum_{j=1}^{s}A_{i}(\mathbf{q}_{j})},\qquad
i=1,\ldots,n,\quad s\leq n.
\end{equation}
Generally, the number~$s$ of samples is smaller than or equal to the number~$n$ of membership functions. The discrete F-transform is applied to recover an approximation~$f_{F,n}$ of the function~$f$ underlying the set of values \mbox{$(f(\mathbf{q}_{i}))_{i=1}^{s}$} through the \emph{inverse F-transform}~\cite{PERFILIEVA2006}, which is defined as \mbox{$f_{F,n}(\mathbf{p}):=\sum_{i=1}^{n}F_{i}A_{i}(\mathbf{p})$}, \mbox{$\mathbf{p}\in\mathbb{R}^{d}$}.

The \emph{inverse F-transform} \mbox{$\mathcal{F}_{-1}:\mathbb{R}^{n}\rightarrow\mathcal{C}^{0}(\Omega)$} identifies any vector of~$\mathbb{R}^{n}$ with a continuous map. Even though~$\mathcal{F}_{-1}$ is not the inverse of~$\mathcal{F}$, the inverse F-transform~$f_{F,n}$ approximates~$f$ up to an arbitrary precision~\cite{PERFILIEVA2006} under mild conditions on the input function values. In particular, discrete data can be transformed to a continuous approximation through the inverse F-transform and computations on~$\mathcal{C}^{0}$ are converted in discrete operations on~$\mathbb{R}^{d}$ through the F-transform. Finally~\cite{PERFILIEVA2006}, for any given approximation accuracy~$\epsilon$ there exists a number~$n_{\epsilon}$ of nodes and a set \mbox{$(A_{i}(\cdot))_{i=1}^{n_{\epsilon}}$} of~$n_{\epsilon}$ membership functions such that the discrepancy \mbox{$\|f-f_{F,n_{\epsilon}}\|_{\infty}$} between~$f$ and its inverse F-transform~$f_{F,n_{\epsilon}}$ is lower than~$\epsilon$.

\textbf{Applications of the F-transform}
The ubiquity of the F-transform is due to its different constructions~\cite{STEPNICKA2009,DIMARTINO2008}, which apply linear algebra in vector spaces of finite dimension, residuated lattice, specialised basis functions (e.g., B-splines, Shepard kernels, Bernstein polynomials), and higher-degree formulations~\cite{PERFILIEVA2016A}, which provide a link between the F-transform and approximation schemes~\cite{PATANE2014} or neural networks~\cite{GAETA2014}.

During the years, fuzzy modelling and the F-transform have been applied to a wide range of applications, such as the construction of fuzzy versions of binary morphological operations (e.g., shape detection, edge detection, and clutter removal) for image processing~\cite{SINHA1997} and fuzzy modelling algorithms, which partition the input space according to the correlation among components of sampled data~\cite{KIM1998}. We also mention the definition and characterisation of the Fourier transform by considering the uncertainty of the transformed function~\cite{ BUTKIEWCZ2008}, and a fuzzy-based paradigm for data compression aimed at reducing the computational burden of data analysis in smart grids for 5G applications~\cite{LOLA2017}.

In particular, the F-transform has been applied to image~\cite{PATERNIAN2017,PERFILIEVA2010,DIMARTINO2008} and data~\cite{SZTYBER2014,ABDELLAL2013,BASHLOVKINA2015,GAETA2014} compression by optimising different aspects, such as bandwidth, allocated memory space, signal curvature, and data reduction~\cite{DIMARTINO2014,PERFILIEVA2017}. Further applications include data mining~\cite{AGRAWAL1993,HONG2003,MITRA2002,ZHANG2006}, knowledge discovery~\cite{FAYYAD1996,SHAPIRO200}, and the analysis of linguistic expressions~\cite{NOVAK2008}. According to~\cite{HURTIK2019}, the compression and reconstruction quality of the F-transform is improved by imposing the monotonicity and Lipschitz continuity of functions. To optimise the compression, the F-transform is combined with quantisation~\cite{HELFROUSH2011}, fuzzy edge detection~\cite{GAMBHIR2015}, coding/decoding schemes~\cite{GAMBHIR2017}, JPEG~\cite{HURTIK2017} through a discrete cosine transform~\cite{PENNEBAKER1992,RAO1990}. The first and second order degree F-transforms have been applied to the solution of the Cauchy problem~\cite{KHASTAN2016}, of two-points boundary value problems~\cite{KHASTAN2017} and of Volterra-Fredholm integral equations~\cite{TOMASIELLO2019}. Finally, splines collocation methods~\cite{ALIJANI2020} have been applied to the solution of a system of fuzzy fractional differential equations.

\subsection{Integral operators\label{sec:INTEGRAL-OPERATOR}}
Given a compact domain~$\Omega$ of~$\mathbb{R}^{d}$, let us consider the space \mbox{$\mathcal{L}^{2}(\Omega)$} of square integrable functions defined on~$\Omega$, endowed with the \mbox{$\mathcal{L}^{2}(\Omega)$} scalar product \mbox{$\langle f,g\rangle_{2}:=\int_{\Omega} f(\mathbf{p})g(\mathbf{p})d\mathbf{p}$} and the corresponding norm \mbox{$\| f\|_{2}^{2}:=\int_{\Omega} \vert f(\mathbf{p})\vert^{2}d\mathbf{p}$}. On the space \mbox{$\mathcal{C}^{0}(\Omega)$} of continuous functions defined on~$\Omega$, we consider the \mbox{$\mathcal{L}^{2}(\Omega)$} and the~$\mathcal{L}^{\infty}(\Omega)$-norm \mbox{$\|f\|_{\infty}
:=\max_{\mathbf{p}\in\Omega}\{\vert f(\mathbf{p})\vert\}$}. Given a measurable kernel \mbox{$K:\Omega\times\Omega\rightarrow\mathbb{R}$}, the corresponding \emph{integral operator} \mbox{$\mathcal{L}_{K}:\mathcal{L}^{2}(\Omega)\rightarrow\mathcal{L}^{2}(\Omega)$} is defined as the linear operator \mbox{$(\mathcal{L}_{K}f)(\mathbf{p})
:=\int_{\Omega}K(\mathbf{p},\mathbf{q})f(\mathbf{q})\textrm{d}\mathbf{q}$}. If \mbox{$K(\cdot,\cdot)$} is a square integrable kernel on \mbox{$\Omega\times\Omega$}, then~$\mathcal{L}_{K}$ is a bounded (i.e., continuous) operator and its norm is \mbox{$\|\mathcal{L}_{K}\|=\|K\|_{2}$}.

\section{Continuous F-transform as integral operator\label{sec:CONTINUOUS-FT}}
We introduce the continuous F-transform, its relation with integral operators (Sect.~\ref{sec:CONTINUOUS-FT-INT-OPER}) and previous work (Sect.~\ref{sec:HILBERT-FT}), the properties of the integral operator (Sect.~\ref{sec:CONTINUOUS-FT-PROPERTIES}) and of the continuous F-transform (Sect.~\ref{sec:PROPERTIES-LKF}).
\begin{figure}[t]
\centering
\begin{tabular}{cc|cc}
\multicolumn{2}{c|}{$t=50$}	&\multicolumn{2}{c}{$t=10$}\\
\includegraphics[height=42pt]{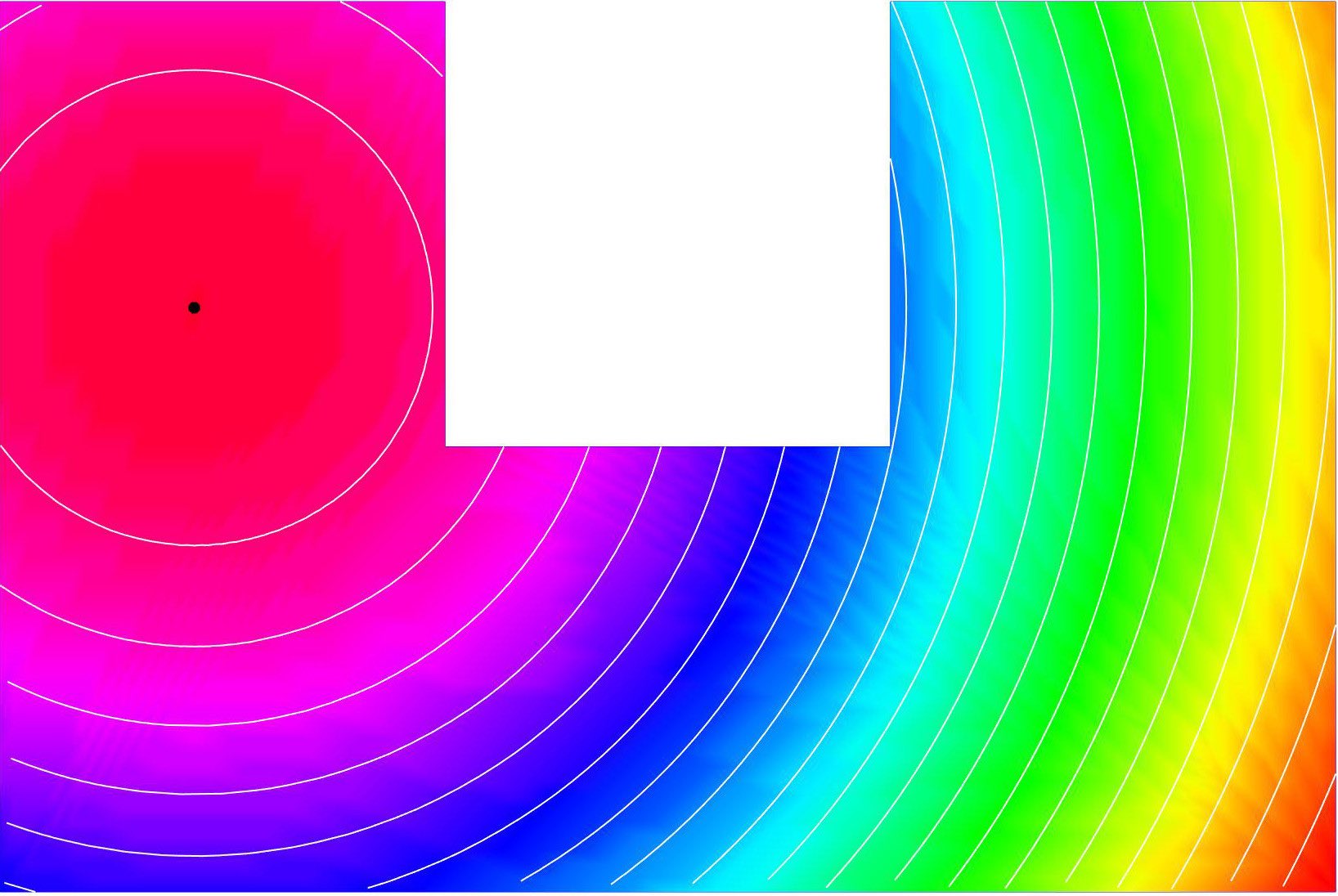}
&\includegraphics[height=42pt]{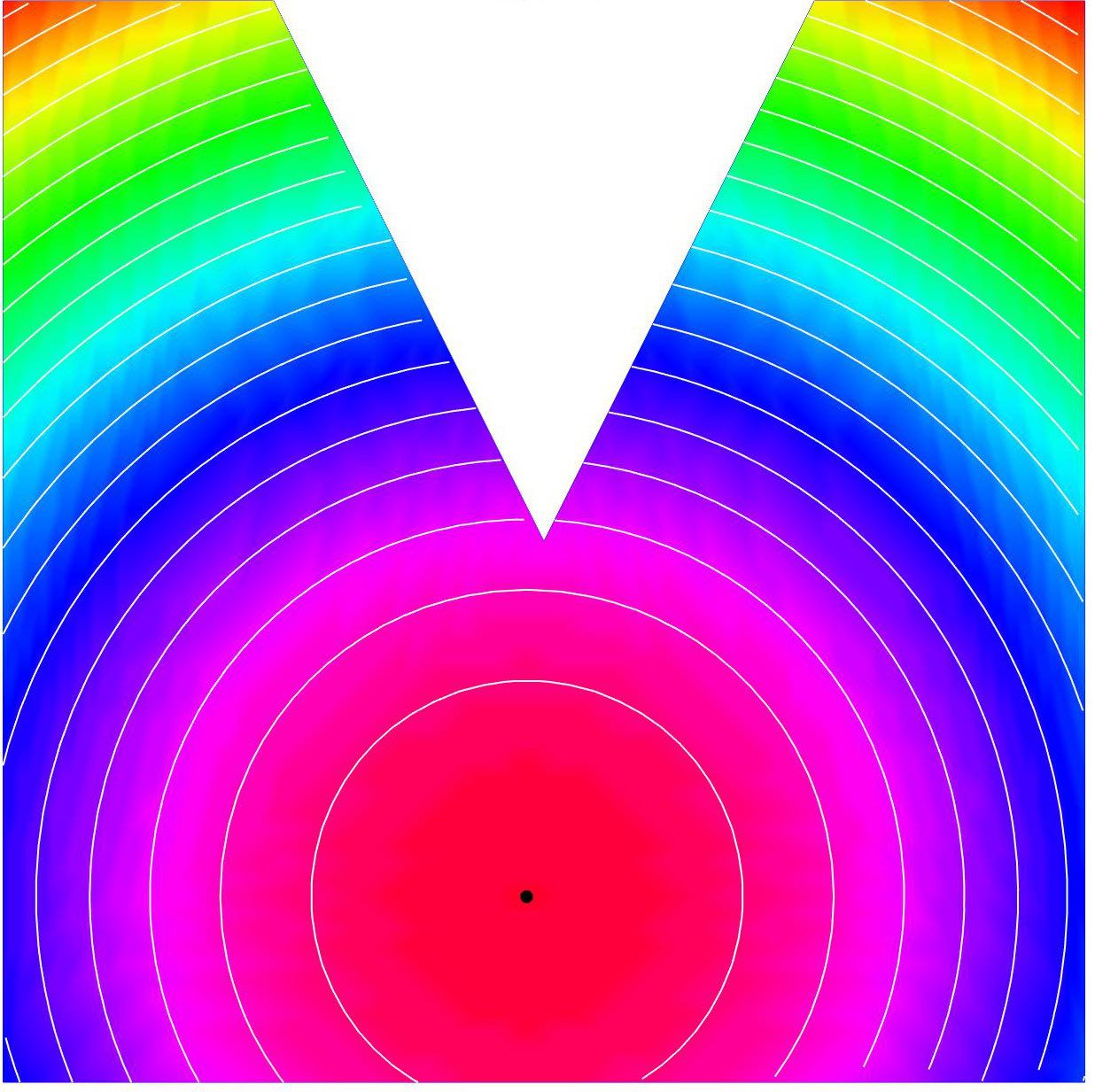}
&\includegraphics[height=42pt]{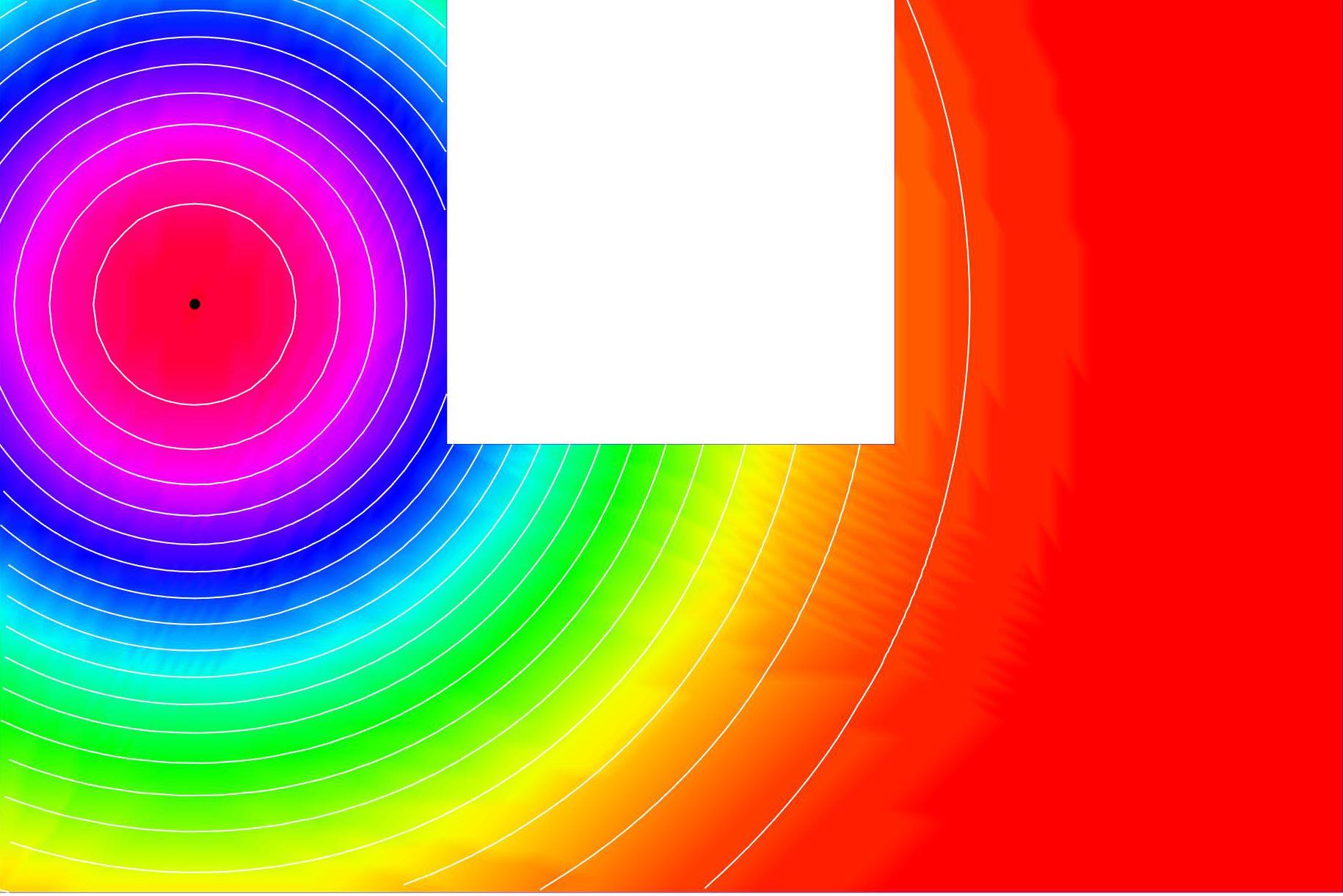}
&\includegraphics[height=42pt]{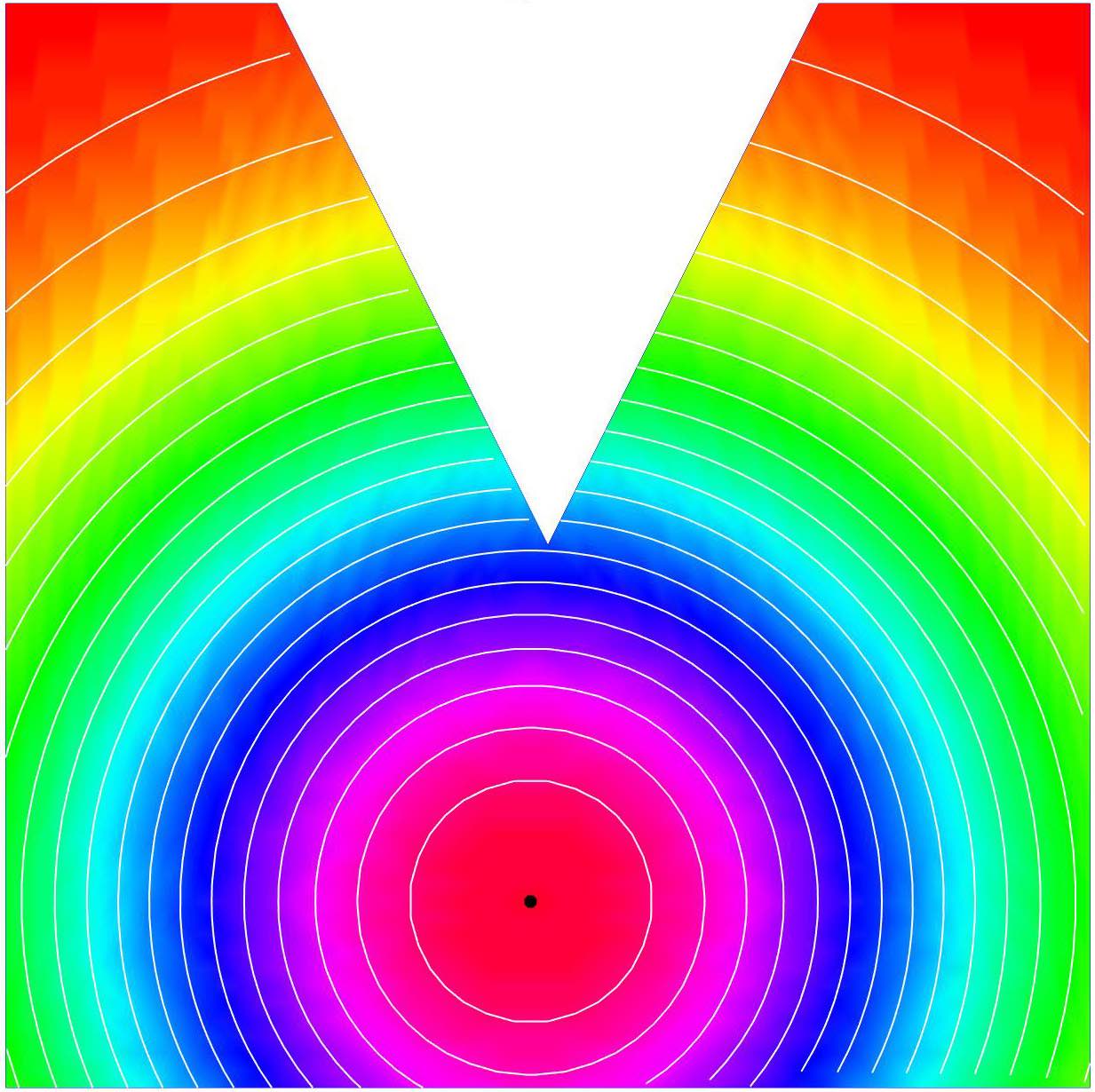}\\
\multicolumn{2}{c|}{$t=1$}	&\multicolumn{2}{c}{$t=10^{-1}$}\\
\includegraphics[height=42pt]{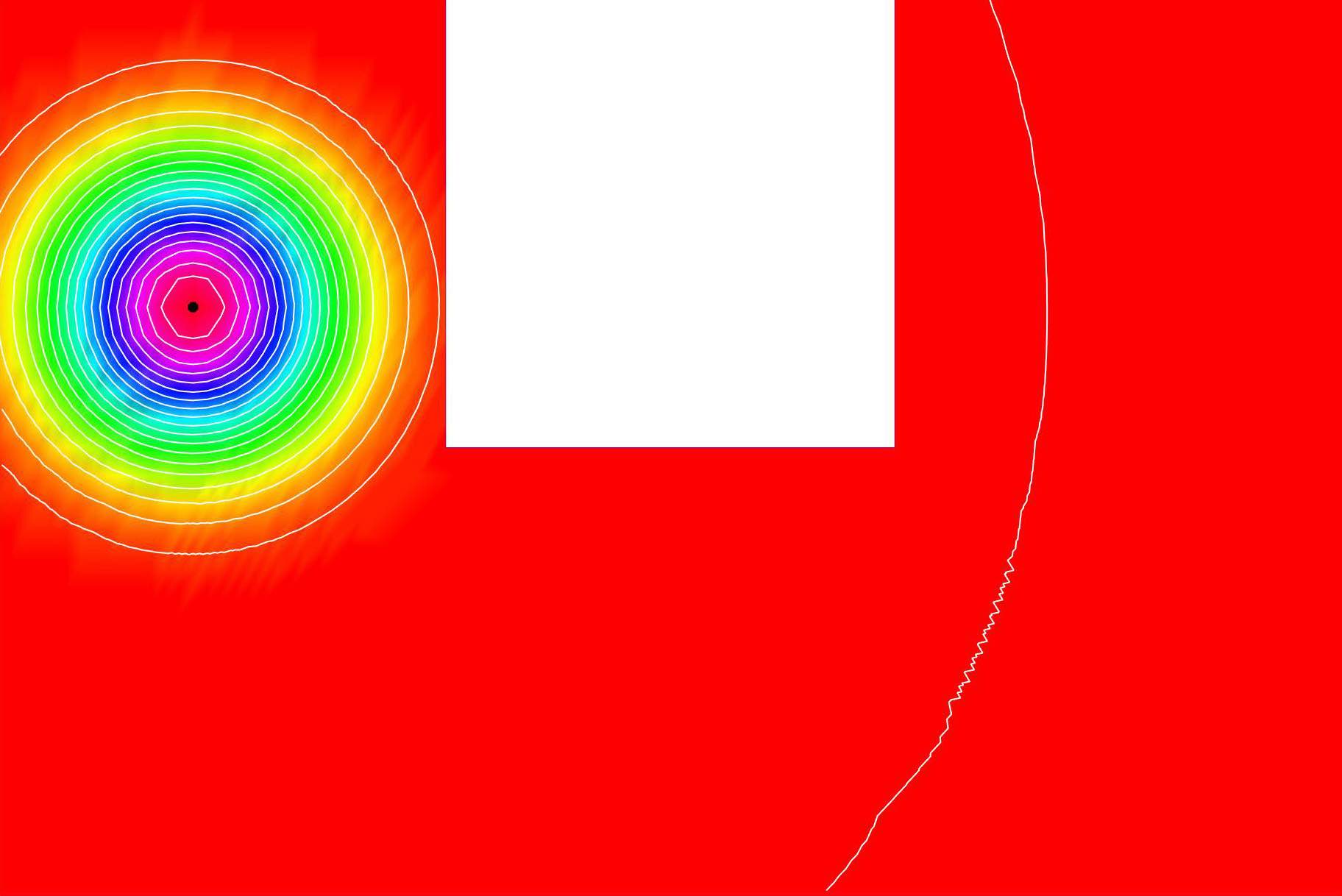}
&\includegraphics[height=42pt]{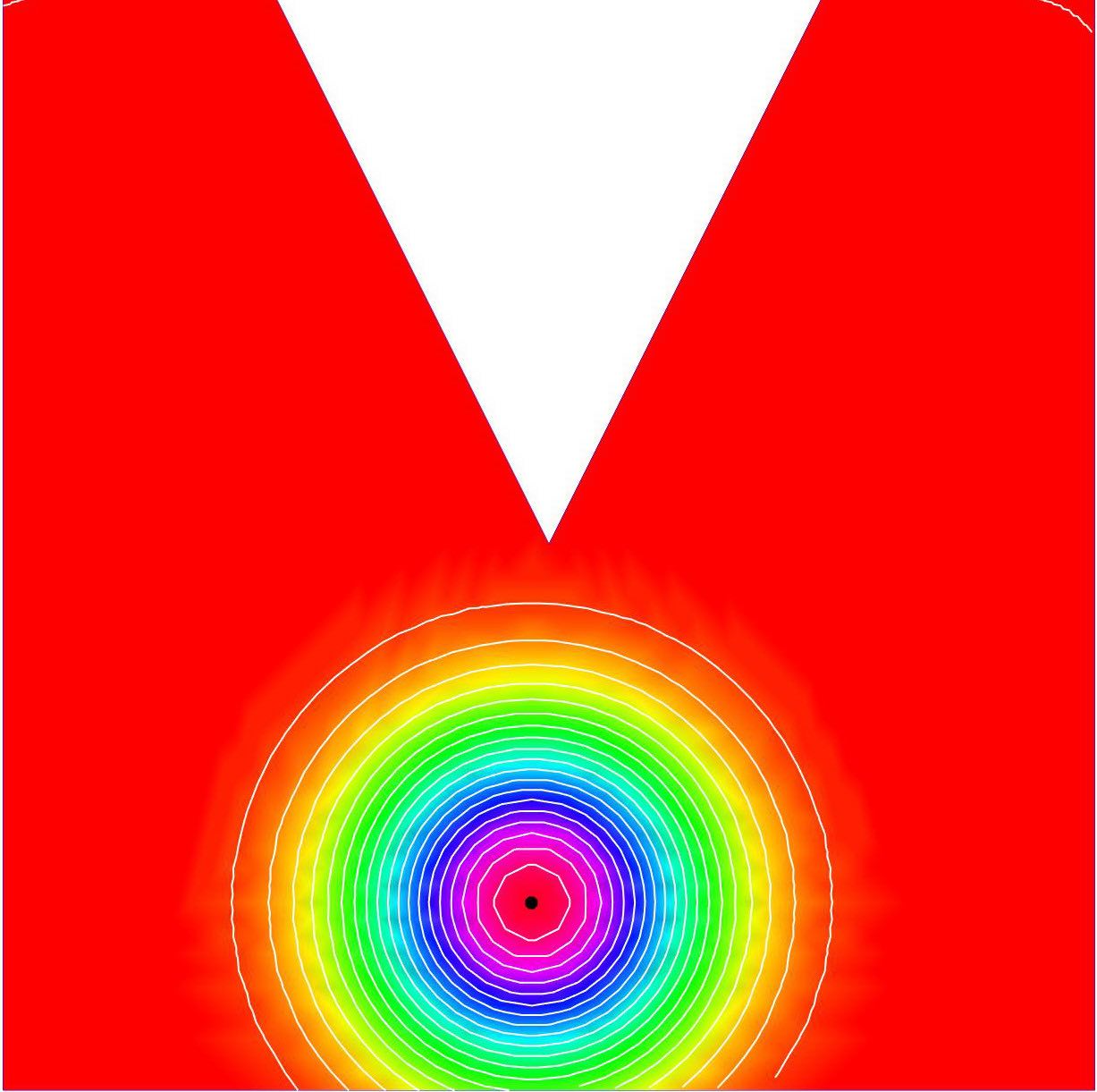}
&\includegraphics[height=42pt]{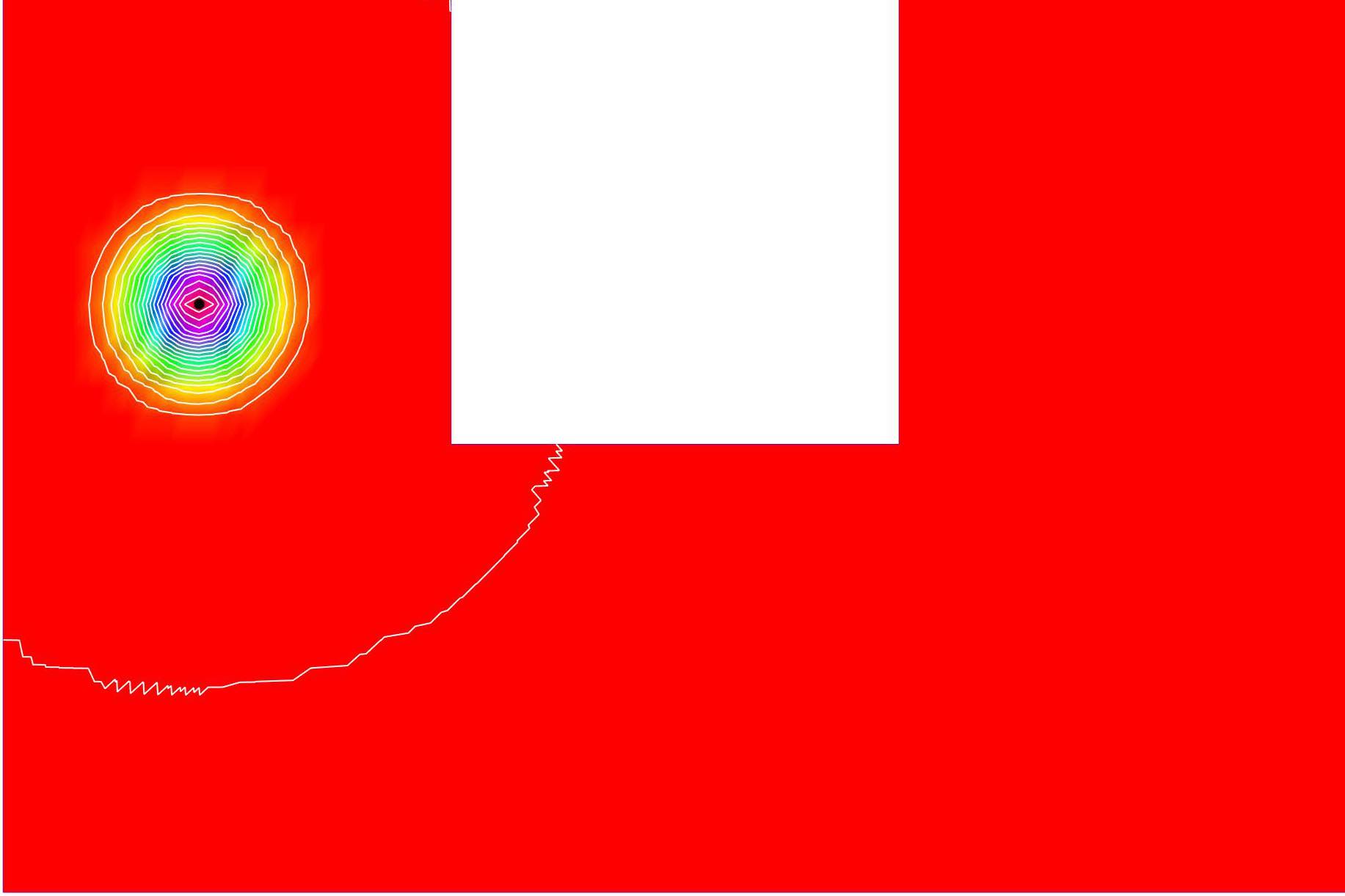}
&\includegraphics[height=42pt]{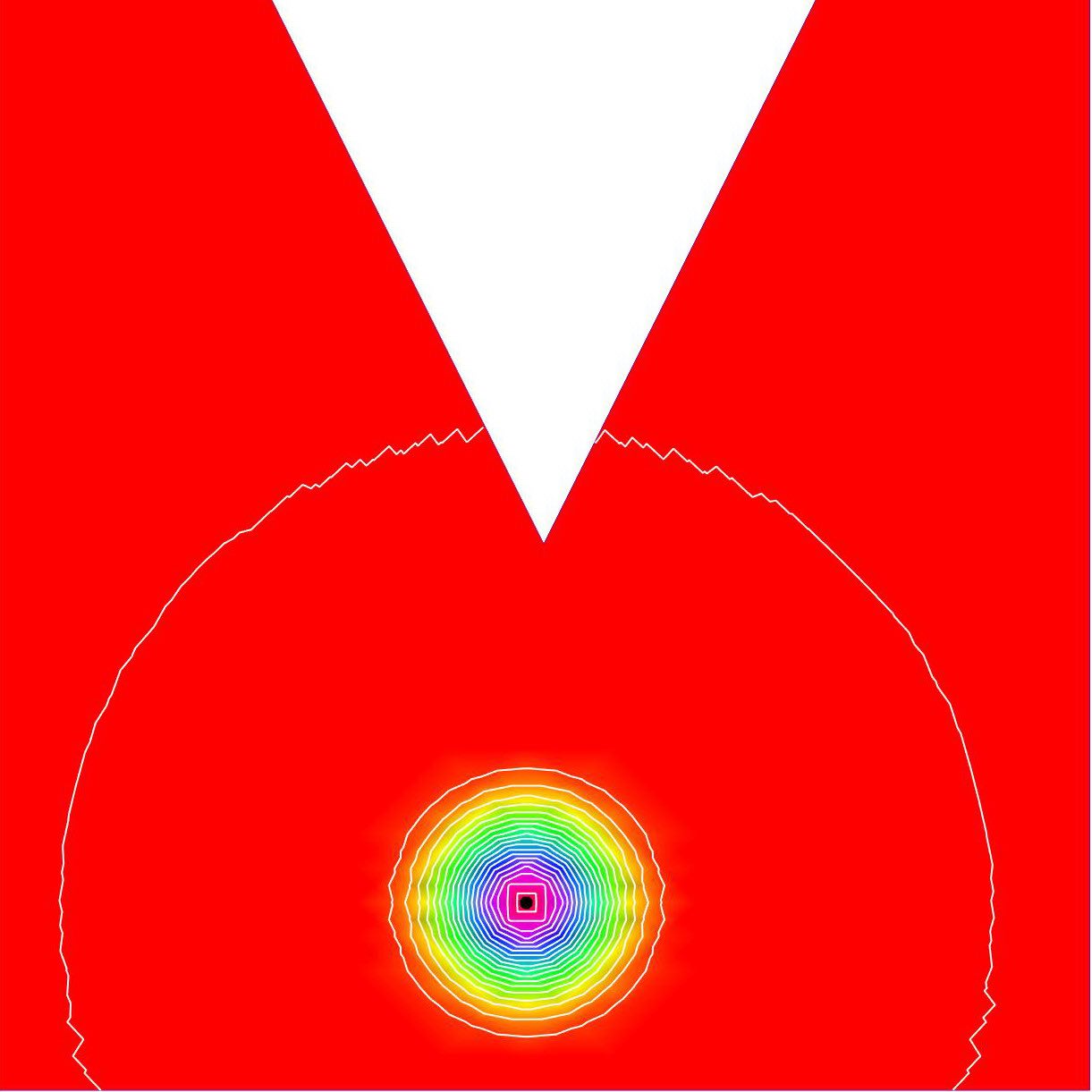}
\end{tabular}
\caption{Level-sets and color-map of multi-scale membership functions at a seed point, induced by the Gaussian kernel.\label{fig:2D-ANALYTIC-GAUSSIAN}}
\end{figure}
\subsection{Continuous F-transform as integral operator\label{sec:CONTINUOUS-FT-INT-OPER}}
We consider a \emph{membership function} \mbox{$A:\Omega\times\Omega\rightarrow\mathbb{R}$},~$\Omega$ compact, which is \emph{positive} and \emph{symmetric} (i.e., \mbox{$A(\mathbf{p},\mathbf{q})=A(\mathbf{q},\mathbf{p})$}, \mbox{$\forall\mathbf{p},\mathbf{q}\in\Omega$}), and introduce its maximum and minimum values
\begin{equation}\label{eq:UPPER-BOUND-MF}
\left\{
\begin{array}{l}
A_{\max}:=\max_{\mathbf{p},\mathbf{q}}\{\vert A(\mathbf{p},\mathbf{q})\vert\},\\
A_{\min}:=\min_{\mathbf{p},\mathbf{q}}\{\vert A(\mathbf{p},\mathbf{q})\vert\},
\end{array}
\right.\quad
0<A_{\min}\leq A_{\max}.
\end{equation}
If the input membership function \mbox{$A(\cdot,\cdot)$} is continuous in \mbox{$\Omega\times\Omega$}, then~$A_{\min}$ and~$A_{\max}$ are well-defined; this assumption is typically satisfied in real applications. We define the \emph{continuous F-transform} of \mbox{$f:\Omega\rightarrow\mathbb{R}$} as the function \mbox{$F:\Omega\rightarrow\mathbb{R}$}
\begin{equation}\label{eq:CONTINUOUS-FT}
\begin{split}
F(\mathbf{p}):&=\frac{\int_{\Omega} A(\mathbf{p},\mathbf{q})f(\mathbf{q})d\mathbf{q}}{S(\mathbf{p})},
\quad S(\mathbf{p}):=\int_{\Omega} A(\mathbf{p},\mathbf{q})d\mathbf{q},\\
&=\int_{\Omega}K(\mathbf{p}, \mathbf{q})f(\mathbf{q})d \mathbf{q},\qquad\mathbf{p}\in\Omega,
\end{split}
\end{equation}
where the \emph{integral kernel} is defined as
\begin{equation}\label{eq:NORMALISED-KERNEL}
K:\Omega\times\Omega\rightarrow\mathbb{R},\qquad
K(\mathbf{p},\mathbf{q})
:=\frac{A(\mathbf{p},\mathbf{q})}{S(\mathbf{p})}.
\end{equation}
Indeed, \mbox{$K(\mathbf{p},\mathbf{q})$} is equal to the value of the membership function \mbox{$A(\mathbf{p},\mathbf{q})$} normalised by the average value \mbox{$S(\mathbf{p})$} of the membership functions at~$\mathbf{p}$. The last equality in Eq. (\ref{eq:CONTINUOUS-FT}) provides the link between the F-transform and the \emph{integral operator} \mbox{$\mathcal{L}_{K}:\mathcal{L}^{2}(\Omega)\rightarrow\mathcal{C}^{0}(\Omega)$}
\begin{equation}\label{eq:INTEGRAL-OPERATOR}
(\mathcal{L}_{K}f)(\mathbf{p}):=\int_{\Omega}K(\mathbf{p},\mathbf{q})f(\mathbf{q})d\mathbf{q}
=\langle K(\mathbf{p},\cdot),f\rangle_{2}.
\end{equation}
Alternatively, we consider the \emph{symmetric kernel}
\begin{equation}\label{eq:NORMALISED-KERNEL-SYMMETRIC}
K:\Omega\times\Omega\rightarrow\mathbb{R},\qquad
K(\mathbf{p},\mathbf{q})
:=\frac{A(\mathbf{p},\mathbf{q})}{(S(\mathbf{p})S(\mathbf{q}))^{1/2}},
\end{equation}
which is equal to the value of the membership function \mbox{$A(\mathbf{p},\mathbf{q})$} normalised by the product \mbox{$(S(\mathbf{p})S(\mathbf{q}))^{1/2}$} of the average values \mbox{$S(\mathbf{p})$}, \mbox{$S(\mathbf{q})$}. Then, the continuous F-transform is still defined as in Eq. (\ref{eq:NORMALISED-KERNEL}). Selecting a set \mbox{$\mathcal{P}:=\{\mathbf{p}_{i}\}_{i=1}^{n}$}, the membership functions introduced in Sect.~\ref{sec:PREVIOUS-WORK} are defined by centring the function \mbox{$A(\cdot,\cdot)$} at the set~$\mathcal{P}$, i.e., \mbox{$A_{i}(\cdot):=A(\mathbf{p}_{i},\cdot)$}. For details on the relation between kernels and membership functions, we refer the reader to Sect.~\ref{KERNEL-EXAMPLE}.

Combining (\ref{eq:CONTINUOUS-FT}) with (\ref{eq:INTEGRAL-OPERATOR}), the continuous F-transform \mbox{$F\equiv\mathcal{L}_{K}f$} of~$f$ is the function defined as the action of the integral operator~$\mathcal{L}_{K}$, induced by the normalised kernel \mbox{$K(\cdot,\cdot)$} in Eq. (\ref{eq:NORMALISED-KERNEL}), on~$f$. Noting that \mbox{$F(\mathbf{p}_{i})$} is the~$i$-th component~$F_{i}$ of the F-transform of~$f$, the continuous F-transform \mbox{$F(\cdot)$} interpolates the values \mbox{$\mathbf{F}_{n}:=(F_{i})_{i=1}^{n}$} of the F-transform associated with the set \mbox{$(f(\mathbf{q}_{i}))_{i=1}^{s}$} of~$f$-values at \mbox{$\mathcal{P}:=\{\mathbf{p}_{i}\}_{i=1}^{n}$}, i.e.,
\begin{equation}\label{eq:DISCRETE-CONTINUOUS-FT}
F_{i}
=F(\mathbf{p}_{i})
=(\mathcal{L}_{K}f)(\mathbf{p}_{i}),\qquad
i=1,\ldots,n.
\end{equation}
Since the function \mbox{$F(\cdot)$} in Eq. (\ref{eq:CONTINUOUS-FT}) is continuous, it can be evaluated at any point and Eq. (\ref{eq:DISCRETE-CONTINUOUS-FT}) is well-defined (Sect.~\ref{sec:PROPERTIES-LKF}). 
\begin{figure}
\centering
\begin{tabular}{cc|cc}
\multicolumn{2}{c|}{Multi-Quad}
&\multicolumn{2}{c}{Inverse Multi-Quad}\\
\includegraphics[height=41pt]{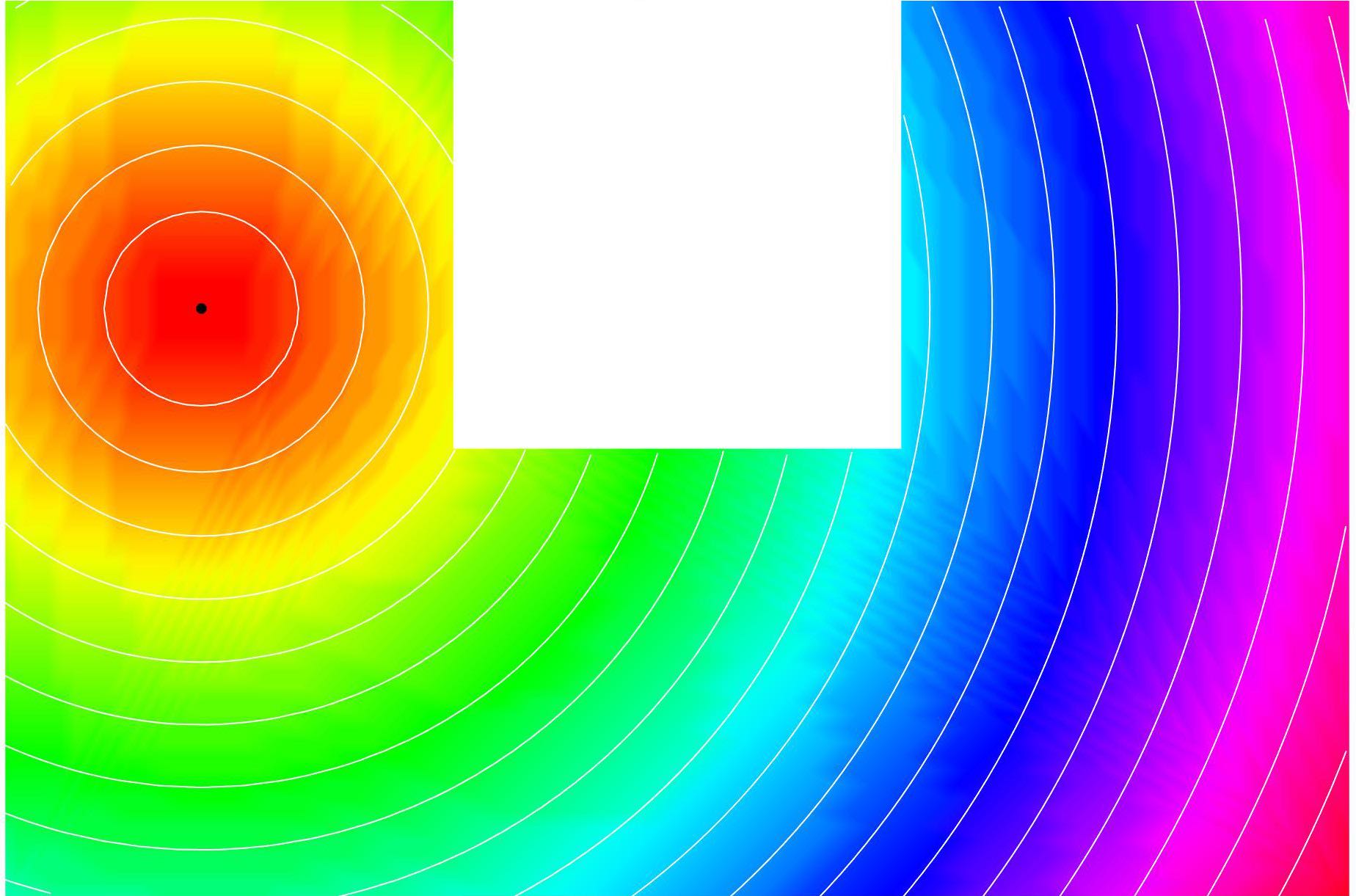}
&\includegraphics[height=41pt]{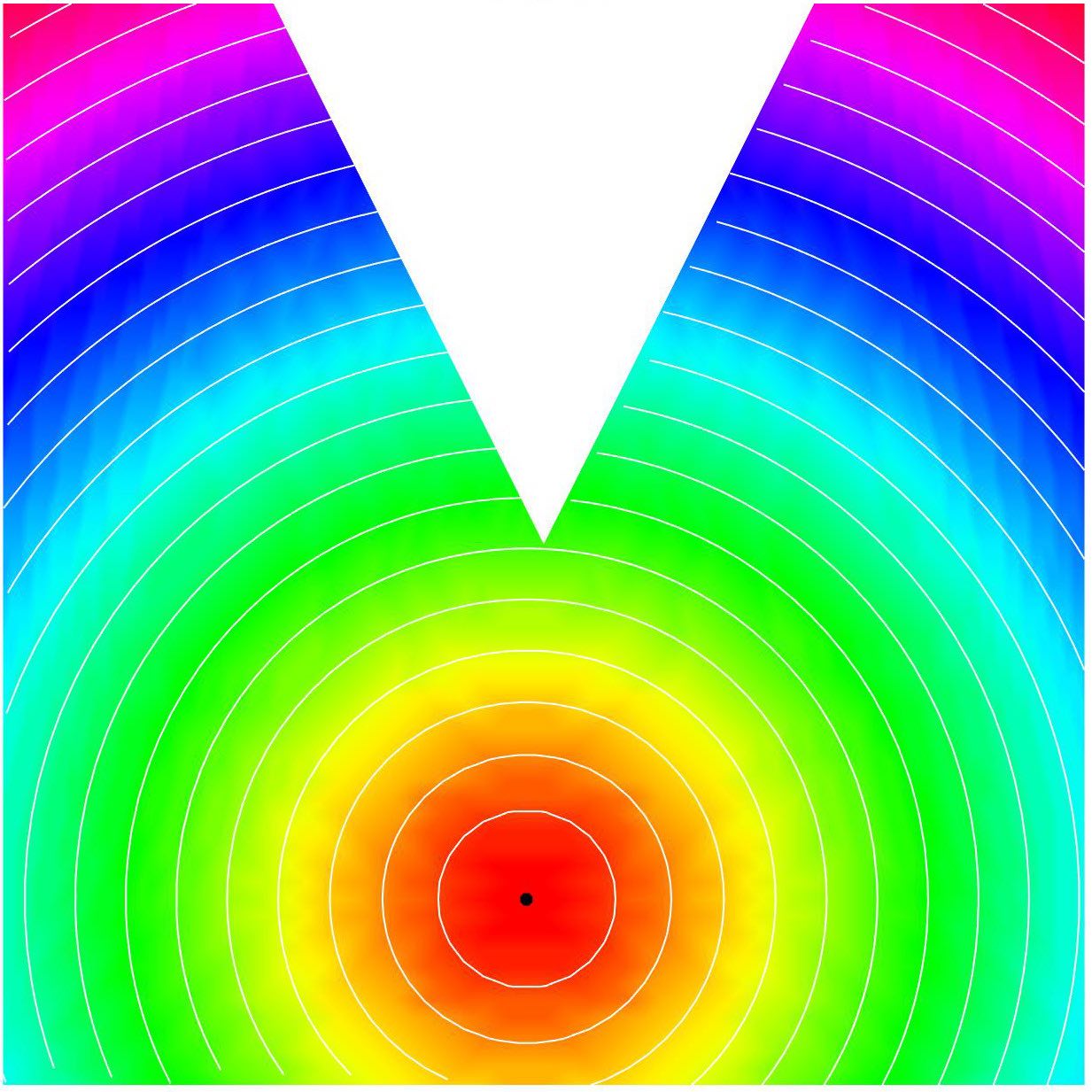}
&\includegraphics[height=41pt]{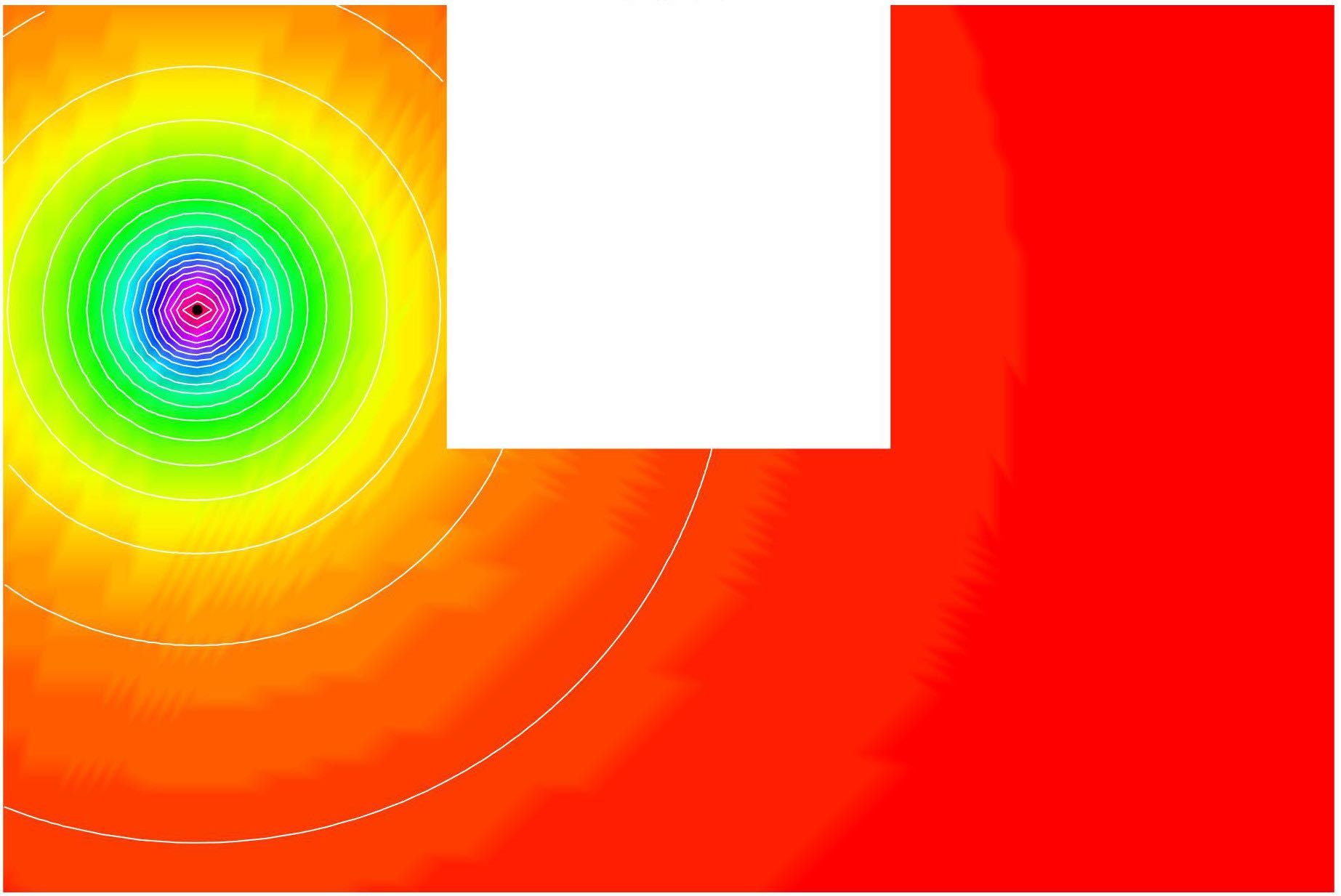}
&\includegraphics[height=41pt]{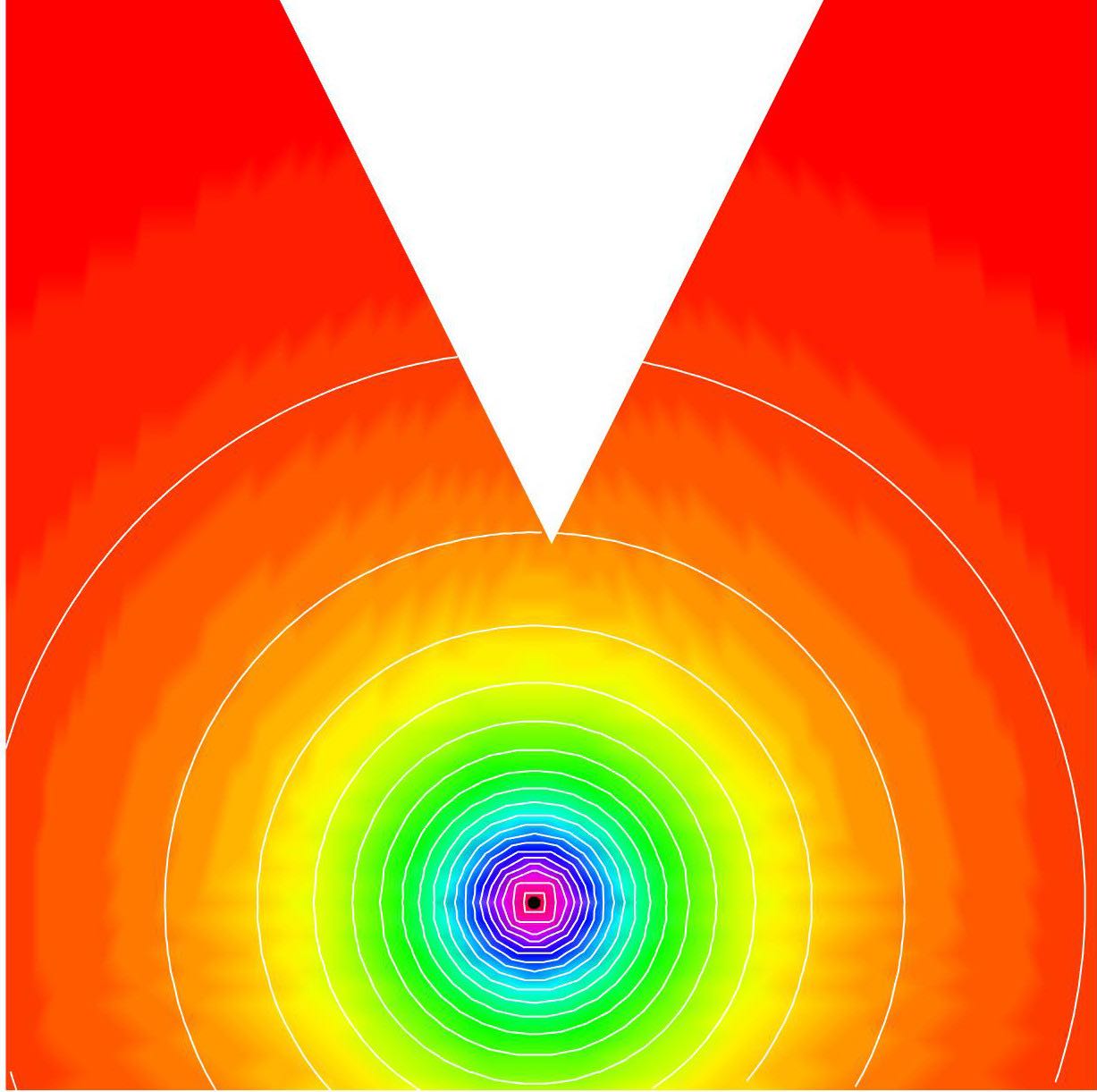}
\end{tabular}
\caption{Level-sets of membership functions at a seed point, induced by multi-quadratic and inverse multi-quadratic kernels.\label{fig:2D-DIFFUSION-MF}}
\end{figure}
\subsection{Properties of the continuous F-transform\label{sec:CONTINUOUS-FT-PROPERTIES}}
Assuming that the kernel is continuous, we define its maximum value as \mbox{$C_{K}:=\max_{\mathbf{p},\mathbf{q}\in\Omega}\{\vert K(\mathbf{p},\mathbf{q})\vert\}$}; in particular, \mbox{$K\in\mathcal{L}^{2}(\Omega\times\Omega)$}. Under this assumption, we discuss the properties of the continuous F-transform~$\mathcal{L}_{K}$, such as well-posedness, linearity, self-adjointness and generalisation. 

\textbf{Well-posedness}
From Eqs. (\ref{eq:UPPER-BOUND-MF}), (\ref{eq:CONTINUOUS-FT}), we get that
\begin{equation}\label{eq:UPPER-BOUND}
S(\mathbf{p})\geq A_{\min}\vert\Omega\vert,\quad
\|A(\mathbf{p},\cdot)\|_{2}
\leq A_{\max}\vert\Omega\vert^{1/2},
\end{equation}
with \mbox{$\vert\Omega\vert$} measure (e.g., area, volume) of~$\Omega$. Then,
\begin{equation*}
\begin{split}
\vert F(\mathbf{p})\vert
:&=\vert\mathcal{L}_{K}f(\mathbf{p})\vert
=\frac{\left\vert \int_{\Omega}A(\mathbf{p},\mathbf{q})f(\mathbf{q})d\mathbf{q}\right\vert}{S(\mathbf{p})}\\
&\leq \frac{\|A(\mathbf{p},\cdot)\|_{2}\|f\|_{2}}{S(\mathbf{p})}
\leq_{\textrm{Eq.}(\ref{eq:UPPER-BOUND})}\frac{A_{\max}}{A_{\min}}\,\frac{\|f\|_{2}}{\vert\Omega\vert^{1/2}},
\quad\mathbf{p}\in\Omega,
\end{split}
\end{equation*}
i.e., \mbox{$\mathcal{L}_{K}f$} is well-defined on~$\Omega$. Choosing the kernel in Eq. (\ref{eq:NORMALISED-KERNEL-SYMMETRIC}), we get an analogous upper bound; in fact,
\begin{equation*}
\begin{split}
\vert F(\mathbf{p})\vert
:&=\vert\mathcal{L}_{K}f(\mathbf{p})\vert
\leq\int_{\Omega}\frac{\left\vert A(\mathbf{p},\mathbf{q})f(\mathbf{q})\right\vert}{(S(\mathbf{p})S(\mathbf{q}))^{1/2}}d\mathbf{q}\\
&\leq \frac{\|A(\mathbf{p},\cdot)\|_{2}\|f\|_{2}}{A_{\min}\vert\Omega\vert}\leq_{\textrm{Eq.}(\ref{eq:UPPER-BOUND})}\frac{A_{\max}}{A_{\min}}\,\frac{\|f\|_{2}}{\vert\Omega\vert^{1/2}},
\quad\mathbf{p}\in\Omega.
\end{split}
\end{equation*}
\textbf{Linearity and self-adjointness of~$\mathcal{L}_{K}$}
The operator~$\mathcal{L}_{K}$ is linear (as~$\mathcal{L}^{2}(\Omega)$-scalar product in Eq. (\ref{eq:INTEGRAL-OPERATOR})) and self-adjoint as a consequence of the symmetry of the kernel (e.g., for symmetric membership functions); in fact, for any \mbox{$f,g,\in\mathcal{L}^{2}(\Omega)$},
\begin{equation*}\label{eq:SELF-ADJOINTNESS}
\langle\mathcal{L}_{K}f,g\rangle_{2}
=\int_{\Omega\times\Omega} K(\mathbf{p},\mathbf{q})f(\mathbf{p})g(\mathbf{q})d\mathbf{p}\,d\mathbf{q}
=\langle f,\mathcal{L}_{K}g\rangle_{2}.
\end{equation*}
\textbf{Positive definiteness of~$\mathcal{L}_{K}$}
If \mbox{$K(\cdot,\cdot)$} is a positive-definite kernel (i.e., \mbox{$\mathbf{K}:=(K(\mathbf{p}_{i},\mathbf{p}_{j}))_{i,j}$} is a positive-definite matrix for any finite set of points in~$\Omega$), then~$\mathcal{L}_{K}$ is a positive-definite operator (i.e., \mbox{$\langle\mathcal{L}_{K}f,f\rangle_{2}\geq 0$}, \mbox{$\forall f$}).

\textbf{Continuity of~$\mathcal{L}_{K}$}
Endowing \mbox{$\mathcal{C}^{0}(\Omega)$} with the~$\mathcal{L}^{2}$-norm,~$\mathcal{L}_{K}$ is a bounded (i.e., continuous) operator; in fact,
\begin{equation*}
\begin{split}
\|\mathcal{L}_{K}f\|_{2}
\leq\left[\int_{\Omega}\| K(\mathbf{p},\cdot)\|_{2}^{2}\,d\mathbf{p}\right]^{1/2}\|f\|_{2}
=\|K\|_{2}\|f\|_{2}.
\end{split}
\end{equation*}

\textbf{Generalisation property\label{sec:ACCURACY-CFT}}
We evaluate the approximation error between \mbox{$\mathcal{L}_{K}f(\mathbf{p})$} and~$F_{i}$ (c.f., Eq. (\ref{eq:DISCRETE-CONTINUOUS-FT})) in a neighbour~$\mathcal{N}_{\mathbf{p}_{i}}$ of~$\mathbf{p}_{i}$, i.e., the discrepancy between the continuous F-transform at \mbox{$\mathbf{p}\in\mathcal{N}_{\mathbf{p}_{i}}$} and the F-transform~$F_{i}$ at~$\mathbf{p}_{i}$. Noting that
\begin{equation*}
\begin{split}
\left\vert
\mathcal{L}_{K}f(\mathbf{p})-\mathcal{L}_{K}f(\mathbf{p}_{i})
\right\vert
&=\left\vert
\int_{\Omega}\left(K(\mathbf{p},\mathbf{q})-K(\mathbf{p}_{i},\mathbf{q})\right)f(\mathbf{q})\textrm{d}\mathbf{q}\right\vert\\
&\leq\|K(\mathbf{p},\cdot)-K(\mathbf{p}_{i},\cdot)\|_{2}\|f\|_{2}\\
&\leq\vert\Omega\vert\,\|K(\mathbf{p},\cdot)-K(\mathbf{p}_{i},\cdot)\|_{\infty}\|f\|_{2},
\end{split}
\end{equation*}
the error is guided by the difference of the membership functions at~$\mathbf{p}$,~$\mathbf{p}_{i}$. Selecting a kernel localised around its center~$\mathbf{p}_{i}$ (e.g., compactly-supported or Gaussian kernels) generally improves the convergence of \mbox{$K(\mathbf{p},\mathbf{q})$} to \mbox{$K(\mathbf{p}_{i},\mathbf{q})$}.
\begin{figure}[t]
\centering
\begin{tabular}{cc|cc}
\multicolumn{2}{c|}{Harmonic}
&\multicolumn{2}{c}{Bi-harmonic}\\
\includegraphics[height=41pt]{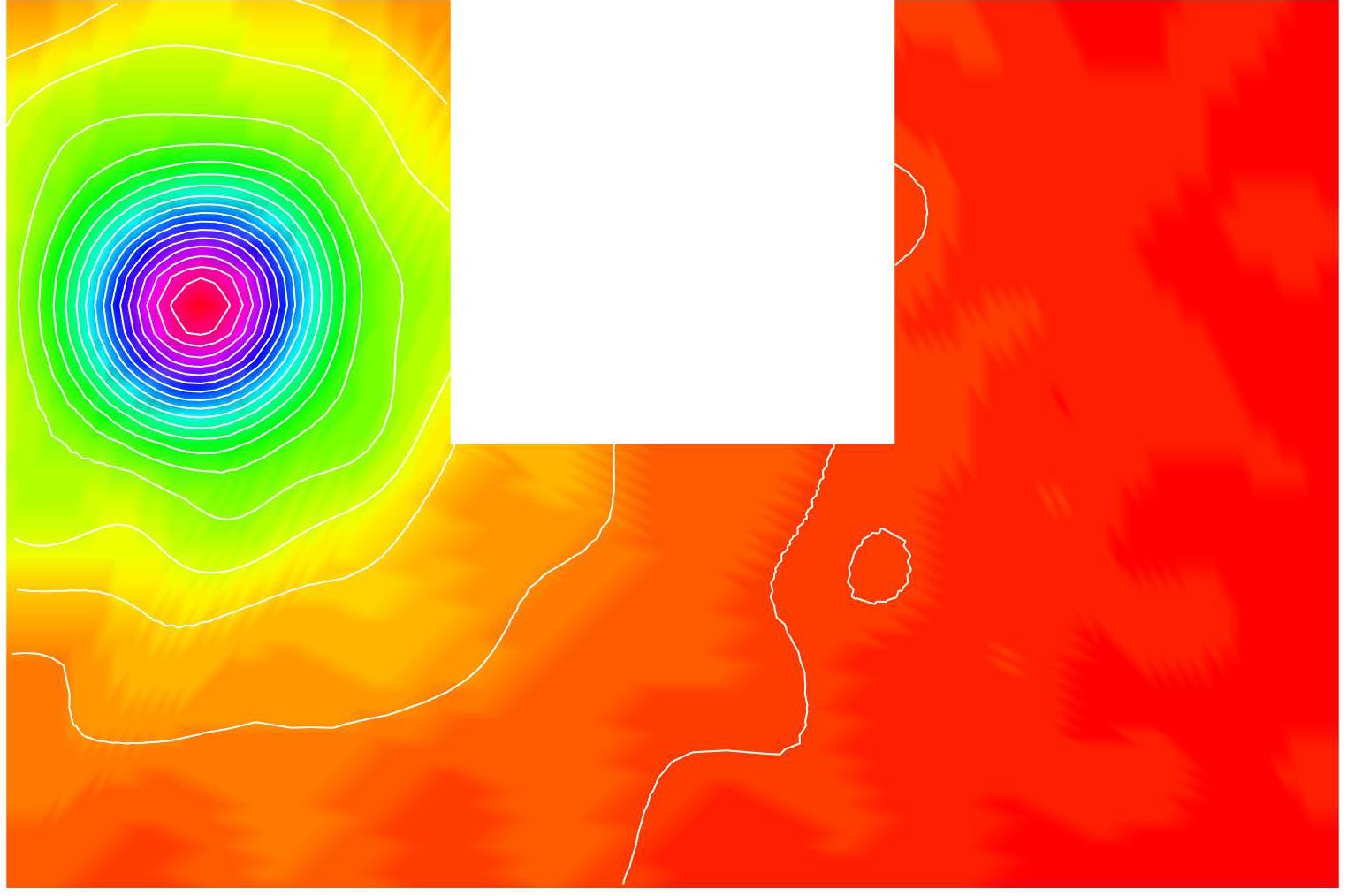}
&\includegraphics[height=41pt]{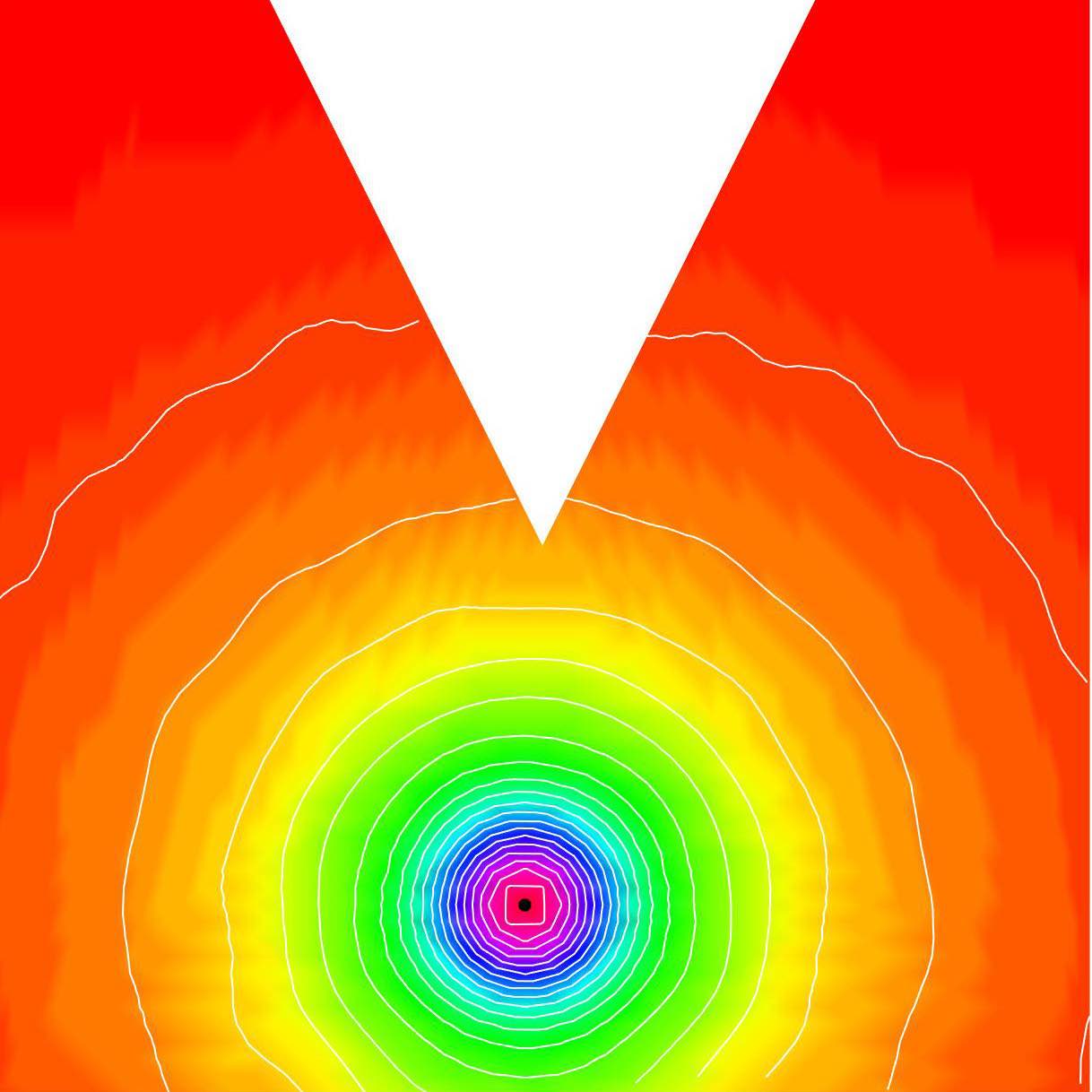}
&\includegraphics[height=41pt]{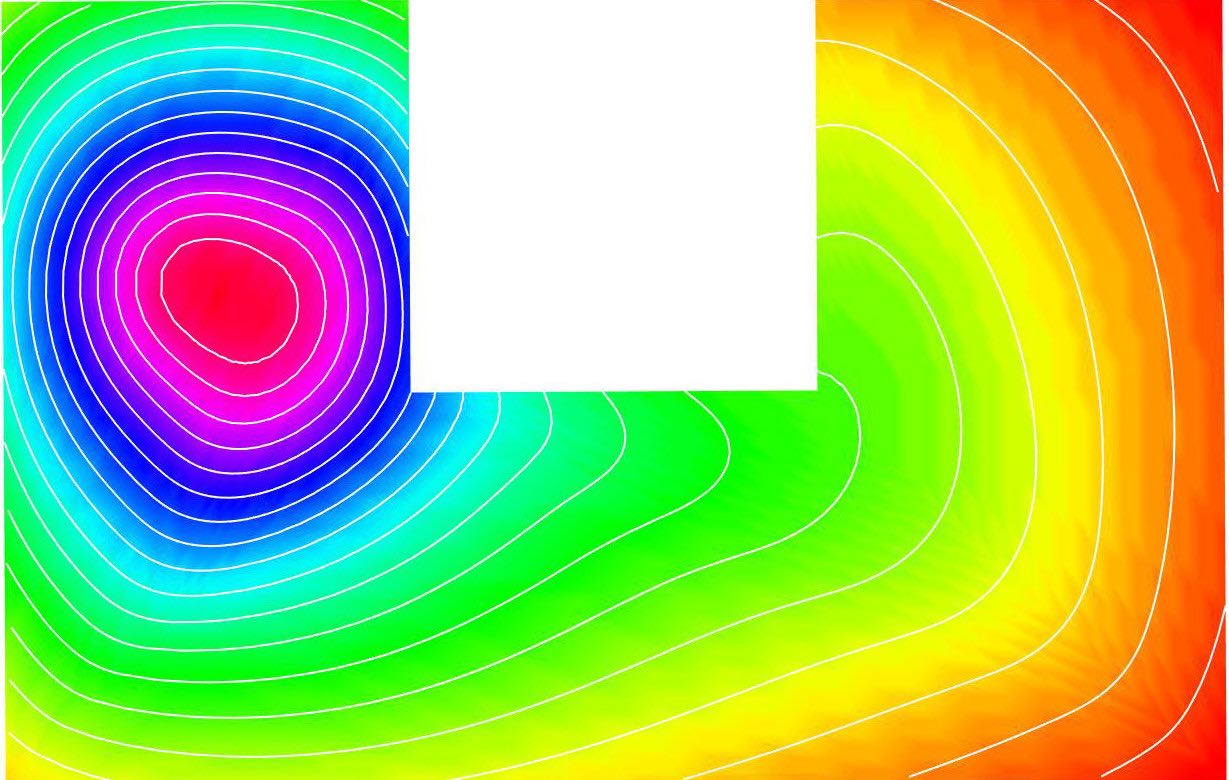}
&\includegraphics[height=40pt]{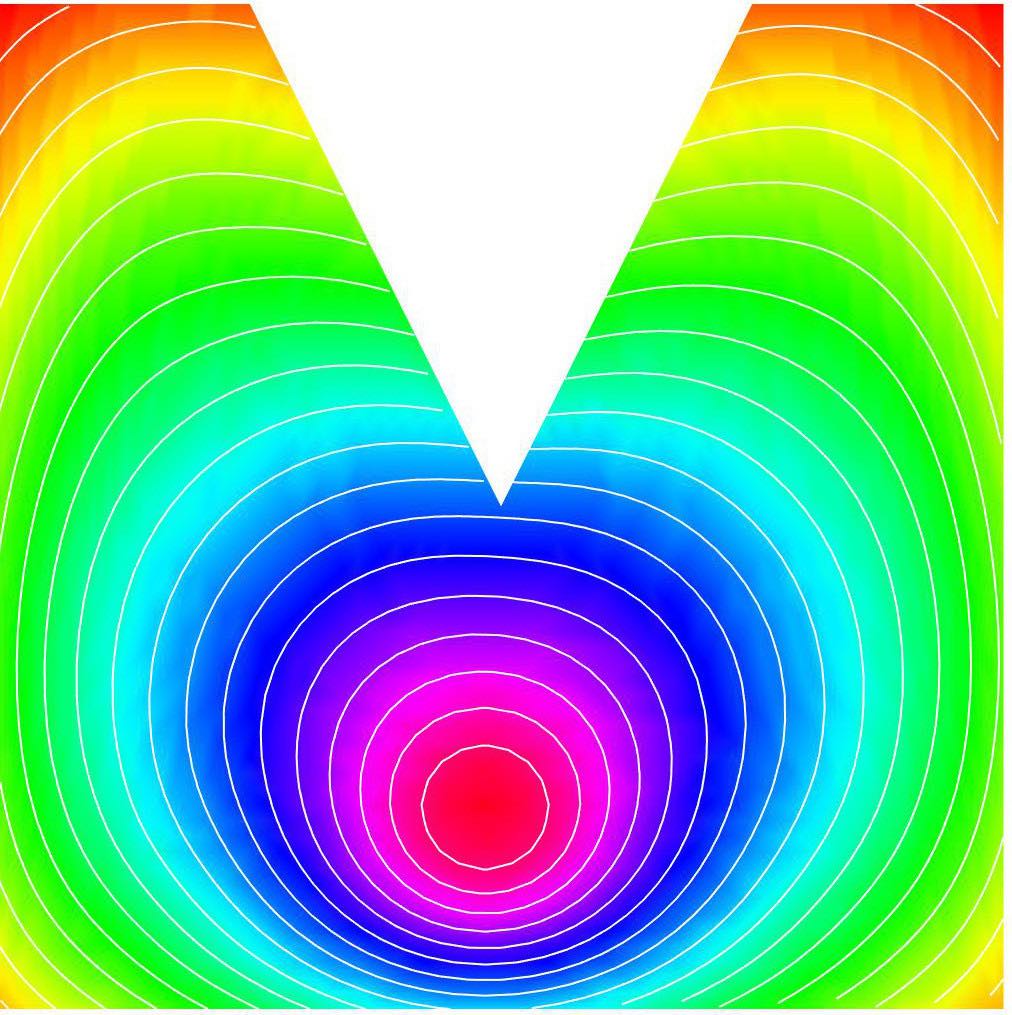}
\end{tabular}
\caption{Level-sets of data-driven harmonic and biharmonic membership functions at a seed point.\label{fig:2D-HAR-BIHAR}}
\end{figure}
\subsection{Properties of the continuous F-transform~$\mathcal{L}_{K}f$\label{sec:PROPERTIES-LKF}}
We now focus on the properties of the functions \mbox{$\mathcal{L}_{K}f$} in the image of the continuous F-transform.

\textbf{Least-squares property}
Analogously to the discrete F-transform, the continuous F-transform~$F$ minimises the \emph{energy function} \mbox{$E(\mathbf{q},t):=\int_{\Omega}K(\mathbf{p},\mathbf{q}) \vert f(\mathbf{p})-t\vert^{2}d\mathbf{p}$}, \mbox{$t\in\mathbb{R}$}. In fact, \mbox{$\partial_{t}E(\mathbf{q},t)=0$} if and only if \mbox{$t=F(\mathbf{p})$}.

\textbf{Continuity of~$\mathcal{L}_{K}f$}
The continuous F-transform of~$f$ is a continuous function; in fact,
\begin{equation*}
\begin{split}
\vert F(\mathbf{p})-F(\mathbf{q})\vert
&=\left\vert\langle K(\mathbf{p},\cdot)-K(\mathbf{q},\cdot),f\rangle_{2}\right\vert\\
&\leq \|K(\mathbf{p},\cdot)-K(\mathbf{q},\cdot)\|_{2}\,\|f\|_{2}\\
&\leq\vert\Omega\vert\|K(\mathbf{p},\cdot)-K(\mathbf{q},\cdot)\|_{\infty}\,\|f\|_{2}.
\end{split}
\end{equation*}
Since the kernel is continuous and~$\Omega$ is compact, it is uniformly continuous (i.e., \mbox{$\lim_{\mathbf{p}\rightarrow\mathbf{q}}\|K(\mathbf{p},\cdot)-K(\mathbf{q},\cdot)\|_{\infty}=0$}); indeed, \mbox{$F=\mathcal{L}_{K}f$} is continuous and is evaluated at any point of~$\Omega$. In particular, \mbox{$F_{i}=\mathcal{L}_{K}f(\mathbf{p}_{i})$} is well-defined.

\textbf{Boundness of~$\mathcal{L}_{K}f$}
The continuous F-transform is bounded, as continuous function on a compact set. The upper bound 
\begin{equation*}
\vert\mathcal{L}_{K}f(\mathbf{p})\vert
\leq\|K(\mathbf{p},\cdot)\|_{2}\|f\|_{2}\\
\leq C_{K}\vert\Omega\vert\,\|f\|_{2},
\end{equation*}
allows us to estimate the maximum variation of the values of the continuous F-transform in terms of the constant~$C_{K}$, the area or volume of \mbox{$\vert\Omega\vert$}, and the~$\mathcal{L}^{2}$ norm of the input function.

\subsection{Relations between discrete and continuous F-transform\label{sec:HILBERT-FT}}
We further study the relation between the discrete and continuous F-transforms through the \emph{sampling operator} \mbox{$\mathcal{R}:\mathcal{C}^{0}(\Omega)\rightarrow\mathbb{R}^{s}$}, \mbox{$f\mapsto\mathcal{R}f:=(f(\mathbf{q}_{i}))_{i=1}^{s}$}, and the \emph{out-of-sample operator} \mbox{$\mathcal{E}:\mathbb{R}^{s}\rightarrow\mathcal{C}^{0}(\Omega)$}, \mbox{$\mathbf{f}=(f_{i})_{i=1}^{s}\mapsto\mathcal{E}\mathbf{f}$}, with \mbox{$(\mathcal{E}\mathbf{f})(\mathbf{q}_{i})=f_{i}$}, \mbox{$\forall i$}. We also require that the out-of-sample operator is linear: i.e., \mbox{$\mathcal{E}(\alpha\mathbf{f}+\beta \mathbf{g})=\alpha\mathcal{E}\mathbf{f}+\beta \mathcal{E}\mathbf{g}$}, \mbox{$\forall\alpha,\beta$}, \mbox{$\forall\mathbf{f},\mathbf{g}$}. To this end, we select a set \mbox{$(\phi_{i})_{i=1}^{n}$} of radial basis functions centred at \mbox{$(\mathbf{p}_{i})_{i=1}^{s}$} (or at any other set of points) and compute the function \mbox{$f(\mathbf{p})=\sum_{i=1}^{n}\alpha_{i}\phi_{i}(\mathbf{p})$} such that \mbox{$f(\mathbf{p}_{i})=f_{i}$}, \mbox{$i=1,\ldots,s$}. These conditions are equivalent to solve the linear system \mbox{$\mathbf{G}\alpha=\mathbf{f}$}, where~$\mathbf{G}$ is the Gram matrix associated with the input RBFs and~$\alpha$ is unknown vector. Then, we apply the out-of-sample operator to the set \mbox{$\mathbf{f}:=(f(\mathbf{q}_{i}))_{i=1}^{s}$} of the~$f$-values at~$\mathcal{Q}$ and consider the diagram
\begin{equation*}
\mathbf{f}\in\mathbb{R}^{s}\mapsto\mathcal{E}\mathbf{f}\in\mathcal{C}^{0}(\Omega)\mapsto(\mathcal{L}_{K}\mathcal{E}\mathbf{f})\in\mathcal{C}^{0}(\Omega).
\end{equation*}
From the upper bound
\begin{equation}\label{eq:UP-EXTENSION}
\|\mathcal{L}_{K}f-\mathcal{L}_{K}\mathcal{E}\mathbf{f}\|_{2}
=\|\mathcal{L}_{K}(f-\mathcal{E}\mathbf{f})\|_{2}\\
\leq\|K\|_{2}\|f-\mathcal{E}\mathbf{f}\|_{2},
\end{equation}
the error (\ref{eq:UP-EXTENSION}) between the continuous F-transform \mbox{$\mathcal{L}_{K}f$} and \mbox{$\mathcal{L}_{K}\mathcal{E}\mathbf{f}$} of~$f$ and \mbox{$\mathcal{E}\mathbf{f}$} is guided mainly by the accuracy of the approximation \mbox{$\mathcal{E}\mathbf{f}$} of~$f$. Through the restriction and extension operators, we introduce the diagram in Fig.~\ref{fig:DIAGRAM}, which summarises the relations between the continuous and discrete F-transforms and their inverse operators.
\begin{figure}[t]
\centering
\begin{tabular}{cc|cc}
\multicolumn{2}{c}{$t=5\times 10^{-3}$}
&\multicolumn{2}{c}{$t=10^{-2}$}\\
\includegraphics[height=41pt]{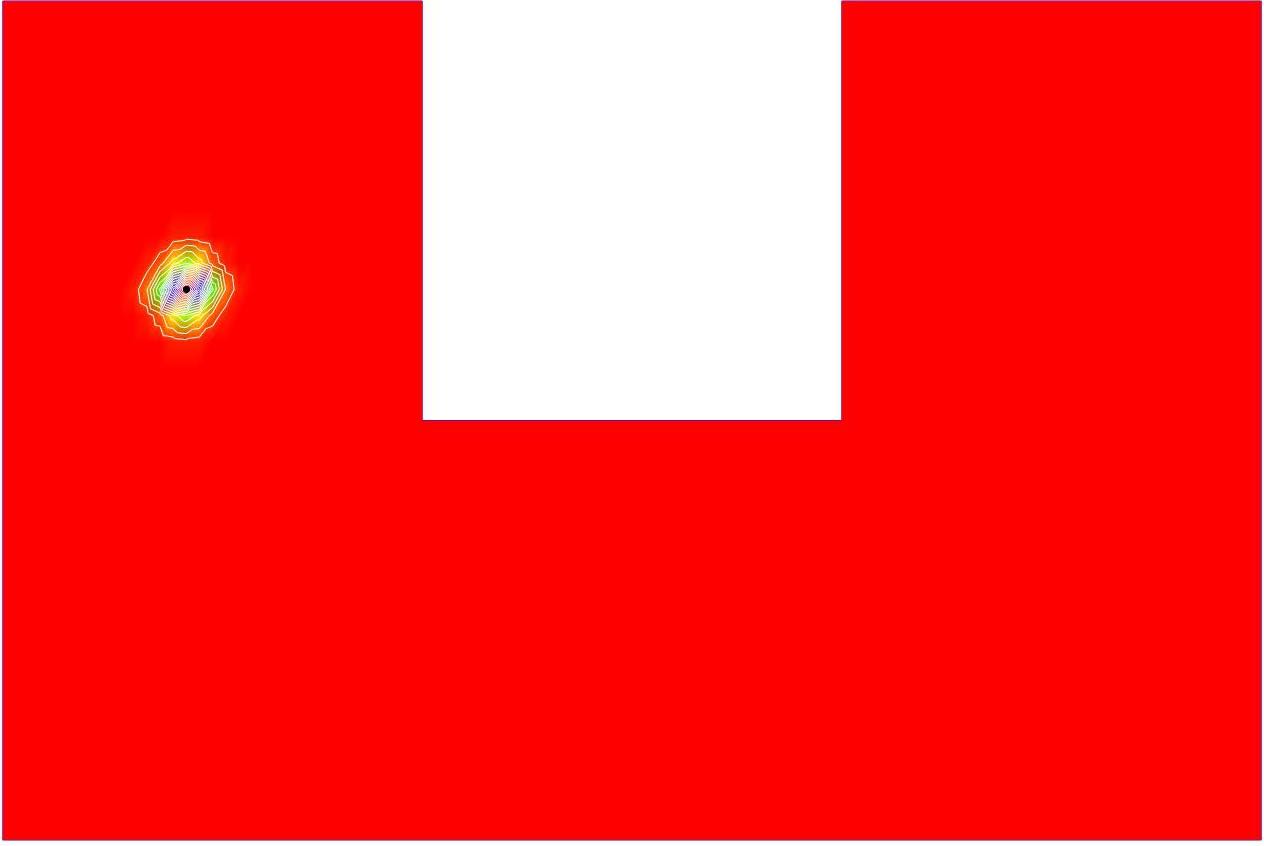}
&\includegraphics[height=41pt]{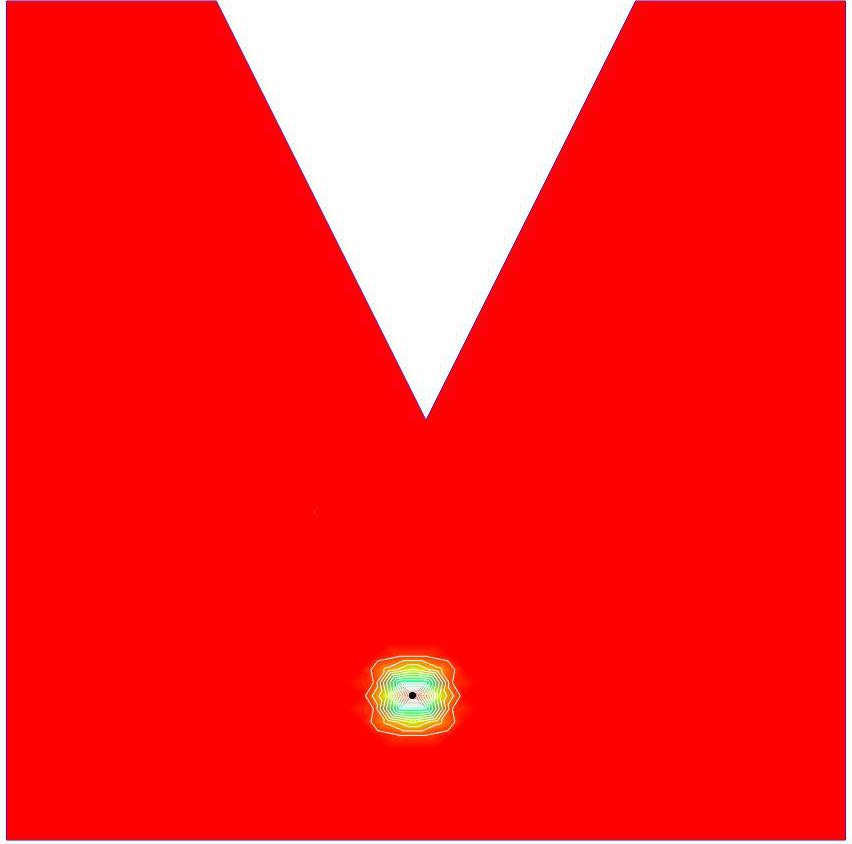}
&\includegraphics[height=41pt]{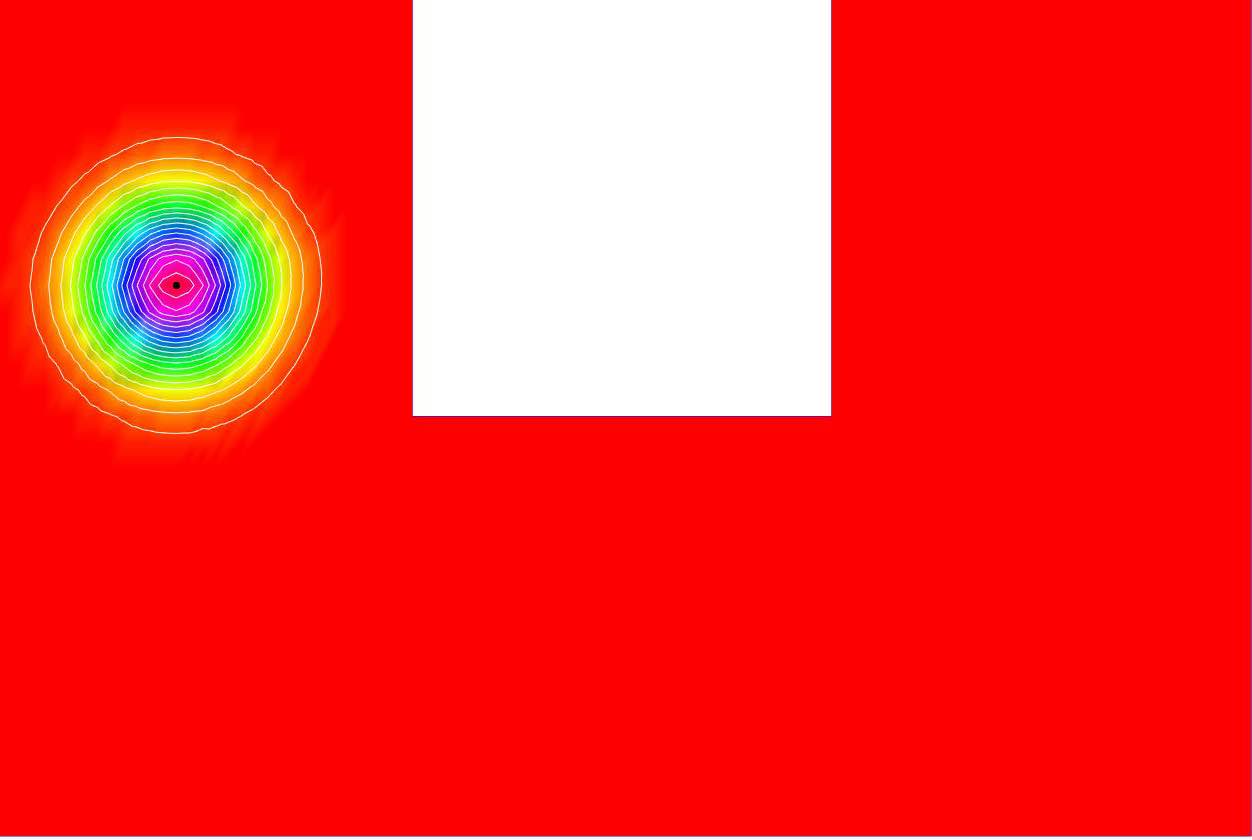}
&\includegraphics[height=41pt]{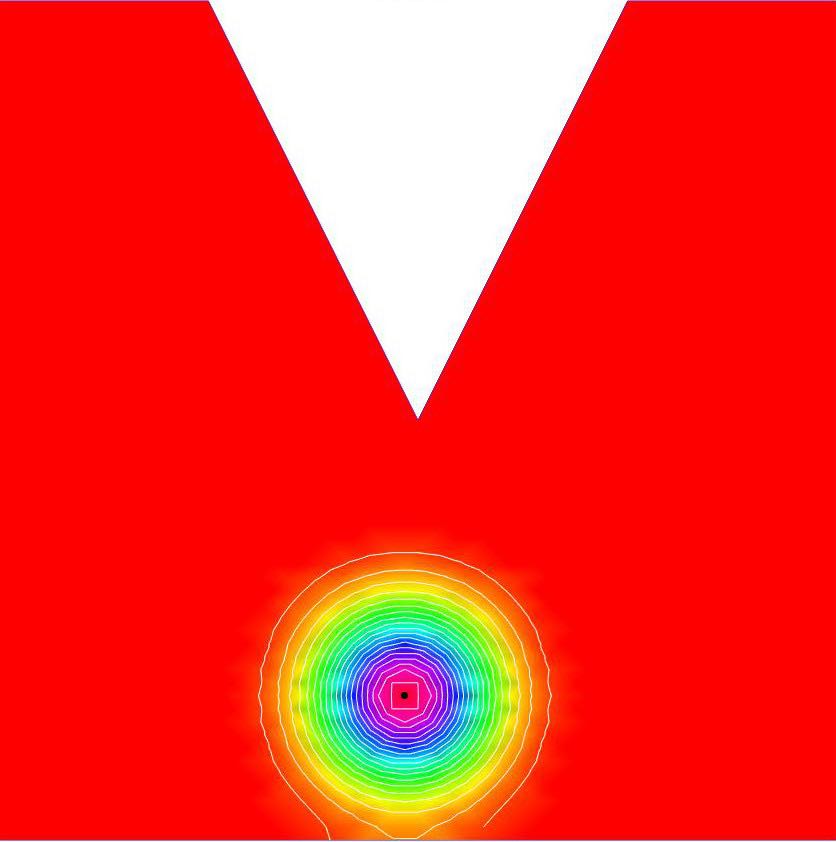}\\
\multicolumn{2}{c|}{$t=5\times 10^{-2}$}
&\multicolumn{2}{c}{$t=5\times 10^{-1}$}\\
\includegraphics[height=41pt]{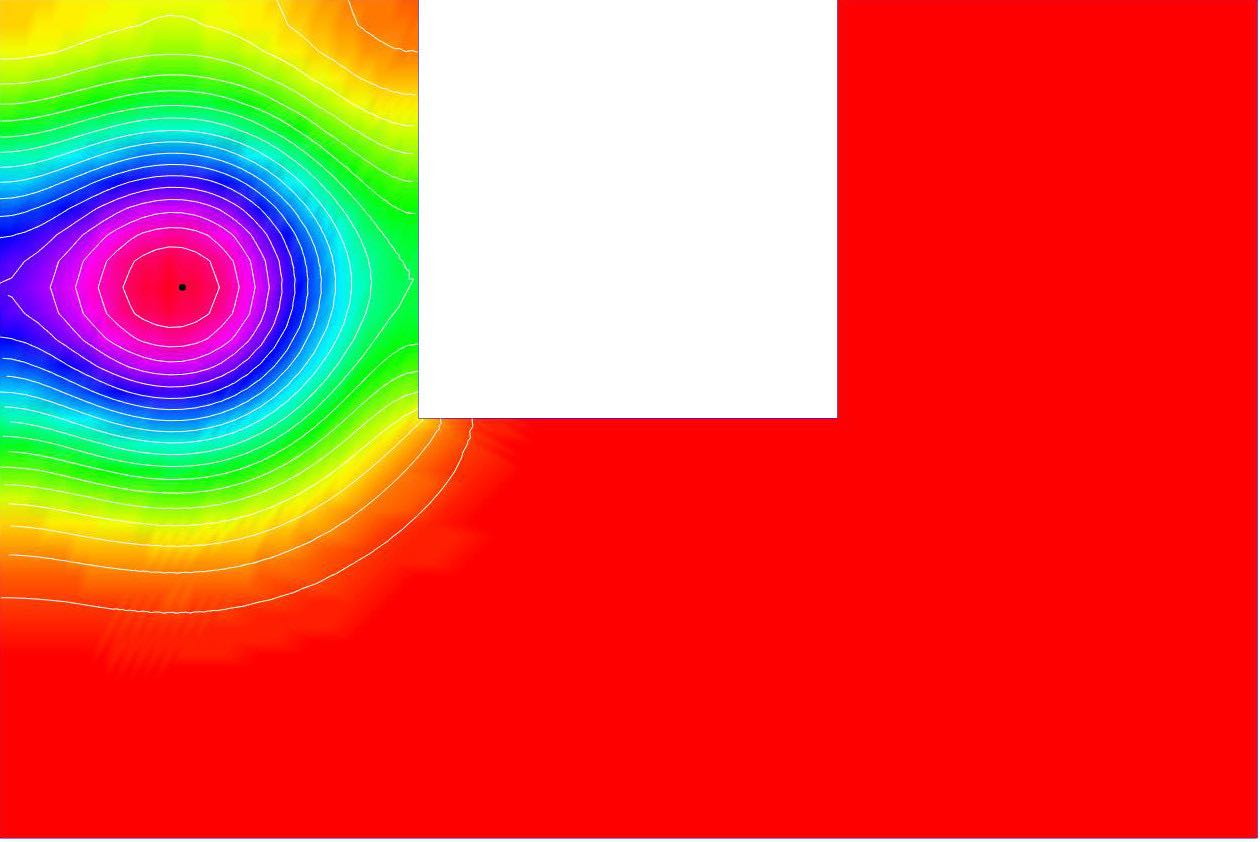}
&\includegraphics[height=41pt]{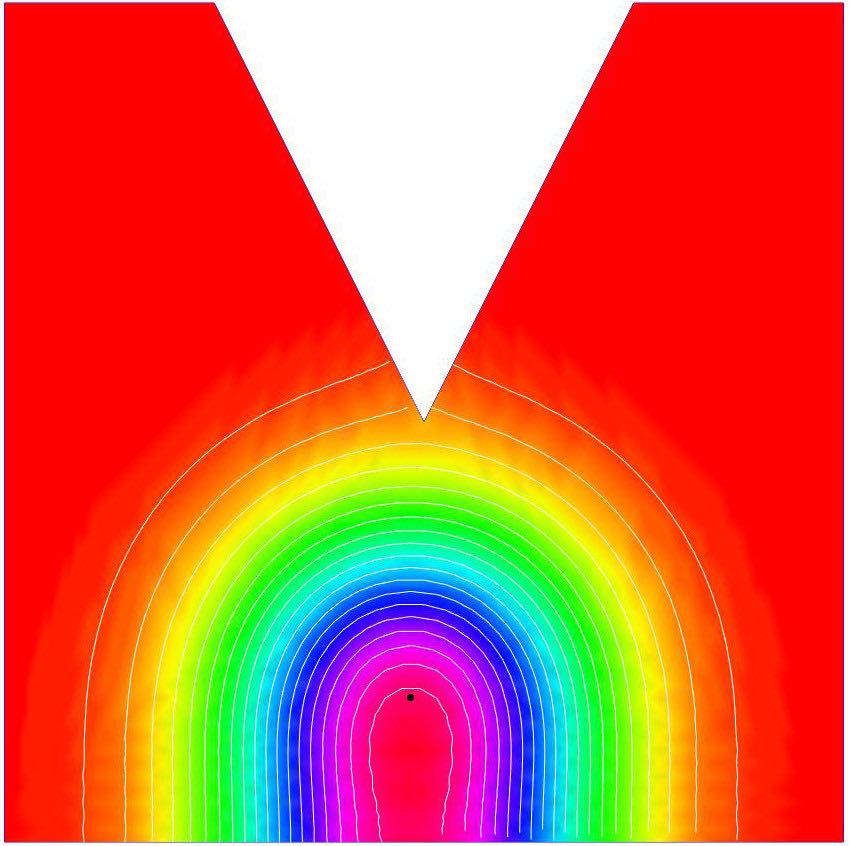}
&\includegraphics[height=41pt]{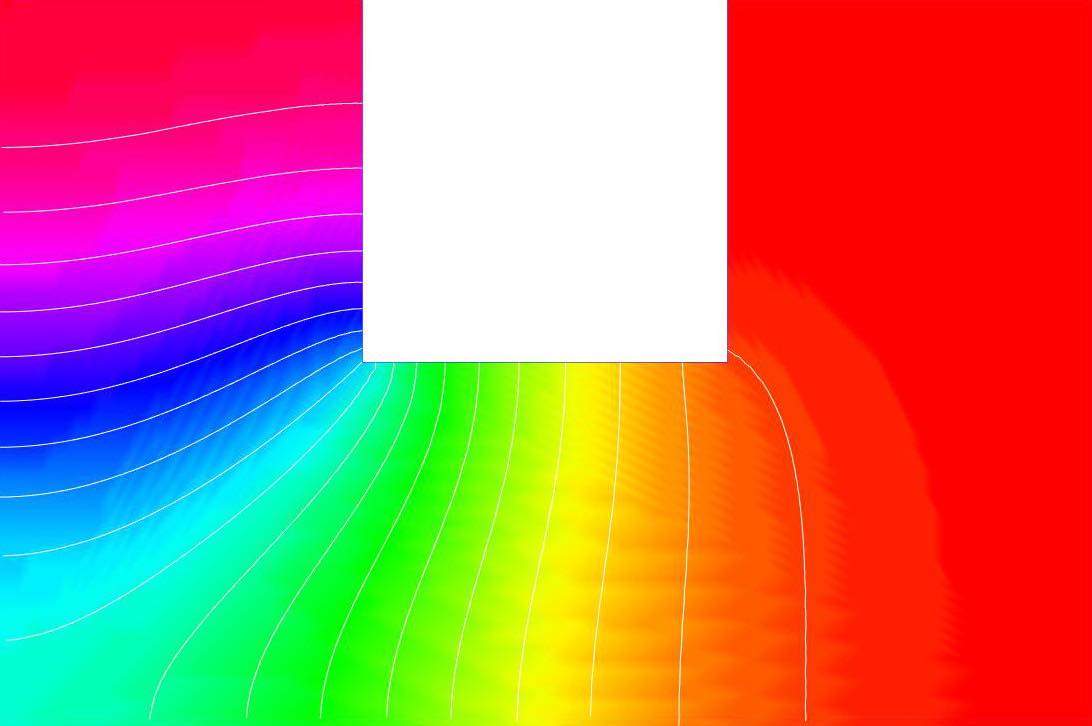}
&\includegraphics[height=41pt]{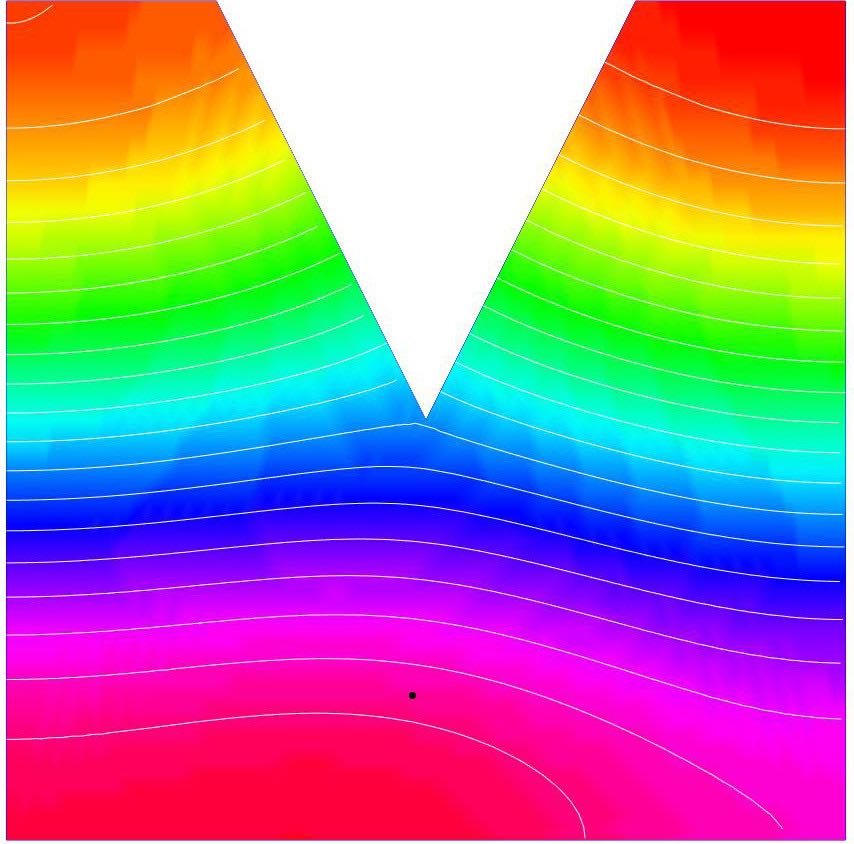}
\end{tabular}
\caption{Level-sets of data-driven diffusive membership functions centred at a seed point and at different scales~$t$.\label{fig:2D-DIFFUSION}}
\end{figure}
\section{Space of continuous F-transforms\label{sec:INT-OPERATOR-SP}}
We define the space of continuous F-transforms (Sect.~\ref{sec:SPACE-FT}) and different classes of membership functions (Sect.~\ref{KERNEL-EXAMPLE}). Data-driven membership functions are introduced in Sect.~\ref{sec:CONT-DATA-FT-SHORT}.
\begin{figure}[t]
\centering
(a)\\
\includegraphics[height=50pt]{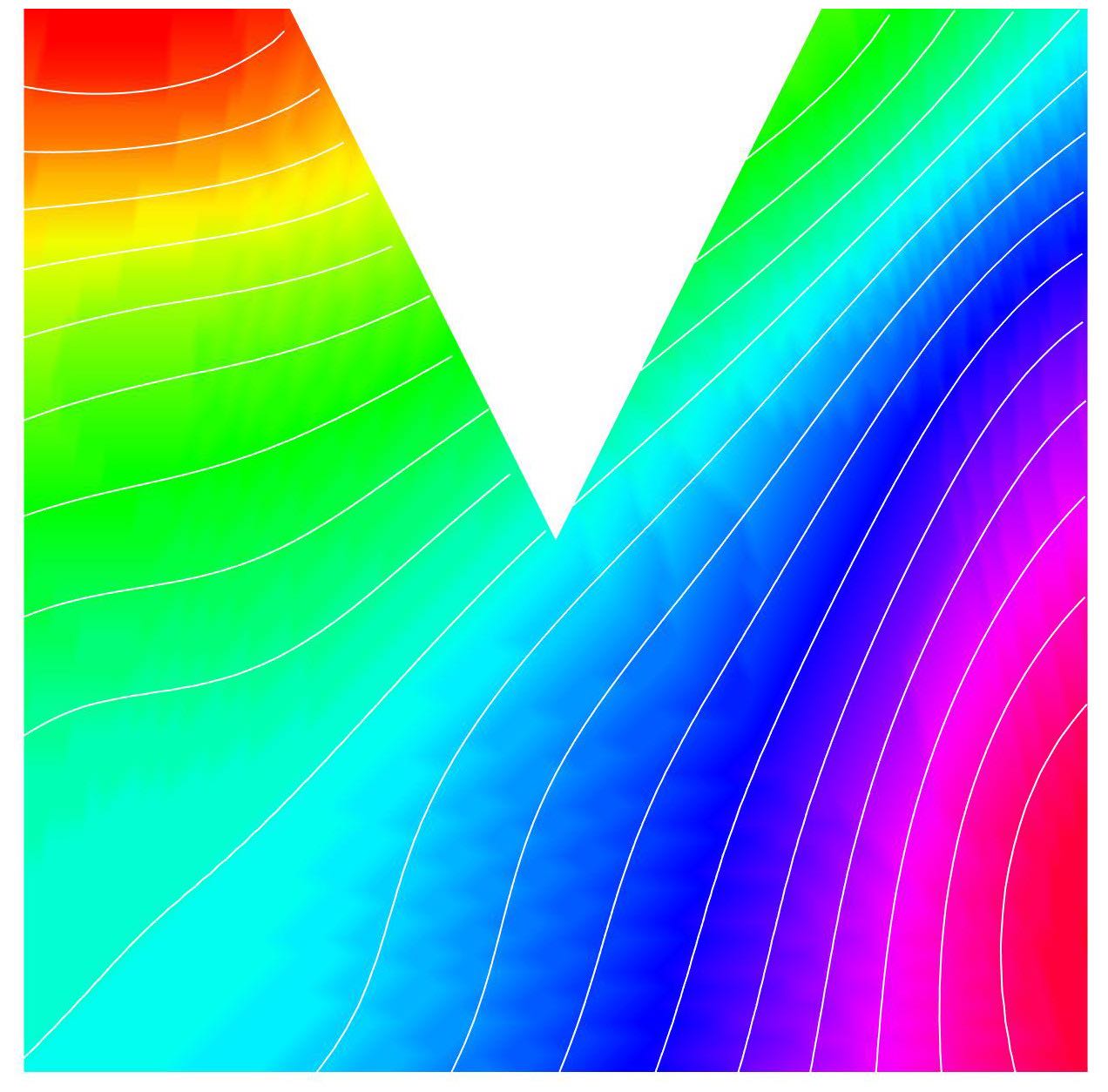}
\begin{tabular}{cccc}
\hline
$\alpha:=5\%$ &$\alpha:=10\%$ &$\alpha:=50\%$ &$\alpha:=100\%$\\
\hline
\multicolumn{4}{c}{$f$}\\
\includegraphics[height=50pt]{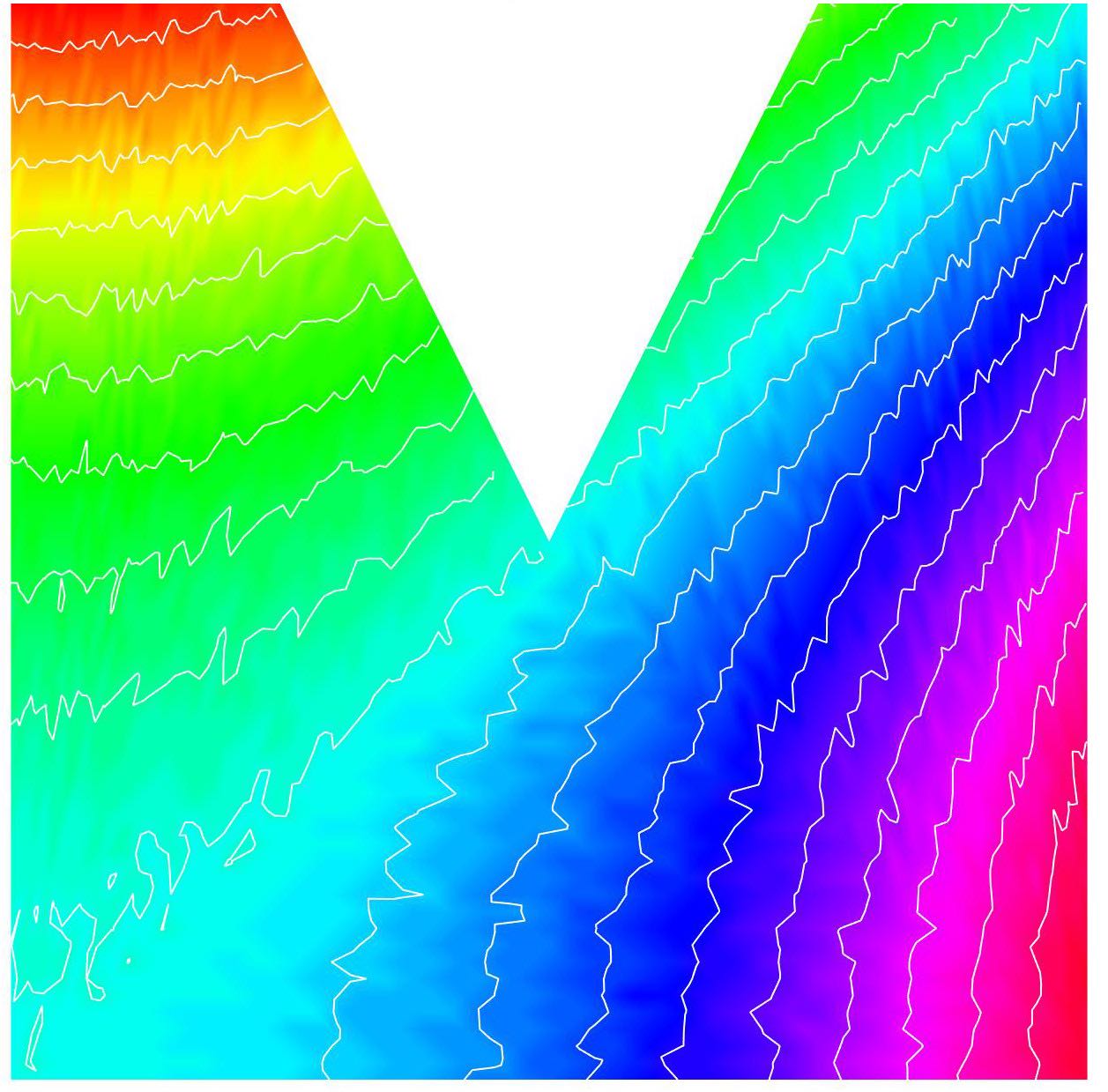}
&\includegraphics[height=50pt]{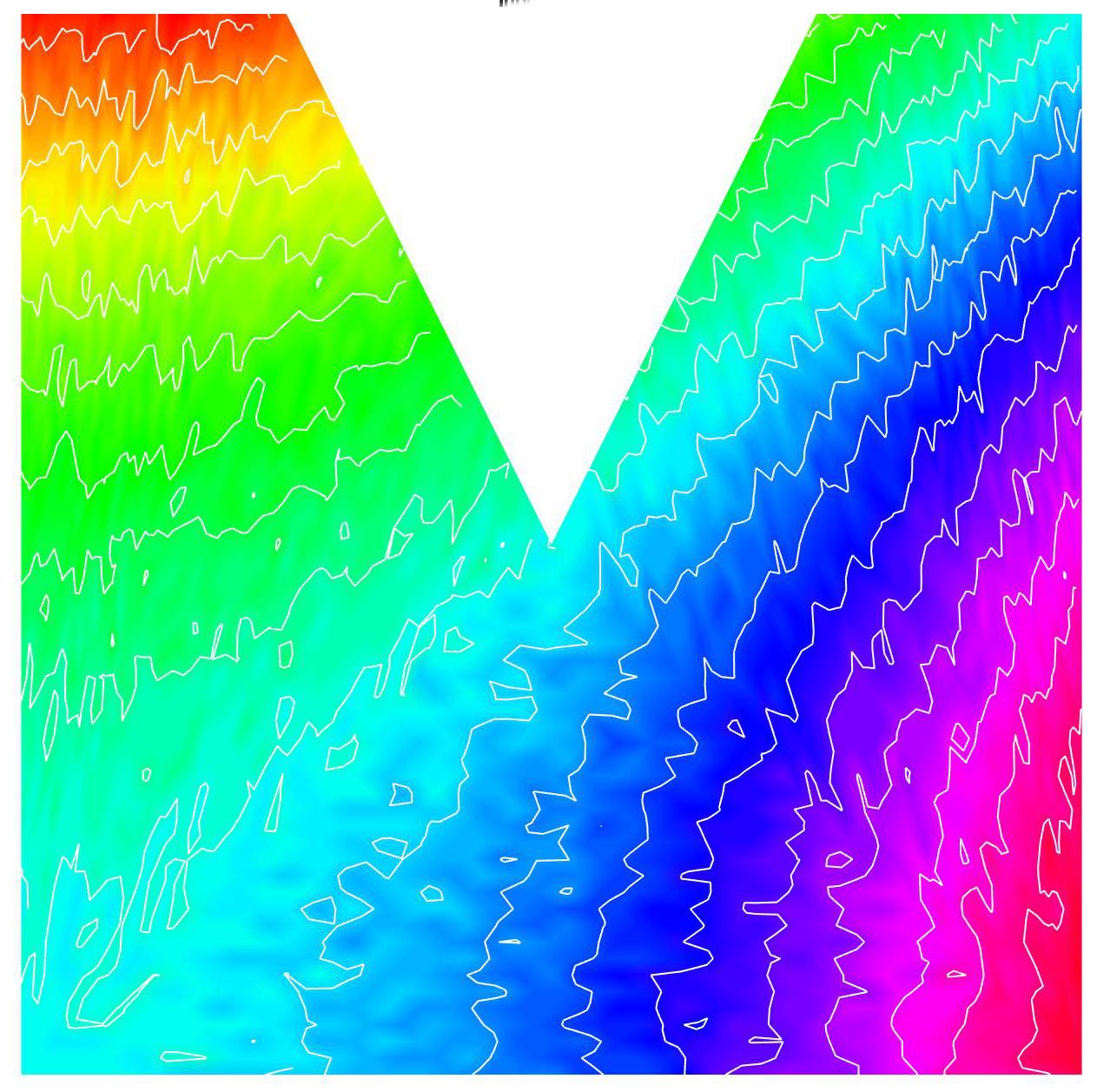}
&\includegraphics[height=50pt]{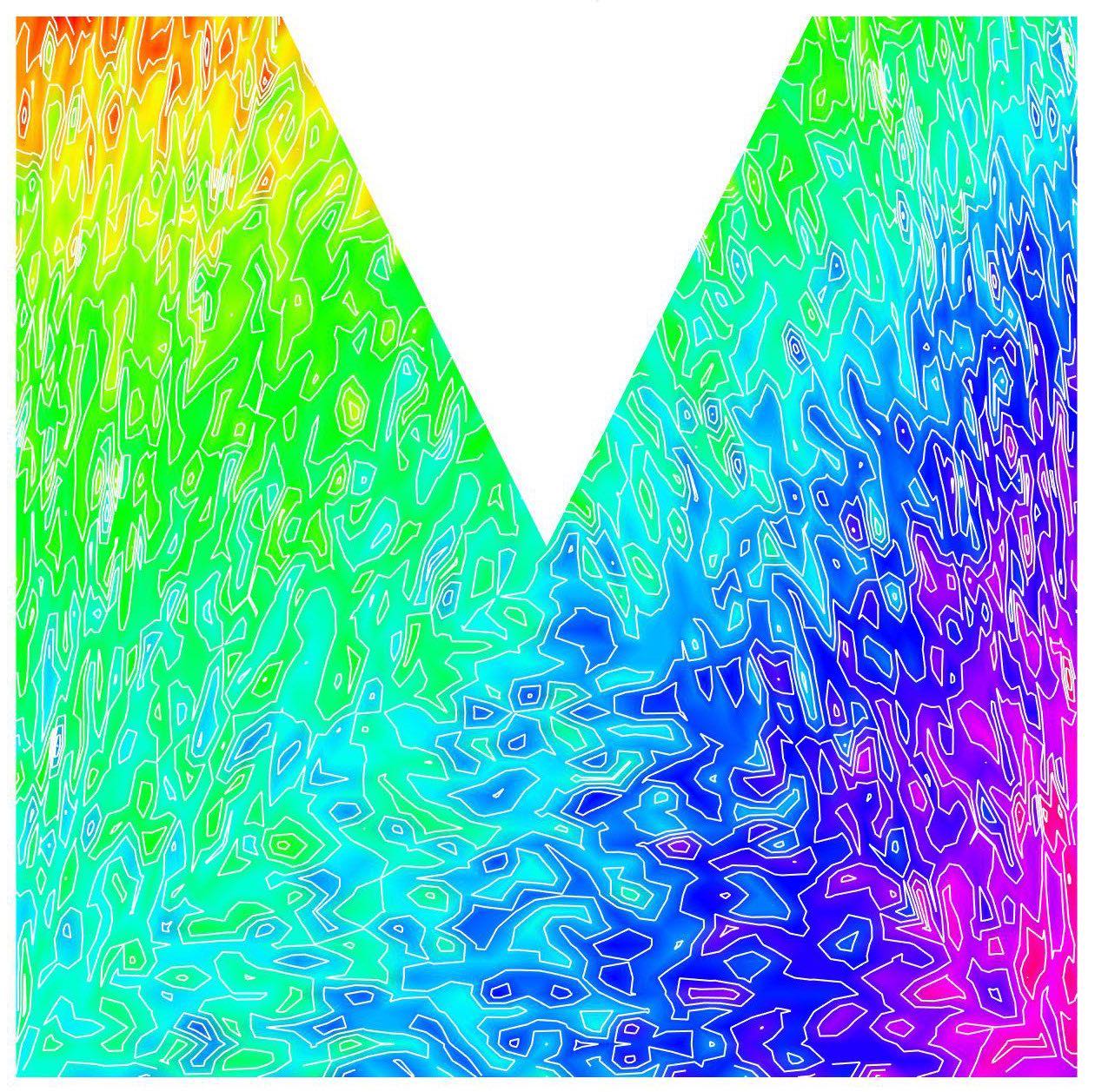}
&\includegraphics[height=50pt]{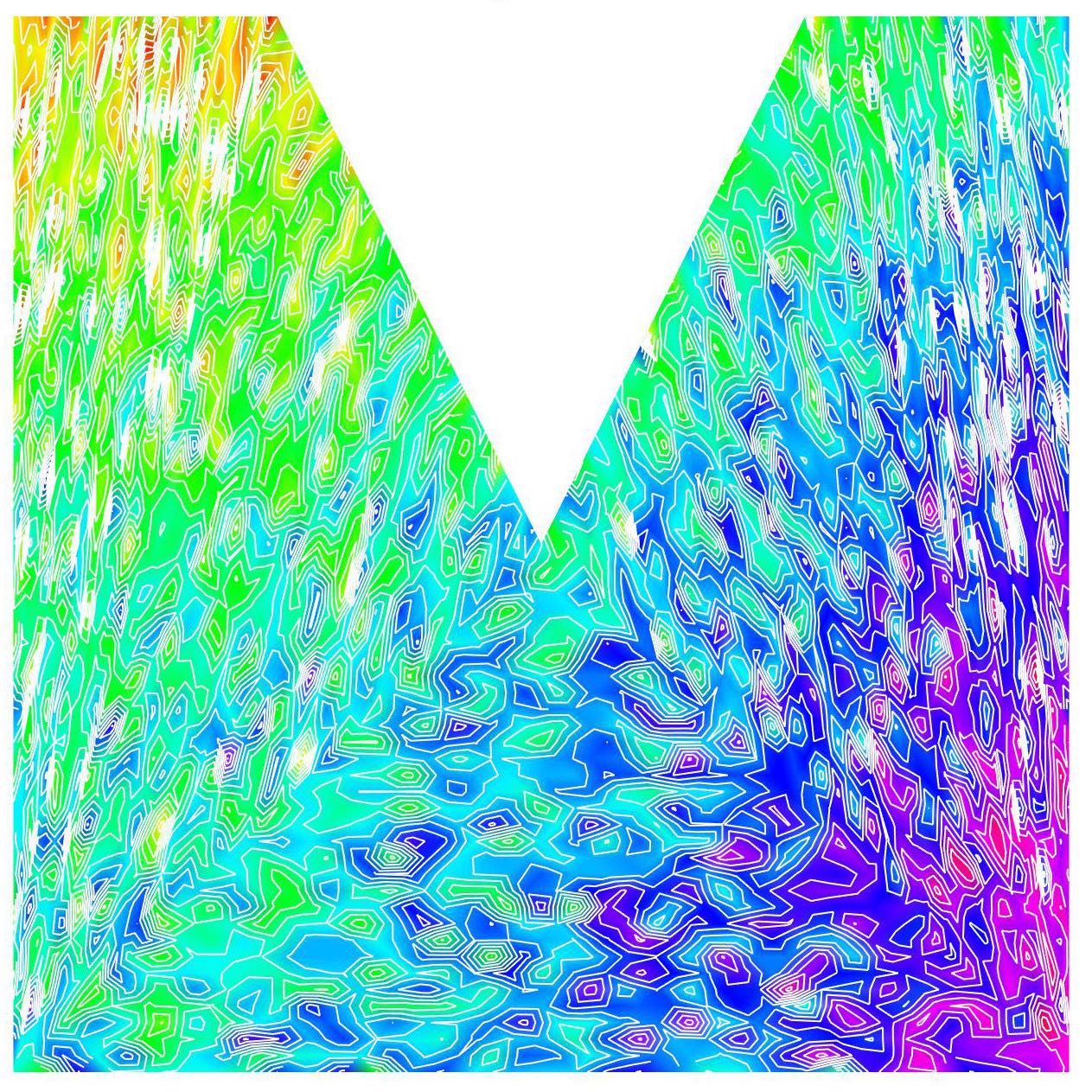}\\
\hline
\multicolumn{4}{c}{$\mathcal{F}f$}\\
\includegraphics[height=50pt]{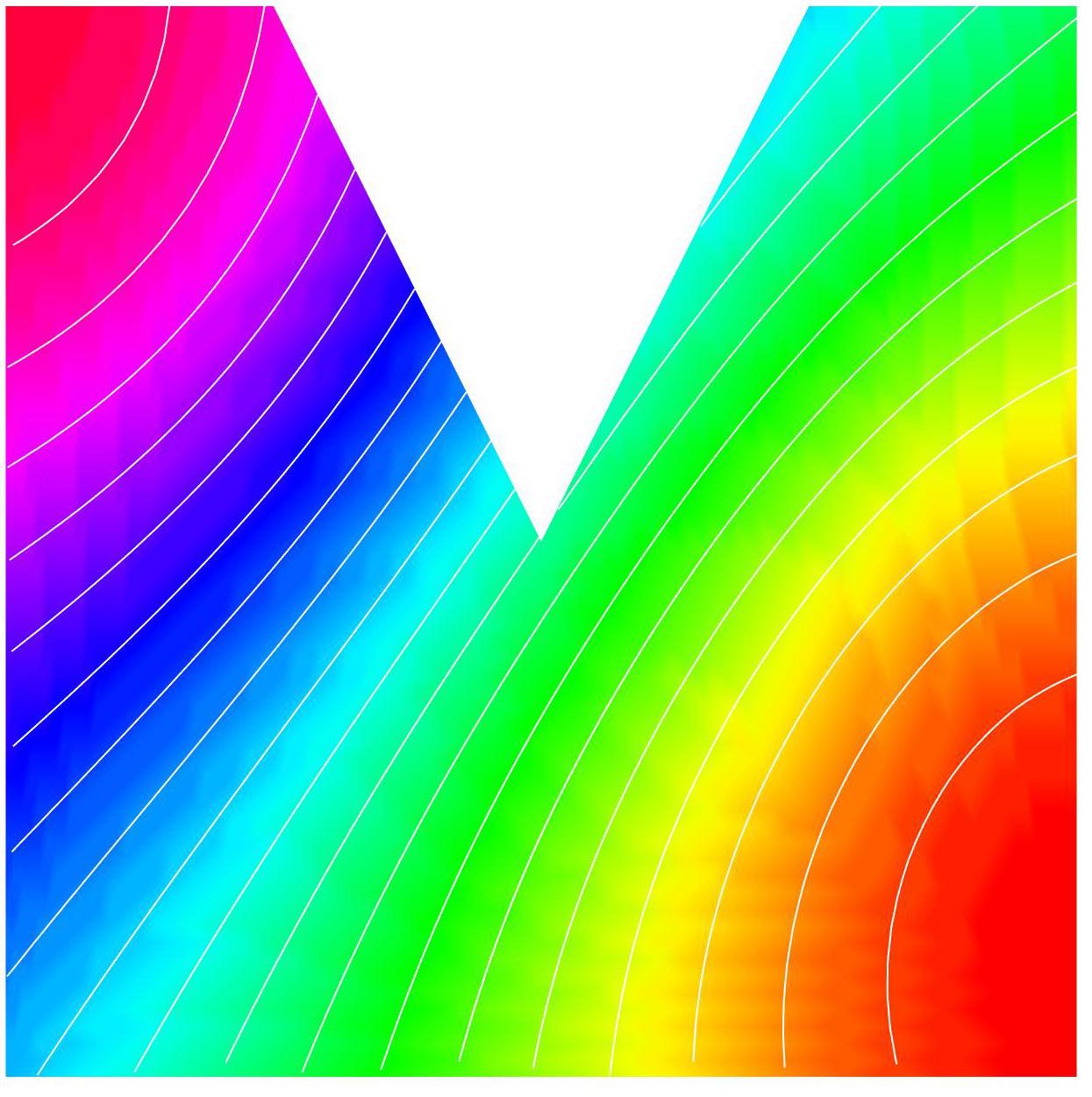}
&\includegraphics[height=50pt]{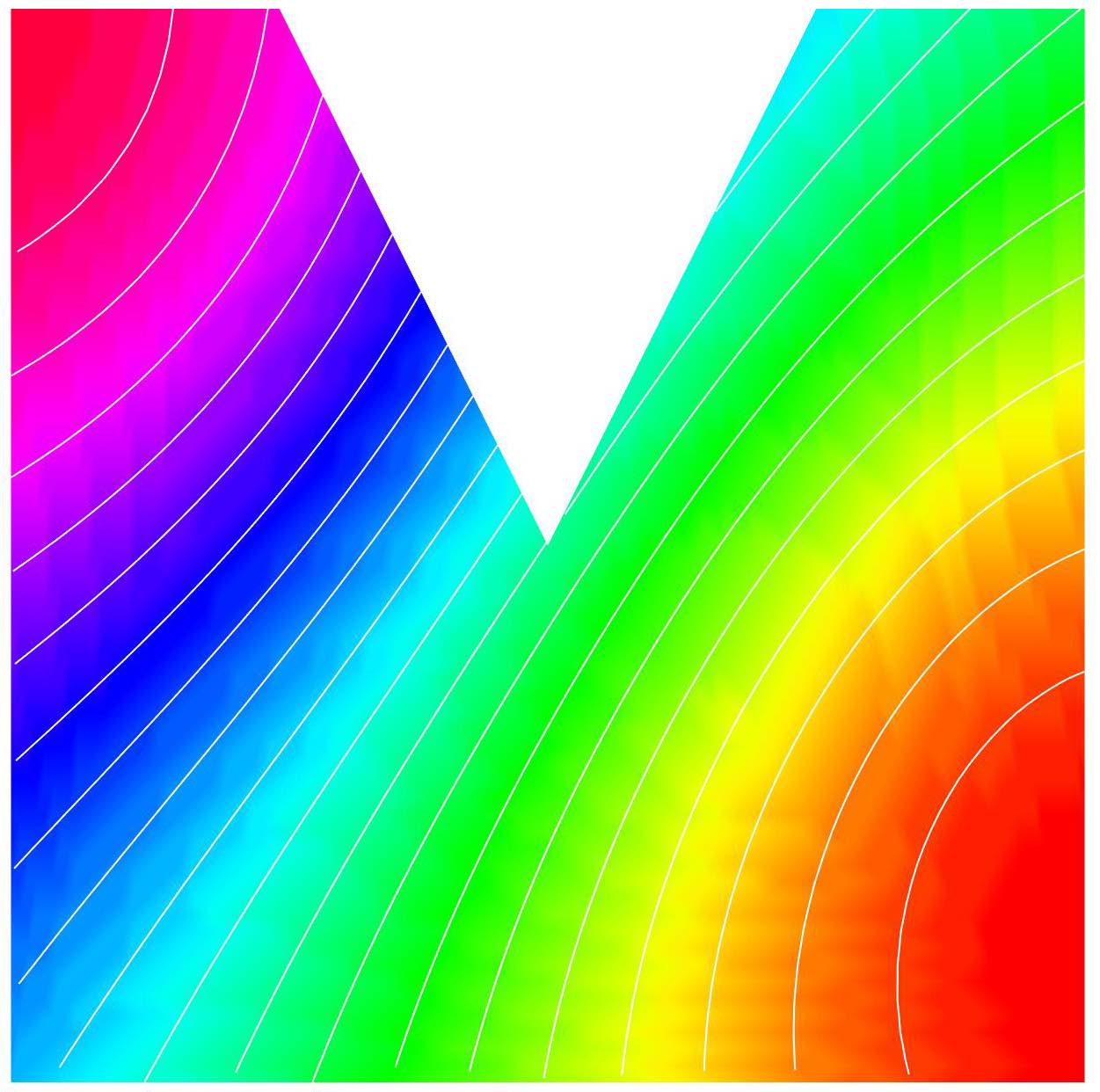}
&\includegraphics[height=50pt]{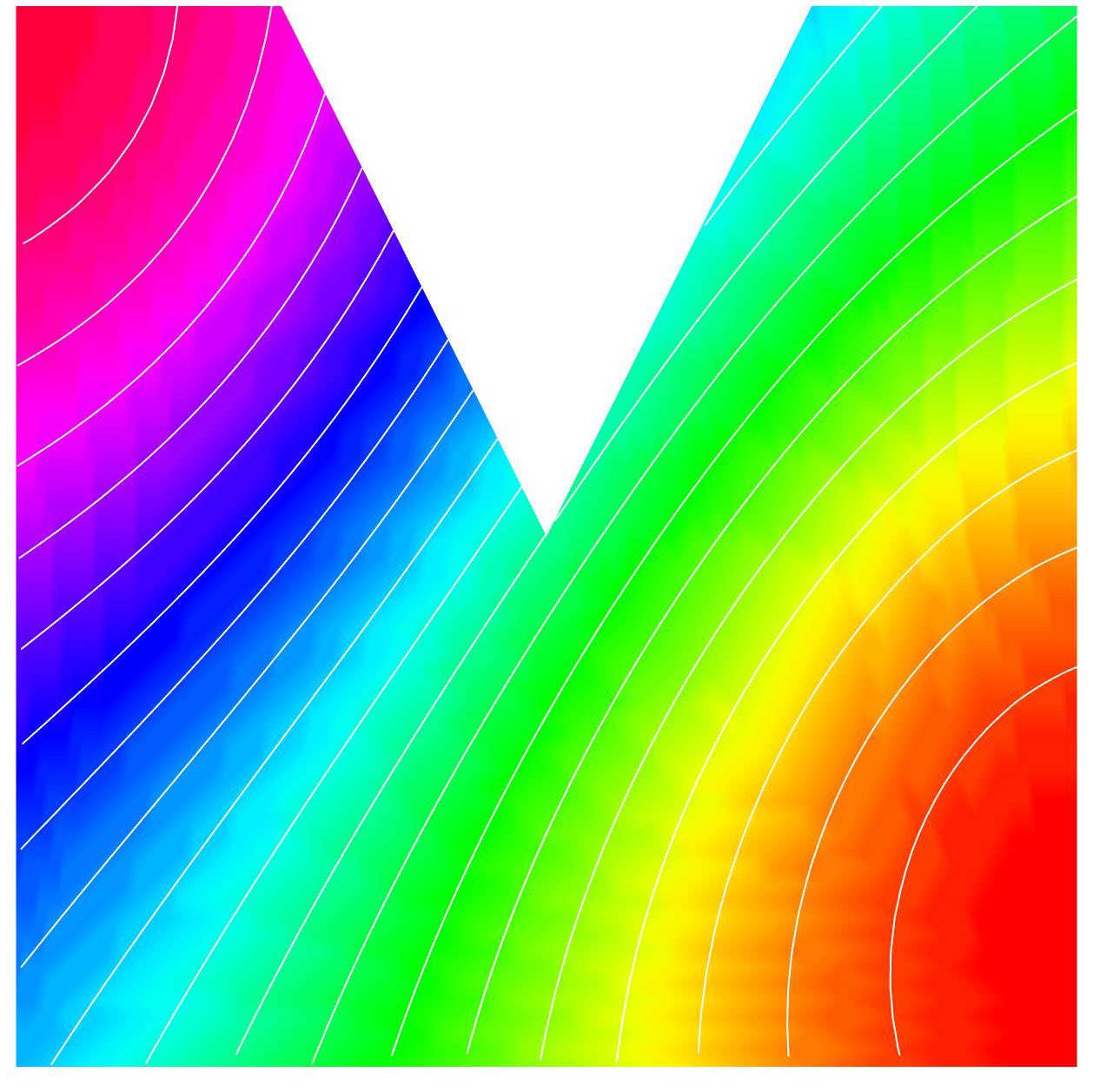}
&\includegraphics[height=50pt]{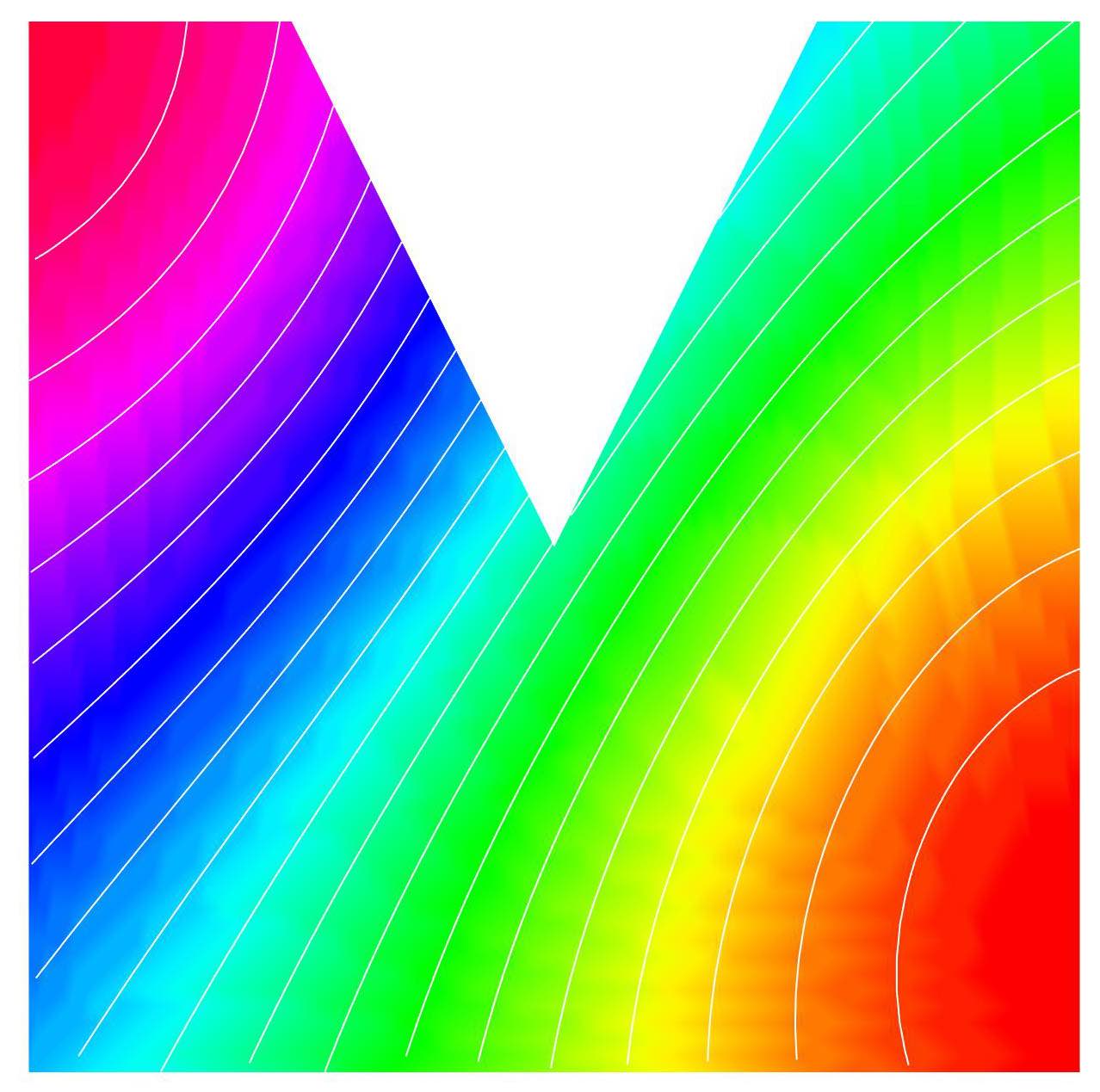}\\
\hline
\multicolumn{4}{c}{$\mathcal{F}^{-1}(\mathcal{F}f)$}\\
\includegraphics[height=50pt]{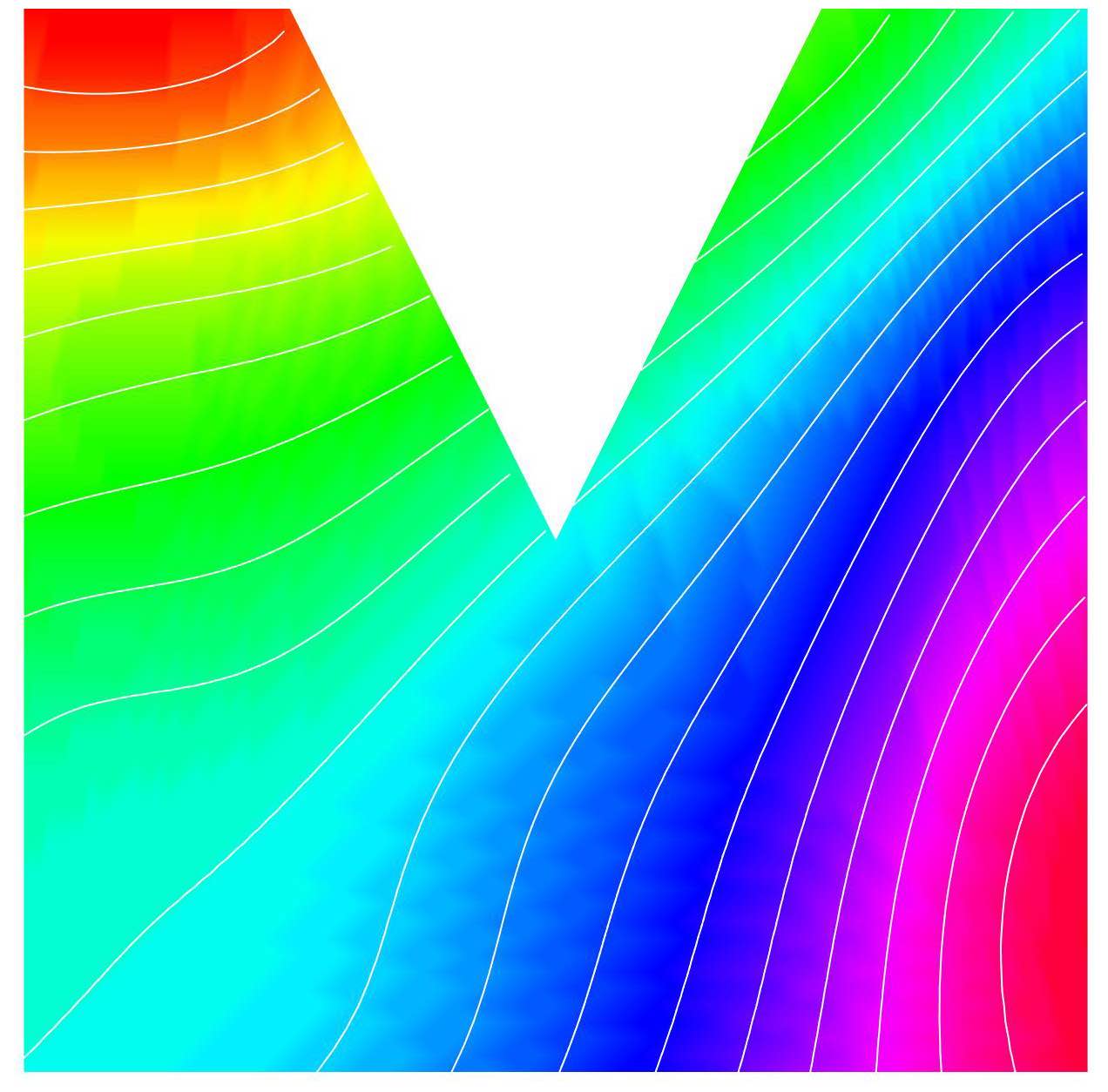}
&\includegraphics[height=50pt]{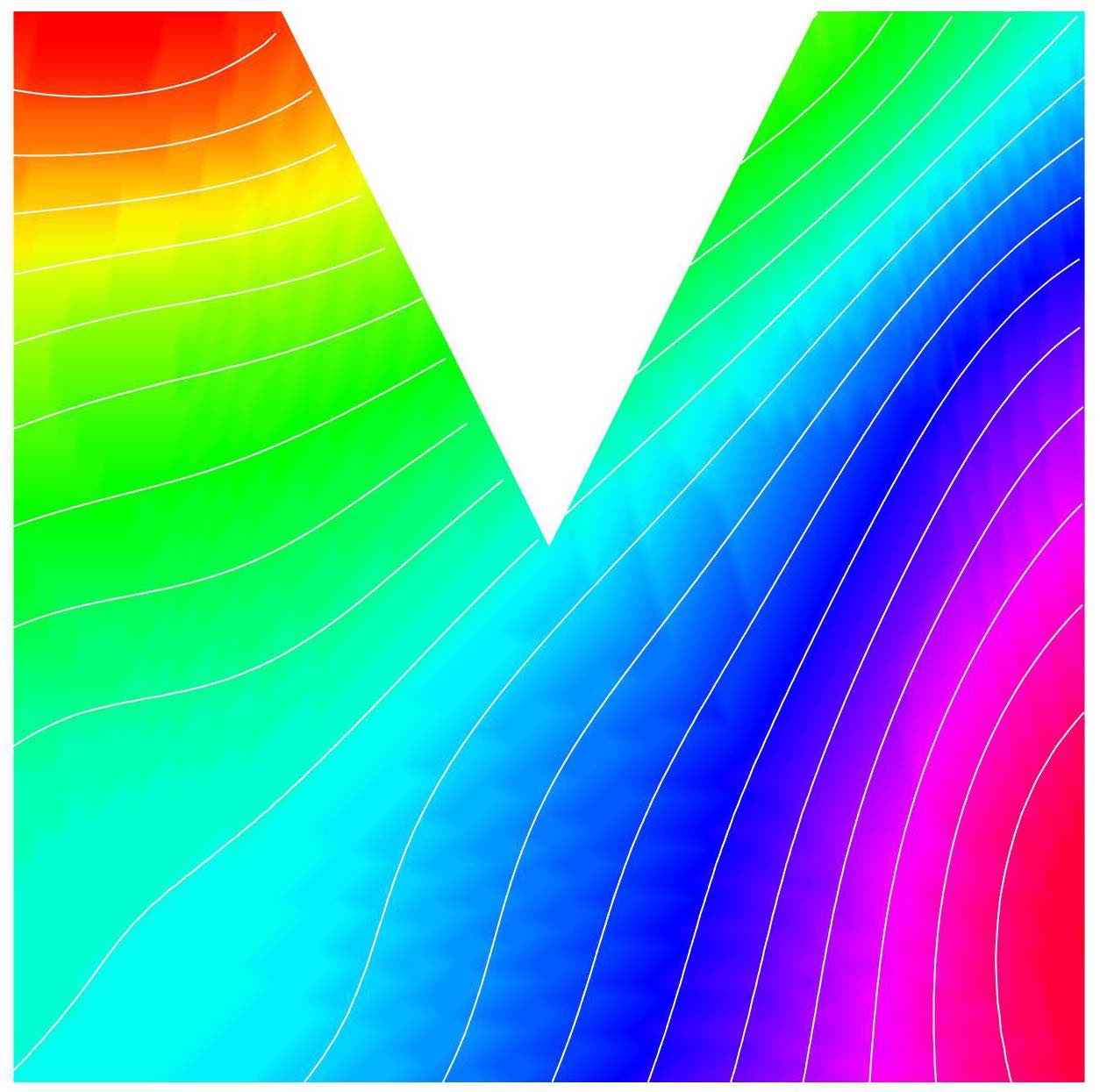}
&\includegraphics[height=50pt]{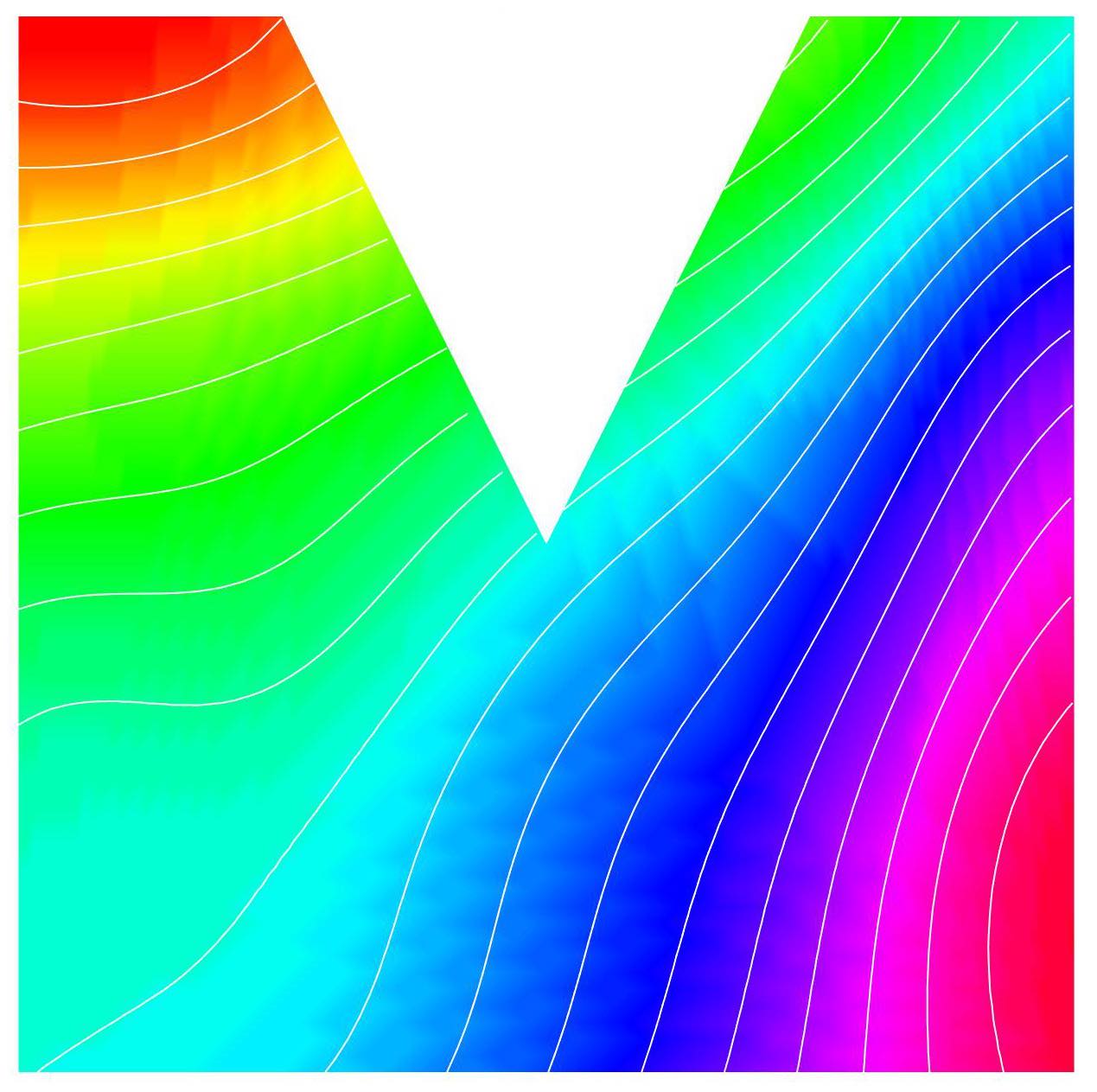}
&\includegraphics[height=50pt]{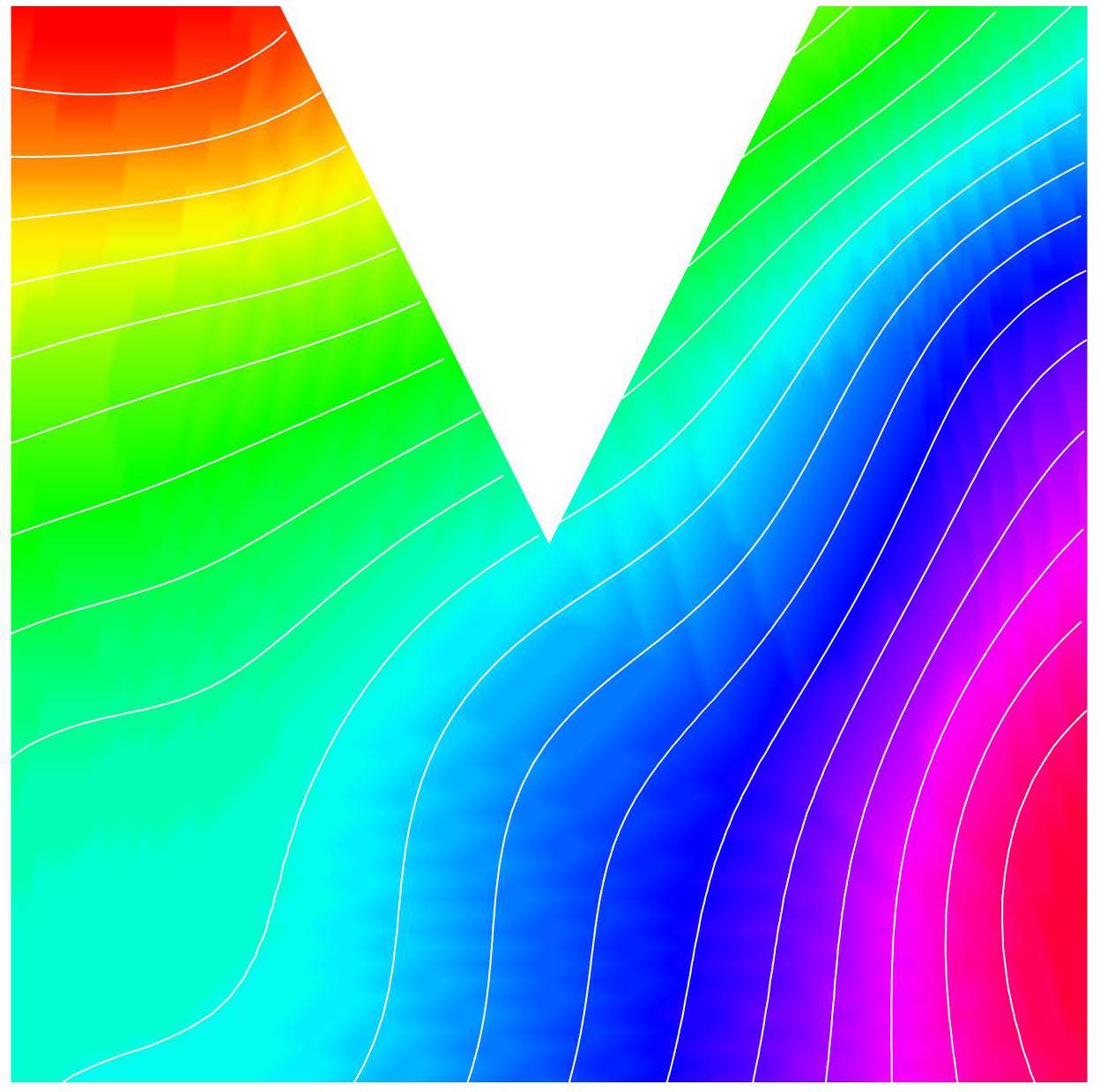}\\
$\epsilon_{\infty}:=2.6\%$		&$\epsilon_{\infty}:=4.1\%$		&$\epsilon_{\infty}:=5.8\%$		&$\epsilon_{\infty}:=7.8\%$\\
\hline
(b) &(c) &(d) &(e)
\end{tabular}
\caption{(First row) Level-sets and color-map of noisy signals with an increasing error magnitude~$\alpha$ (from (b) to (e)) achieved by adding a Gaussian noise to (a) an input signal~$f$. Level-sets (second row) of the F-transform \mbox{$\mathcal{F}f$} induced by the multi-quadratic kernel and (third row) of the reconstructed signal \mbox{$\mathcal{F}^{-1}(\mathcal{F}f)$}. The reconstruction error is defined as \mbox{$\epsilon_{\infty}:=\|f-\mathcal{F}^{-1}(\mathcal{F}f)\|_{\infty}/\|f\|_{\infty}$}.\label{fig:2D-FT-IFT-KERNELS-NOISE}}
\end{figure}
\subsection{Space of continuous F-transforms\label{sec:SPACE-FT}}
Given \mbox{$\mathcal{H}(\Omega):=\mathcal{C}^{0}(\Omega)
\cap\mathcal{L}^{2}(\Omega)$}, we define the \emph{linear space of continuous F-transforms} as
\begin{equation*}
\mathcal{F}:=
\{\mathcal{L}_{K}:\mathcal{H}(\Omega)\rightarrow\mathcal{C}^{0}(\Omega),\quad K\in\mathcal{L}^{2}(\Omega\times\Omega)\}.
\end{equation*}
In fact, \mbox{$\alpha\mathcal{L}_{K_{1}}+\beta\mathcal{L}_{K_{2}}=\mathcal{L}_{\alpha K_{1}+\beta K_{2}}$} and any couple of kernels \mbox{$K_{1},K_{2}:\Omega\times\Omega\rightarrow\mathbb{R}$} allows us to generate a new continuous F-transform induced by a linear combination of these kernels. The spase~$\mathcal{F}$ inherits a Hilbert structure with respect to the scalar product \mbox{$\langle\mathcal{L}_{H},\mathcal{L}_{K}\rangle_{2}=\langle H,K\rangle_{2}$}. Through  the norm \mbox{$\|\mathcal{L}_{H}-\mathcal{L}_{K}\|_{2}=\|H-K\|_{2}$}, we can compare two continuous F-transforms and express the convergence of~$(\mathcal{L}_{K_{n}})_{n=0}^{+\infty}$ to~$\mathcal{L}_{K}$ in terms of the convergence of \mbox{$(K_{n})_{n=0}^{+\infty}$} to~$K$ in \mbox{$\mathcal{L}^{2}(\Omega\times\Omega)$}.

\textbf{Relation between membership function and integral kernel\label{sec:MEMBERSHIP-KERNEL}}
The properties of the membership functions (e.g., square integrability, continuity, symmetry, positiveness) are inherited by the corresponding integral kernel. The normalised membership function, or equivalently the integral kernel (\ref{eq:NORMALISED-KERNEL}), is positive (by definition) and square integrable, according to the upper bound
\begin{equation*}
\|K\|_{2}
=\left\|A(\cdot,\cdot)/S(\cdot)\right\|_{2}
\leq_{\textrm{Eq.}(\ref{eq:UPPER-BOUND})}(A_{\max}\vert\Omega\vert^{1/2})/A_{\min}.
\end{equation*}
Since the upper bound
\begin{equation*}
\begin{split}
\vert K(\mathbf{p},\mathbf{q})&-K(\mathbf{p}_{0},\mathbf{q}_{0})\vert\\
&=\frac{\left\vert A(\mathbf{p},\mathbf{q})S(\mathbf{p}_{0})^{1/2}-A(\mathbf{p}_{0},\mathbf{q}_{0})S(\mathbf{p})^{1/2}\right\vert}{(S(\mathbf{p})S(\mathbf{p}_{0}))^{1/2}}\\
&\leq\frac{\left\vert A(\mathbf{p},\mathbf{q})S(\mathbf{p}_{0})^{1/2}-A(\mathbf{p}_{0},\mathbf{q}_{0})S(\mathbf{p})^{1/2}\right\vert}{A_{\min}\vert\Omega\vert}
\end{split}
\end{equation*}
converges to zero as \mbox{$(\mathbf{p},\mathbf{q})\rightarrow(\mathbf{p}_{0},\mathbf{q}_{0})$}, the integral kernel is continuous.  Analogous results apply to the kernel in Eq. (\ref{eq:NORMALISED-KERNEL-SYMMETRIC}).
\begin{figure*}[t]
\centering
\begin{tabular}{c|cccc}
\multicolumn{1}{c}{} &\multicolumn{4}{c}{}\\
(a)~$\alpha:=0$ &(b)~$\alpha:=5\%$ &(c)~$\alpha:=10\%$ &(d)~$\alpha:=50\%$ &(e)~$\alpha:=100\%$\\
\hline
$f$\includegraphics[height=60pt]{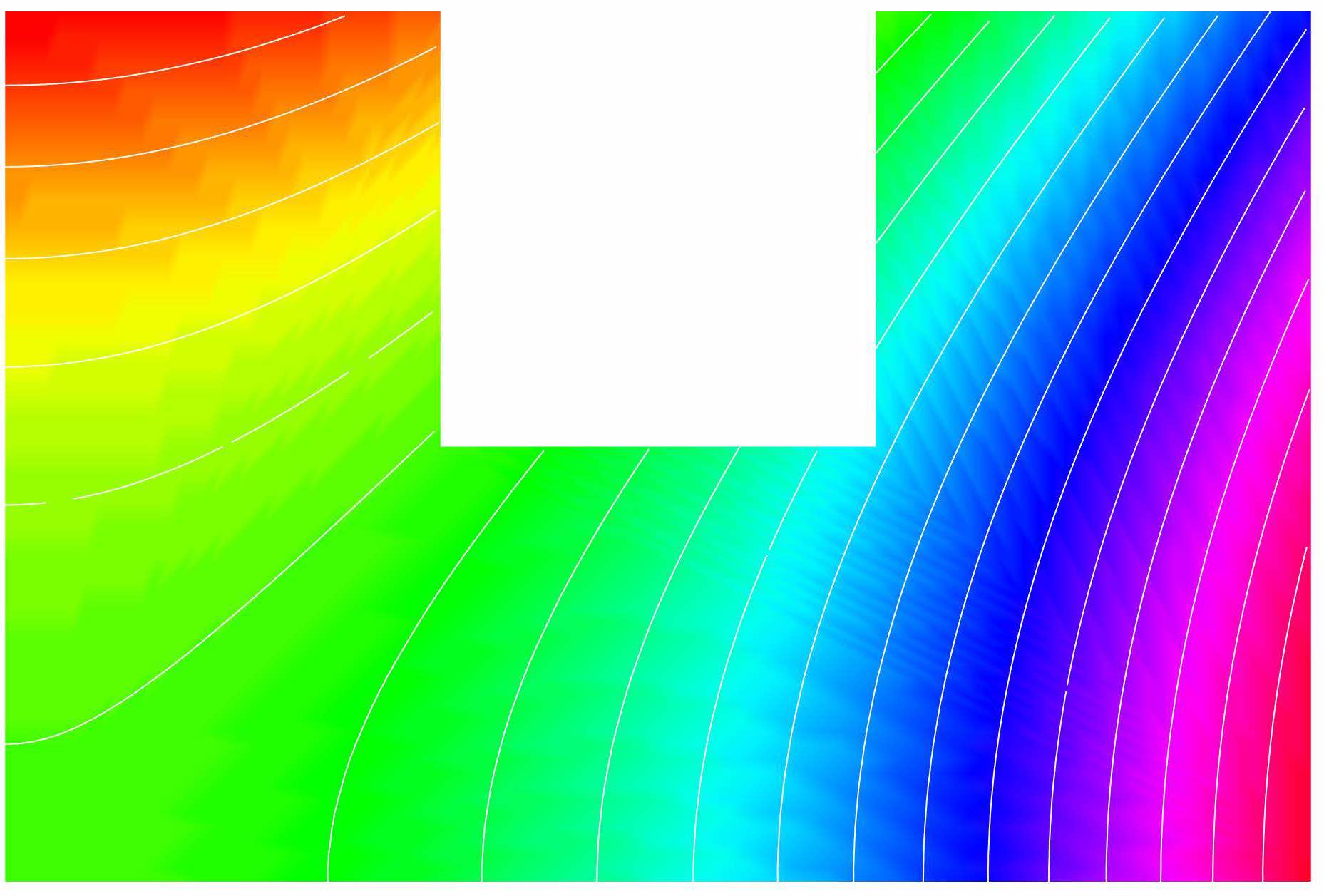}
&\includegraphics[height=60pt]{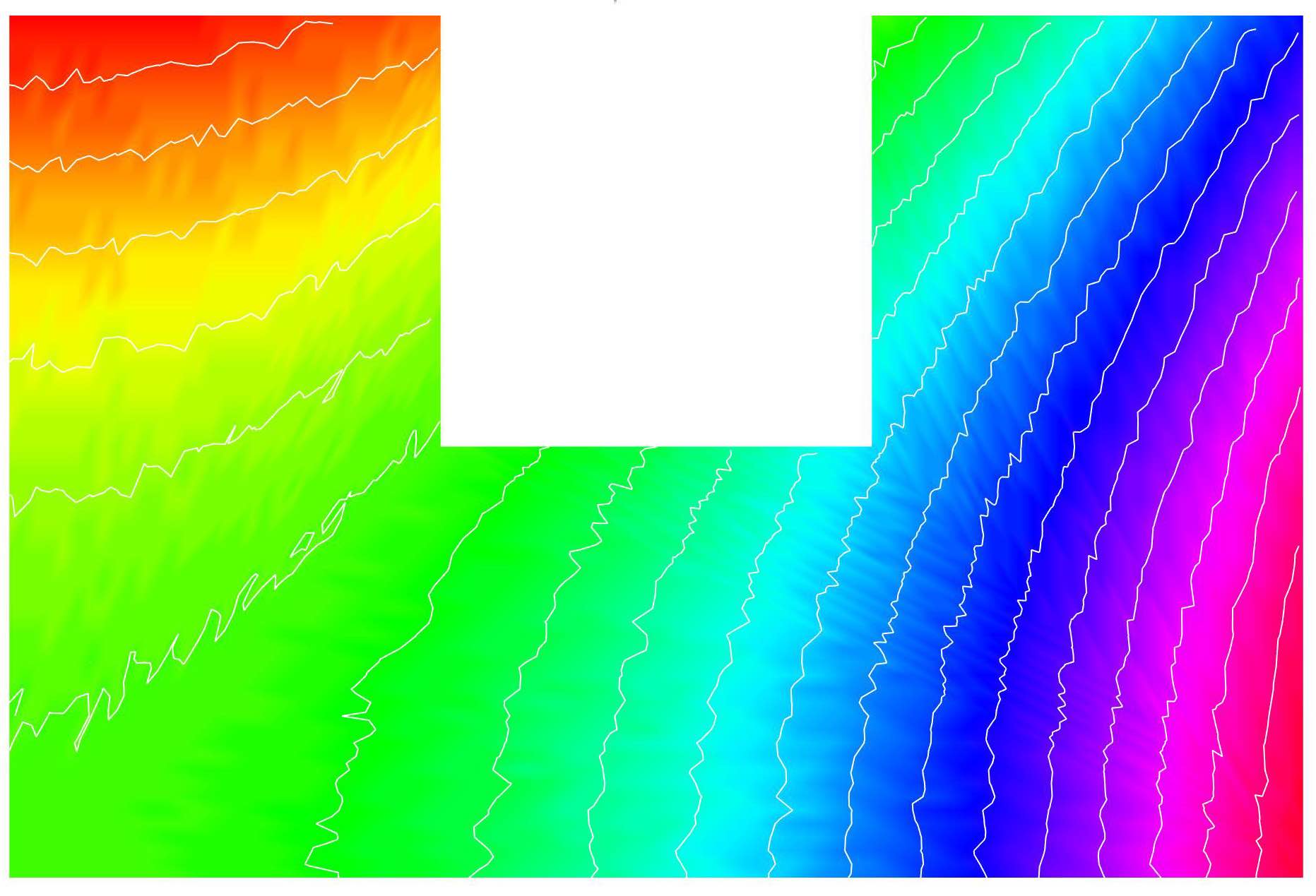}
&\includegraphics[height=60pt]{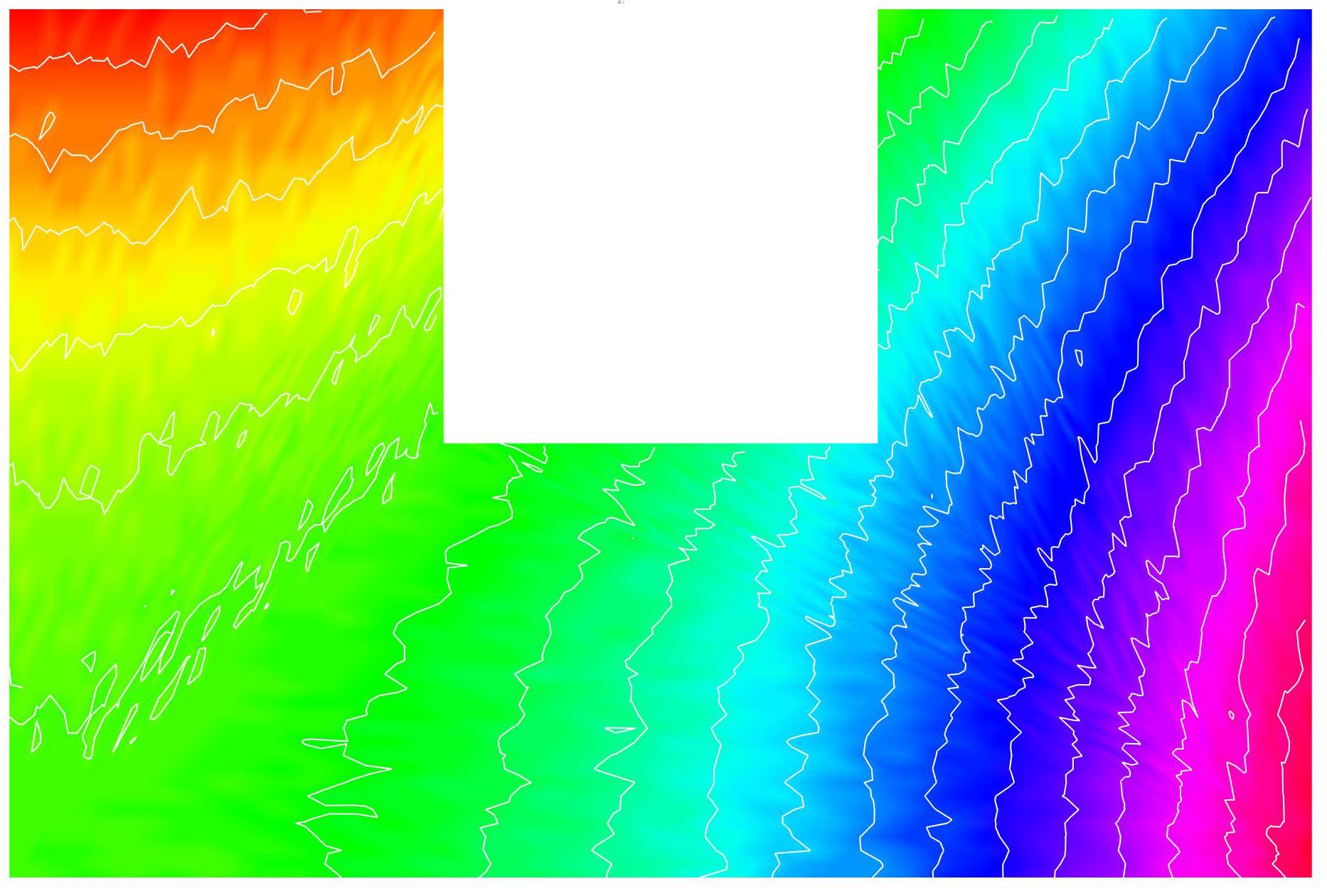}
&\includegraphics[height=60pt]{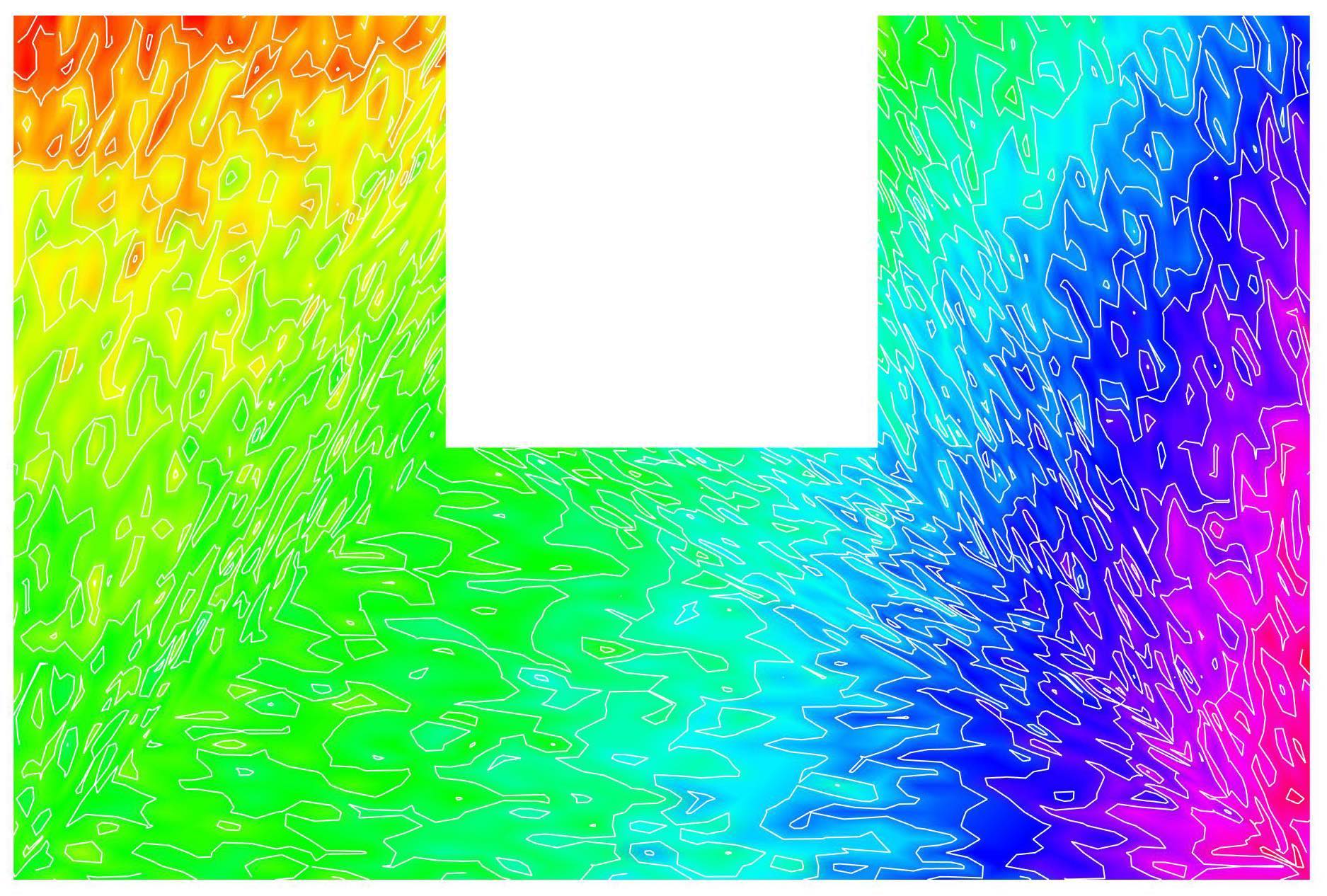}
&\includegraphics[height=60pt]{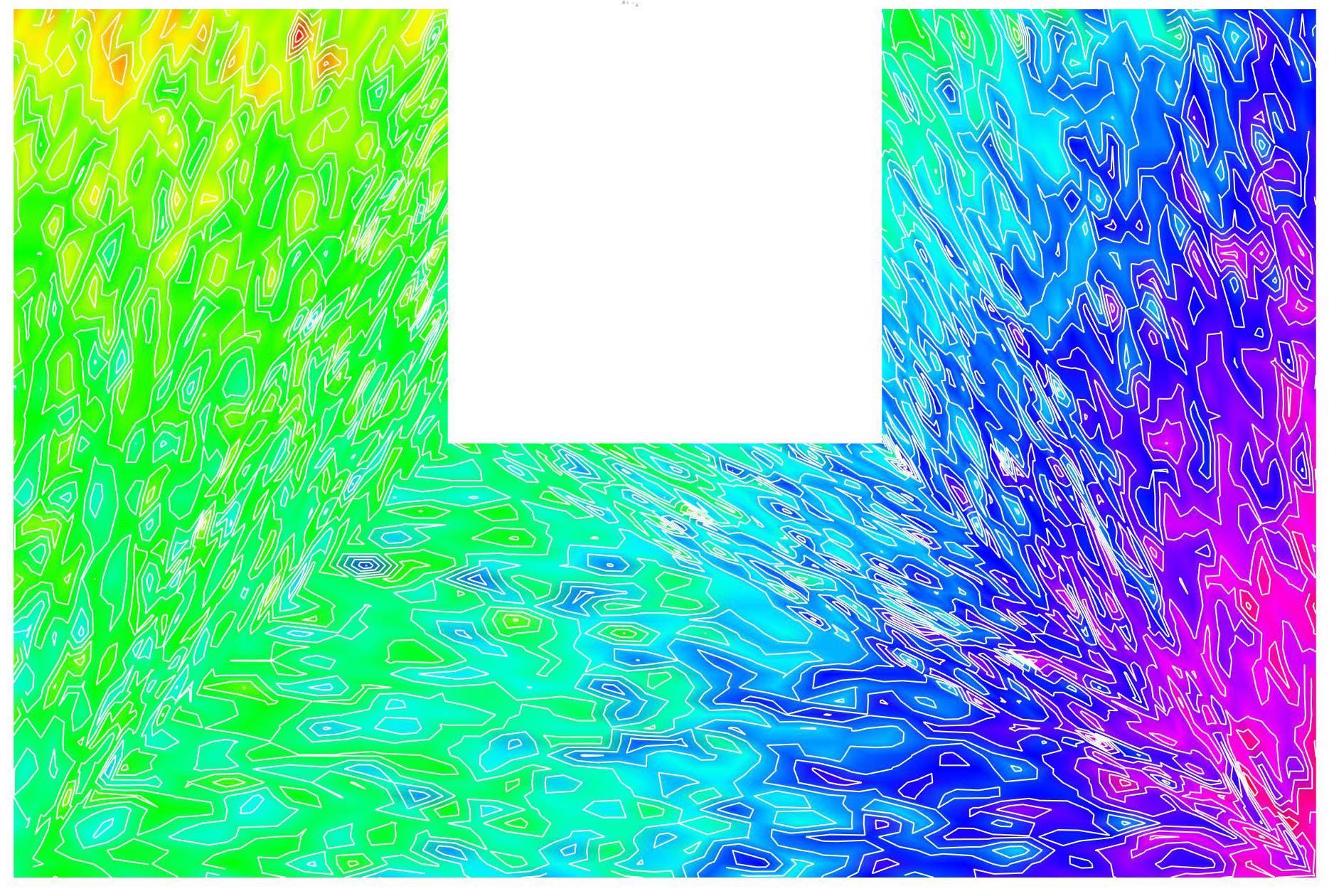}\\
\hline
&\includegraphics[height=60pt]{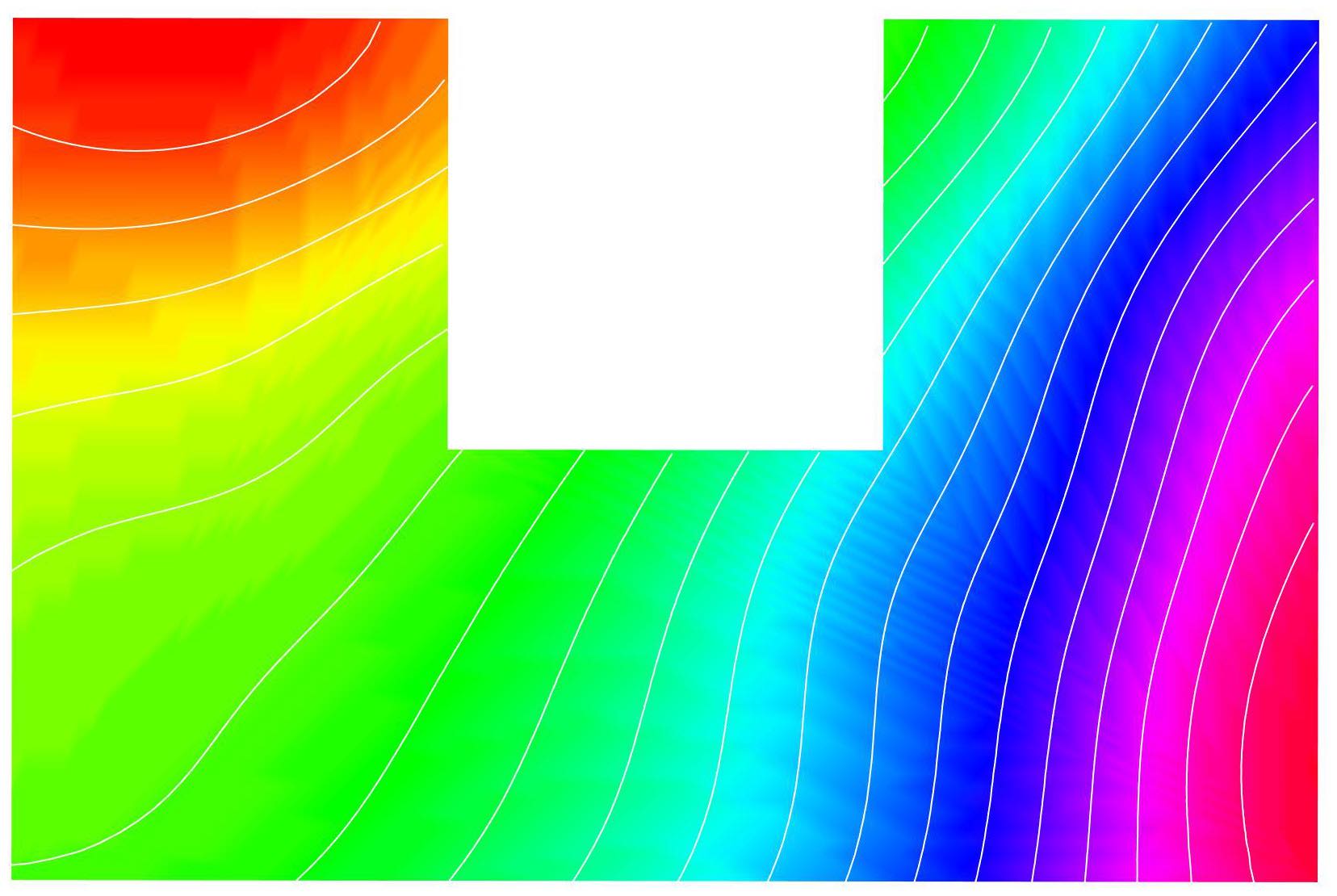}
&\includegraphics[height=60pt]{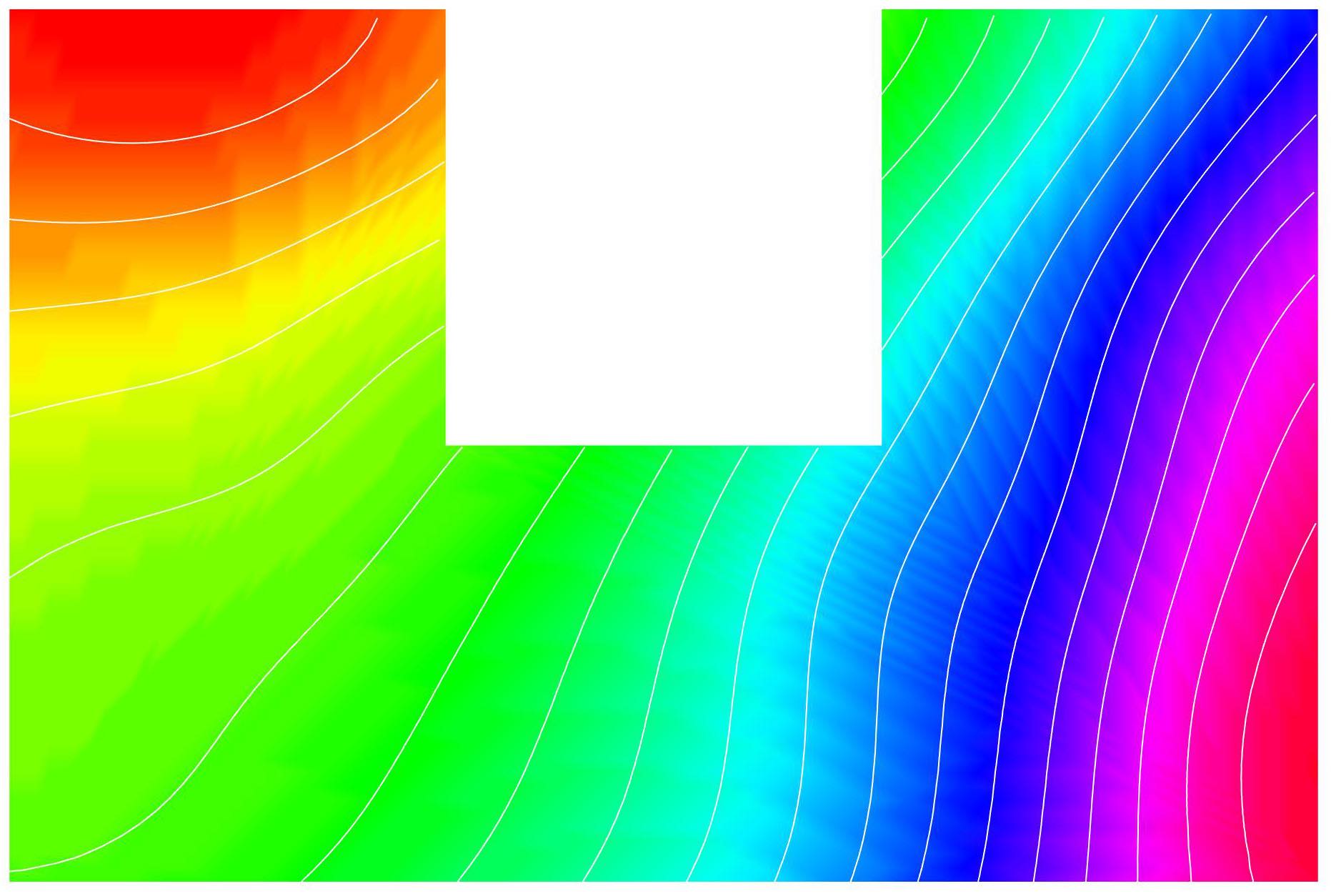}
&\includegraphics[height=60pt]{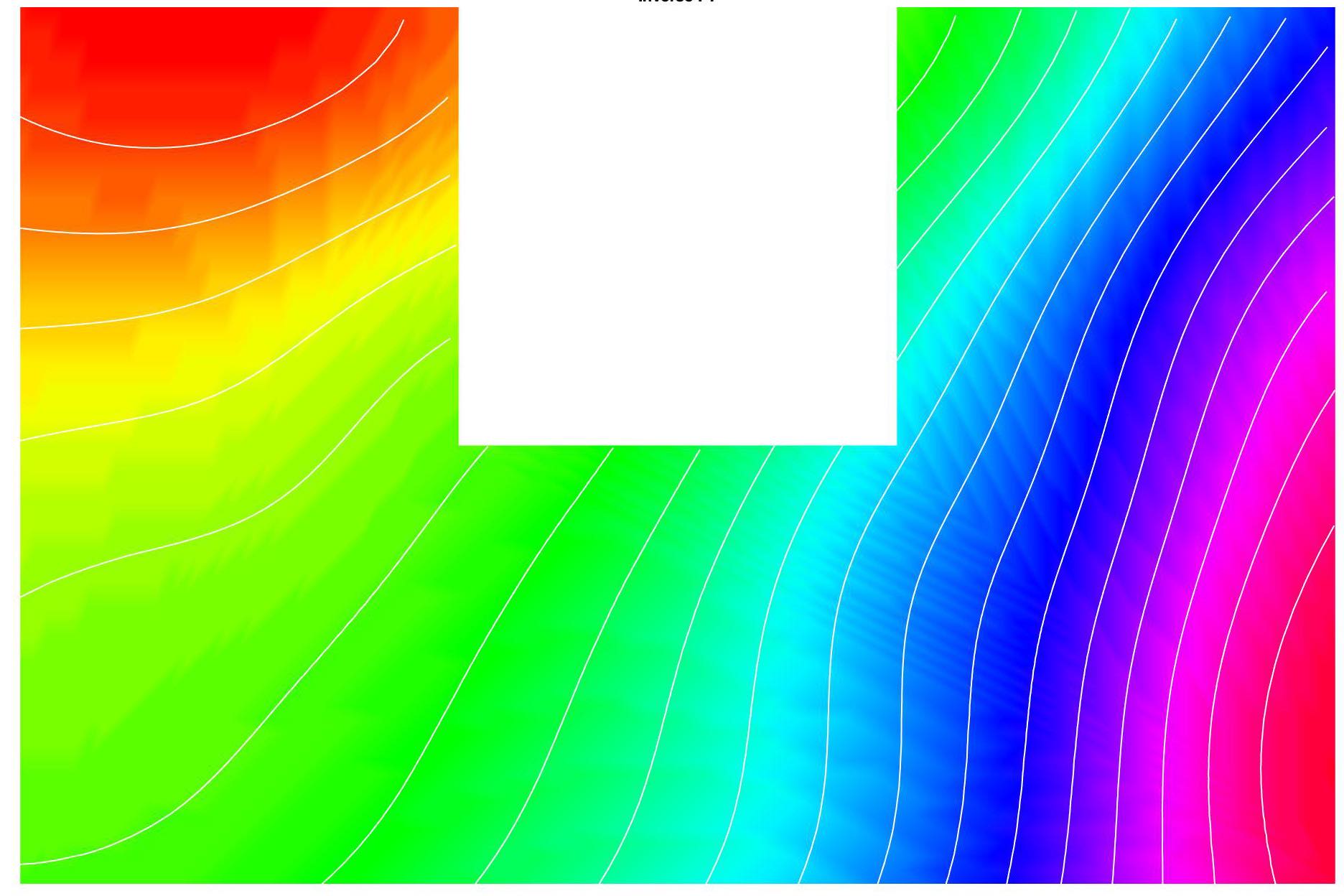}
&\includegraphics[height=60pt]{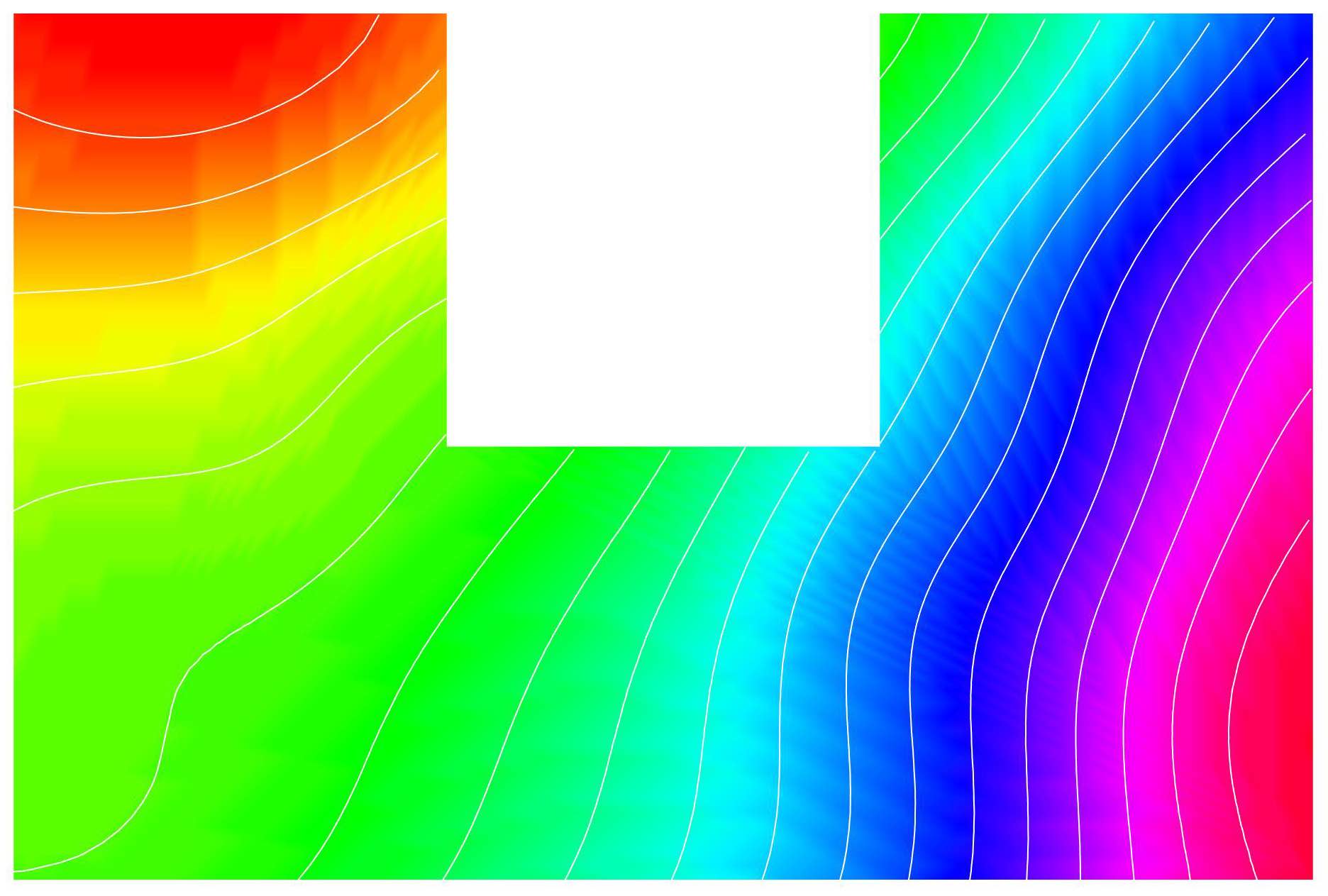}\\
$\mathcal{F}^{-1}(\mathcal{F}f)$ &$\epsilon_{\infty}:=3.5\%$		&$\epsilon_{\infty}:=6.3\%$		&$\epsilon_{\infty}:=8.2\%$		&$\epsilon_{\infty}:=9.1\%$\\
\hline
\end{tabular}
\caption{Level-sets and color-map of (first row) noisy signals with an increasing error magnitude~$\alpha$ (from (b) to (e)) achieved by adding a Gaussian noise to (a) an input signal~$f$, (second row) reconstructed signal \mbox{$\mathcal{F}^{-1}(\mathcal{F}f)$}. Here, the F-transform \mbox{$\mathcal{F}f$} is induced by the multi-quadratic kernel.\label{fig:2D-FT-IFT-KERNELS}}
\end{figure*}
\subsection{Selection and generation of membership functions\label{KERNEL-EXAMPLE}}
In the continuous and discrete settings, we consider kernels commonly used for the integral operators in machine learning~\cite{SCHOELKOPF02}, and normalised as in Eq. (\ref{eq:NORMALISED-KERNEL}). We represent the behaviour of a function \mbox{$f:\Omega\rightarrow\mathbb{R}$} through the corresponding level-sets \mbox{$\gamma_{\alpha}:=\{\mathbf{p}:\,f(\mathbf{p})=\alpha\}$} and colour-map, which begins with red, passes through yellow, green, cyan, blue, and magenta, and returns to red. In our experiments, we compute the membership functions (Figs.~\ref{fig:2D-ANALYTIC-GAUSSIAN},~\ref{fig:2D-DIFFUSION-MF},~\ref{fig:2D-HAR-BIHAR},~\ref{fig:2D-DIFFUSION},~\ref{fig:3D-DIFFUSION-MF}), the F-transform of the Dirac function~$\delta_{\mathbf{p}}$ at a seed point~$\mathbf{p}$ (Fig.~\ref{fig:3D-FT-IFT}), and the smoothing of noisy signals (Figs.~\ref{fig:2D-FT-IFT-KERNELS-NOISE},~\ref{fig:2D-FT-IFT-KERNELS},~\ref{fig:3D-FT-IFT-KERNELS-NOISE}). Then, we analyse their main properties, such as locality, encoding of local properties of the input domain, and smoothness, through the locality of the level-sets of \mbox{$\mathcal{F}\delta_{\mathbf{p}}$}, their alignment with geometric features around the seed point (e.g., the tubular features of the legs of the elephant), the smooth shape and regular distribution of the level-sets. 

\textbf{Membership functions as radial kernels}
According to~\cite{PATANE2014}, we select a radial membership function (Figs.~\ref{fig:2D-ANALYTIC-GAUSSIAN},~\ref{fig:2D-DIFFUSION-MF}) induced by a kernel \mbox{$K_{i}(\mathbf{p}):=\varphi(\|\mathbf{p}-\mathbf{p}_{i}\|_{2})$}, \mbox{$i=1,\ldots,n$}, centred at~$\mathbf{p}_{i}$ and generated by a \emph{kernel} \mbox{$\varphi:\mathbb{R}^{+}\rightarrow\mathbb{R}$}. Main examples include the Gaussian kernel \mbox{$\varphi(t):=\exp(-t/h)$}, where~$h$ is the kernel support, the triangular and sinusoidal shaped basis kernels induced by \mbox{$\varphi_{1}(t):=\frac{1-t}{h}$}, \mbox{$\varphi_{2}(t):=\frac{t}{h}$} and \mbox{$\varphi(t):=\frac{\cos t}{h}$}, \mbox{$t\in\mathbb{R}$}, respectively. Further options are the \emph{polynomial kernel}: \mbox{$K(\mathbf{p},\mathbf{q}):=\langle\mathbf{p},\mathbf{q}\rangle_{2}^{d}$}, the \emph{sigmoid kernel}: \mbox{$K(\mathbf{p},\mathbf{q}):=\tanh(\kappa \langle\mathbf{p},\mathbf{q}\rangle_{2})+\theta$}, \mbox{$\kappa>0$}, \mbox{$\theta<0$}, and the \emph{inhomogeneous polynomial}: \mbox{$K(\mathbf{p},\mathbf{q}):=(\langle\mathbf{p},\mathbf{q}\rangle_{2}+c)^{d}$}. All the previous kernels are invariant with respect to rotations, i.e., \mbox{$K(\mathbf{p},\mathbf{q})=K(\mathbf{U}^{\top}\mathbf{p}, \mathbf{U}^{\top}\mathbf{q})$}, with~$\mathbf{U}$ orthogonal matrix. 

\textbf{Tensor kernels\label{sec:TENSOR-KERNEL}}
Given two functions \mbox{$f,g:\Omega\rightarrow\mathbb{R}$}, let us introduce the \emph{tensor kernel} \mbox{$K(\mathbf{p},\mathbf{q}):=f(\mathbf{p})g(\mathbf{q})$}, whose integral operator \mbox{$(\mathcal{L}_{K}h)(\mathbf{p})=g(\mathbf{p})\langle h,f\rangle_{2}$} maps any function~$h$ to a multiple of~$g$ and its norm is \mbox{$\|\mathcal{L}_{K}\|:=\|f\|_{2}\|g\|_{2}$}. For the Gaussian function \mbox{$f(s):=\exp(-s^{2}/\sigma^{2})$}, the tensor kernel in~$\mathbb{R}^{2}$ is \mbox{$K(x,y):=(f\otimes f)(x,y)=\exp(-(x^{2}+y^{2})/\sigma^{2})$}. Finally, the \emph{tensor kernel} \mbox{$\tilde{K}(\mathbf{p},\mathbf{q})=f(\mathbf{p})K(\mathbf{p},\mathbf{q})f(\mathbf{q})$}, \mbox{$\forall f:\Omega\rightarrow\mathbb{R}$}, and the \emph{normalised kernel} \mbox{$\tilde{K}(\mathbf{p},\mathbf{q}):=\frac{K(\mathbf{p},\mathbf{q})}{\left[K(\mathbf{p},\mathbf{p})K(\mathbf{q},\mathbf{q})\right]^{1/2}}$} are positive-definite, with \mbox{$K(\mathbf{p},\mathbf{p})>0$}, \mbox{$\forall\mathbf{p}\in\Omega$}.

\textbf{Generating membership functions\label{sec:KERNEL-GENERATION}}
Given two positive-definite kernels \mbox{$K_{1}, K_{2}:\Omega\times\Omega\rightarrow\mathbb{R}$}, and the sequence of positive-definite kernels \mbox{$(K_{n})_{n}$}, the functions
\begin{equation*}\label{eq:KERNEL-DEFINITION}
\left\{
\begin{array}{ll}
\alpha K_{1}+\beta K_{2}			&\textrm{(positive linear combination)},\,\alpha,\beta\geq 0;\\
K_{1}K_{2}							&\textrm{(pointwise product)};\\
\lim_{n\rightarrow+\infty}K_{n}		&\textrm{(limit of a sequence~$(K_{n})_{n}$ of kernels)};
\end{array}
\right.
\end{equation*}
are positive-definite kernels. Indeed, the set of positive-definite kernels is closed with respect to the linear combination with positive coefficients, pointwise product, and limit (if it exists). 

\section{Data-driven continuous F-transform\label{sec:CONT-DATA-FT-SHORT}}
The relation between the F-transform and integral operators is used to introduce a \emph{data-driven F-transform} (Sect.~\ref{sec:SPECTRAL-FUZZY-OPERATOR}) through a family of \emph{data-driven membership functions} (Sect.~\ref{sec:SPECTRAL-EXAMPLES}), defined by filtering the Laplacian spectrum and encode \emph{intrinsic} information about the input data (e.g., structure, geometry, sampling density). This choice is motivated by the intrinsic definition of the Laplace-Beltrami operator, which is uniquely determined by the metrics on the input domain and encodes its geometric and topological properties. 

\subsection{Spectral data-driven F-transform\label{sec:SPECTRAL-FUZZY-OPERATOR}}
Recalling that the Laplace-Beltrami operator~$\Delta$ is self-adjoint and positive semi-definite, it has an \emph{orthonormal eigensystem} \mbox{$(\lambda_{n},\phi_{n})_{n=0}^{+\infty}$}, \mbox{$\Delta\phi_{n}=\lambda_{n}\phi_{n}$}, \mbox{$0=\lambda_{0}<\lambda_{1}\leq\lambda_{2}\leq\cdots$}, in \mbox{$\mathcal{L}^{2}(\Omega)$}. Given a strictly positive and square integrable filter function \mbox{$\varphi:\mathbb{R}^{+}\rightarrow\mathbb{R}$}, the \emph{spectral data-driven kernel}
\begin{equation}\label{eq:FUNCT-OPER}
K_{\varphi}(\mathbf{p},\mathbf{q})=\sum_{n=0}^{+\infty}\varphi(\lambda_{n})\phi_{n}(\mathbf{p})\phi_{n}(\mathbf{q})
\end{equation}
is well-posed according to the relation \mbox{$\|K_{\varphi}\|_{2}=\|\varphi\|_{2}$}, as a consequence of the orthonormality of the Laplacian eigenfunctions. Then, the \emph{spectral data-driven F-transform} \mbox{$\mathcal{L}_{\varphi}:=\mathcal{L}_{K_{\varphi}}$} is defined as the integral operator induced by~$K_{\varphi}$, i.e., \mbox{$\mathcal{L}_{\varphi}f=\langle K_{\varphi},f\rangle_{2}=\sum_{n=0}^{+\infty}\varphi(\lambda_{n})\langle f,\phi_{n}\rangle_{2}\phi_{n}$}.
From the identity \mbox{$\mathcal{L}_{\varphi_{1}}\circ \mathcal{L}_{\varphi_{2}}= \mathcal{L}_{\varphi_{1}\varphi_{2}}$},~$\mathcal{L}_{\varphi}$ is invertibile if and only if~$\varphi$ is not null; in this case, \mbox{$\mathcal{L}_{\varphi}^{-1}=\mathcal{L}_{1/\varphi}$}. Indeed, the inverse of the continuous F-transform associated with the spectral kernel~$K_{\varphi}$ is the integral operator induced by~$K_{1/\varphi}$.

\subsection{Data-driven membership functions\label{sec:SPECTRAL-EXAMPLES}}
Main examples of data-driven kernels (\ref{eq:FUNCT-OPER}) are the \emph{commute-time kernel} \mbox{$K_{\Delta}(\mathbf{p},\mathbf{q})=\sum_{n=1}^{+\infty}\lambda_{n}^{-1}\phi_{n}(\mathbf{p})\phi_{n}(\mathbf{q})$}, induced by the filter \mbox{$\varphi(s):=s^{-1}$}, and the \emph{bi-harmonic kernel} \mbox{$K_{\Delta^{2}}(\mathbf{p},\mathbf{q}):=\sum_{n=1}^{+\infty}\lambda_{n}^{-2}\phi_{n}(\mathbf{p})\phi_{n}(\mathbf{q})$}, induced by the filter \mbox{$\varphi(s):=s^{-2}$}. The commute-time and bi-harmonic membership functions are globally-supported (Fig.~\ref{fig:2D-HAR-BIHAR}). 

The \emph{diffusion kernel} \mbox{$K_{t}(\mathbf{p},\mathbf{q})=\sum_{n=0}^{+\infty}\exp(-t\lambda_{n})\phi_{n}(\mathbf{p})\phi_{n}(\mathbf{q})$} is associated with the filter \mbox{$\varphi(s):=\exp(-st)$}. Increasing or reducing the time scale~$t$ of the diffusion membership functions, we easily enlarge or reduce their support. In fact, as~$t$ becomes smaller the support of the corresponding diffusion function centred at a seed point reduces until it degenerates to the seed itself (Fig.~\ref{fig:2D-DIFFUSION}). In this case, we avoid the Gibbs phenomenon (i.e., small undulations as we move far from the seed point) through the Pad\`e-Chebyshev approximation of the heat kernel~\cite{PATANE-STAR2016}.

The selection of data-driven membership functions and integral kernels for the definition of the F-transform allows us to efficiently encode local and global information about the input data in a \emph{multi-scale} manner, thus expressing complex dependencies among variables for large data sets.  Data-driven membership functions avoid coarse fuzzy partitions, which group data into large clusters that do not adapt to their local behaviour, or a too dense fuzzy partition, which generally has cells that are not covered by the data, thus being redundant and resulting in a higher computational cost. Indeed, data-driven membership functions provide an alternative to previous work~\cite{DIMARTINO2011}, which typically controls that the resolution of the fuzzy partition is not too dense with respect to the data sampling through the Wang-Mendel method~\cite{WANG1992}. 

The aforementioned properties of the membership functions are then inherited by the continuous F-transform and are important in case of \emph{structured} (e.g., regular/irregular) \emph{data} and \emph{sparse} or \emph{dense data}. Finally, the spectral representation is applied to discretise the continuous F-transform and its inverse in terms of the singular value decomposition of the Gram matrix associated with the input kernel (Sect.~\ref{sec:DISCRETISATION}). 

\section{Spectral continuous F-transform and inverse\label{sec:SPECTRAL-REPRESENTATION-FT}}
The properties of the continuous F-transform as integral operator are used to represent the F-transform (Sect.~\ref{sec:SPECTRAL-CFT}) and its pseudo-inverse (Sect.~\ref{sec:PSEUDO-INVERSE}) in terms of the spectrum of the integral operator. Then, we specialise these results to Reproducing Kernel Hilbert Spaces (RHKS) (Sect.~\ref{sec:CONT-FT-RKHS}).
\begin{figure}[t]
\centering
\begin{tabular}{cccc}
$t=10^{-3}$		&$t=10^{-2}$ 	&$t=10^{-1}$	&$t=1$\\
\includegraphics[height=60pt]{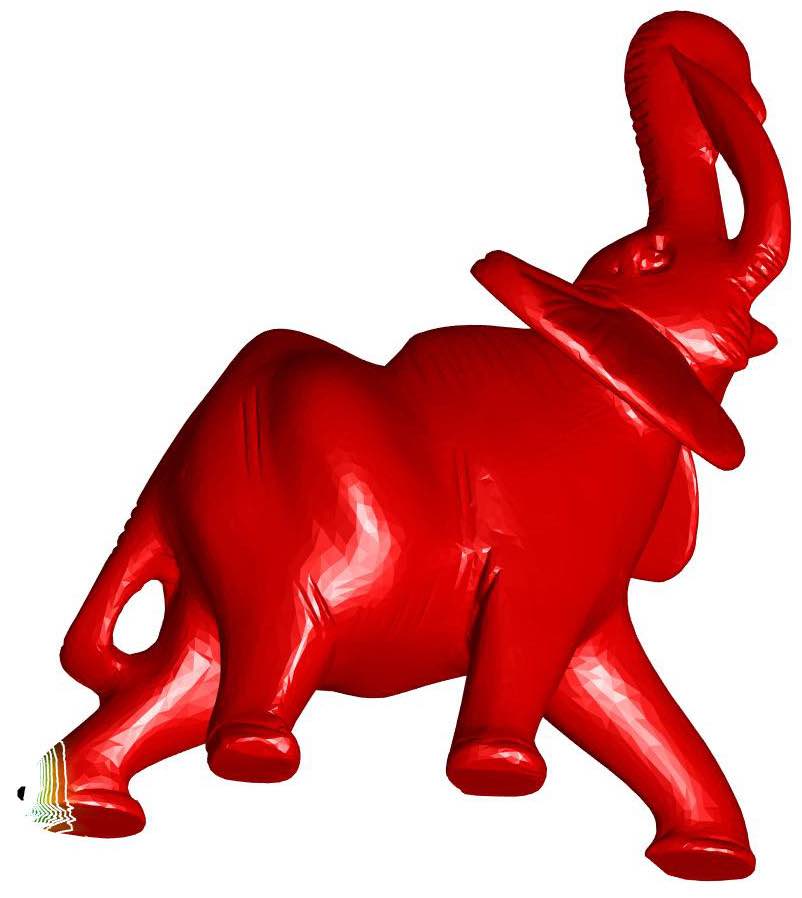}
&\includegraphics[height=60pt]{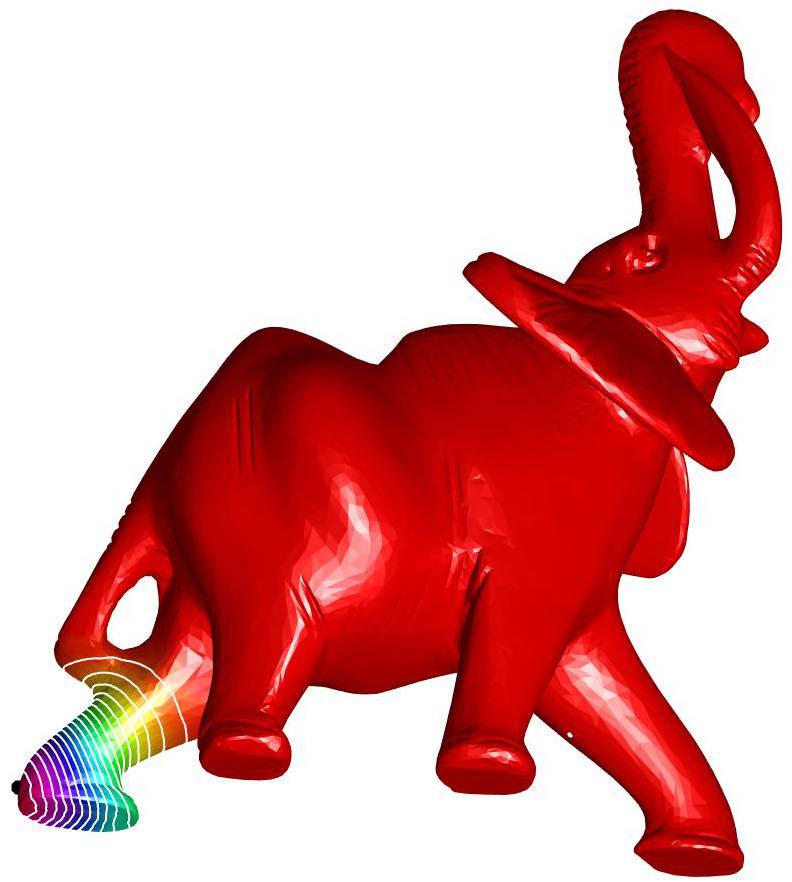}
&\includegraphics[height=60pt]{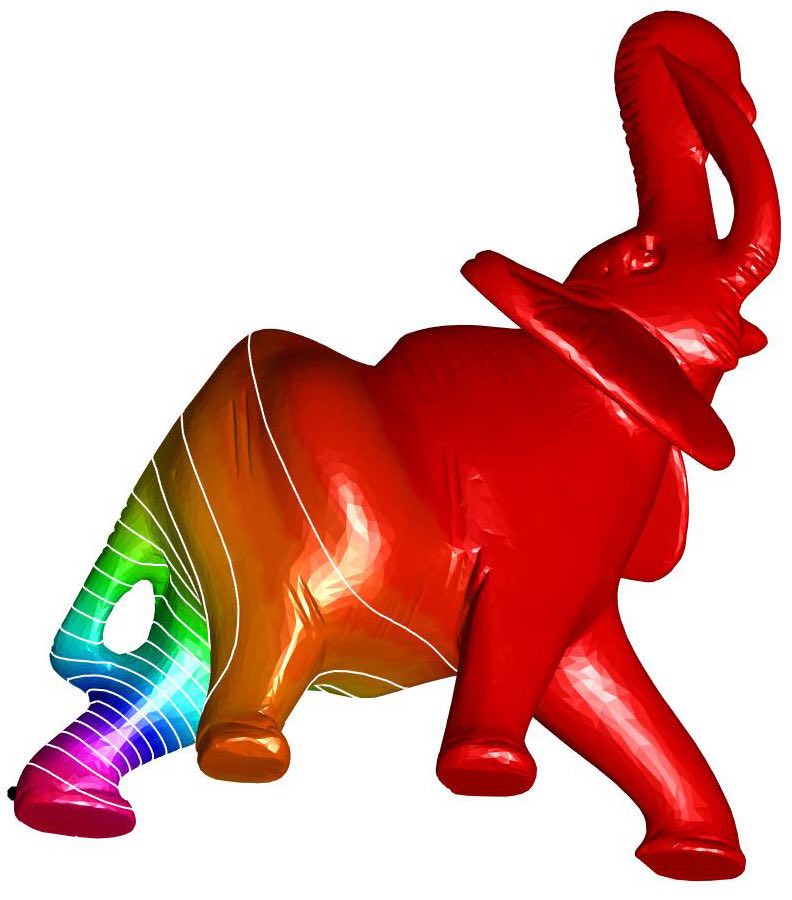}
&\includegraphics[height=60pt]{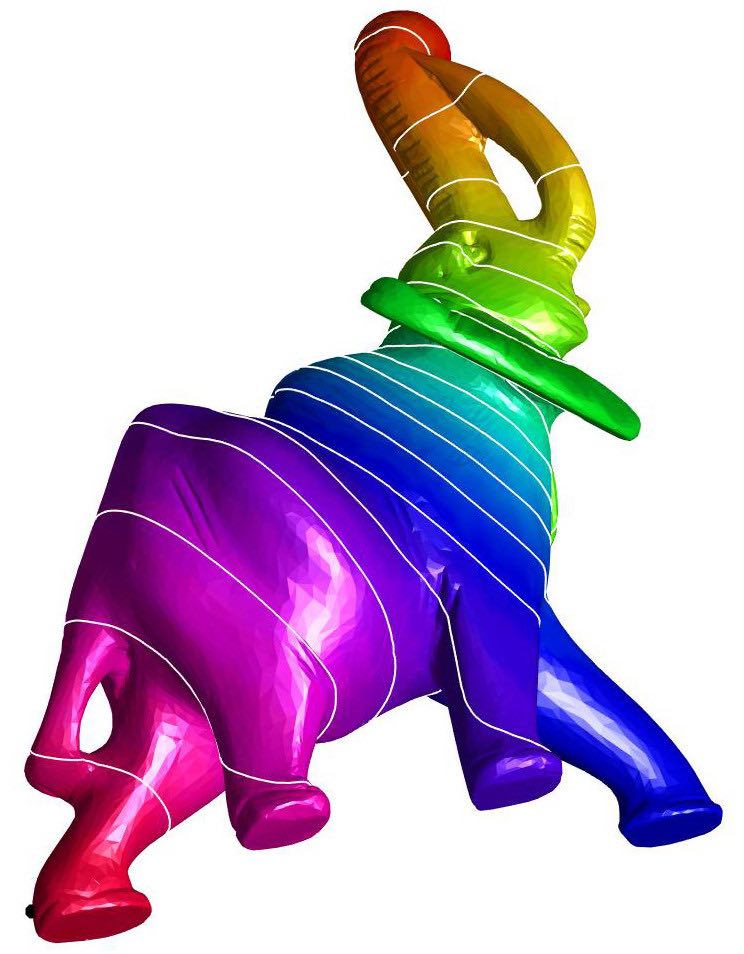}
\end{tabular}
\caption{Level-sets of diffusive membership functions centred at a seed point and at increasing scales~$t$.\label{fig:3D-DIFFUSION-MF}}
\end{figure}
\subsection{Spectral continuous F-transform in Hilbert Spaces\label{sec:SPECTRAL-CFT}}
Since the continuous F-transform~$\mathcal{L}_{K}$ is a linear, continuous, self-adjoint operator (Sect.~\ref{sec:CONTINUOUS-FT-PROPERTIES}), it admits a real eigensystem \mbox{$(\lambda_{n},\phi_{n})_{n=0}^{+\infty}$} in \mbox{$\mathcal{L}^{2}(\Omega)$} such that \mbox{$\mathcal{L}_{K}\phi_{n}=\lambda_{n}\phi_{n}$}, with \mbox{$\langle\phi_{n},\phi_{m}\rangle_{2}=\delta_{mn}$} and \mbox{$\lambda_{n}\leq\lambda_{n+1}$}. In particular,~$\phi_{n}$ is continuous; in fact, \mbox{$\phi_{n}=\lambda_{n}^{-1}\mathcal{L}_{K}\phi_{n}$}, \mbox{$\lambda_{n}\neq 0$}.

According to the Mercer theorem~\cite{SCHOELKOPF02}, we represent the kernel in terms of the spectrum of the integral operator as \mbox{$K(\mathbf{p},\mathbf{q})
=\sum_{n=0}^{+\infty}\lambda_{n}\phi_{n}(\mathbf{p})\phi_{n}(\mathbf{q})$}, and the integral operator is rewritten as \mbox{$F=\mathcal{L}_{K}f=\sum_{n=0}^{+\infty}\lambda_{n}\langle f,\phi_{n}\rangle_{2}\phi_{n}$}. From the previous relations, we get that
\begin{equation}\label{eq:ENERGY-INTEGRAL-ACTION}
\|\mathcal{L}_{K}f\|_{2}^{2}
=\sum_{n=0}^{+\infty}\lambda_{n}^{2}\vert\langle f,\phi_{n}\rangle_{2}\vert^{2}
\leq\|f\|_{2}\sum_{n=0}^{+\infty}\lambda_{n}^{2};
\end{equation}
indeed, the energy \mbox{$\|\mathcal{L}_{K}\|_{2}^{2}=\sum_{n=0}^{+\infty}\lambda_{n}^{2}$} of the continuous F-transform is equal to the~$\ell_{2}$-norm of the eigenvalues. In fact, one upper bound is given by the inequality (\ref{eq:ENERGY-INTEGRAL-ACTION}) and selecting \mbox{$f=\sum_{n=0}^{+\infty}\phi_{n}$}, we get the opposite inequality. 
\begin{figure*}[t]
\centering
\begin{tabular}{cc|cc|cc}
\hline
\multicolumn{2}{c|}{\textbf{Gaussian memb. funct.}}
&\multicolumn{2}{c|}{\textbf{Hyperbolic Tangent memb. funct.}}
&\multicolumn{2}{c}{\textbf{Polynomial memb. funct.}}\\
\hline
FT		&Inverse FT			&FT		&Inverse FT			&FT		&Inverse FT\\	
\hline
\includegraphics[height=80pt]{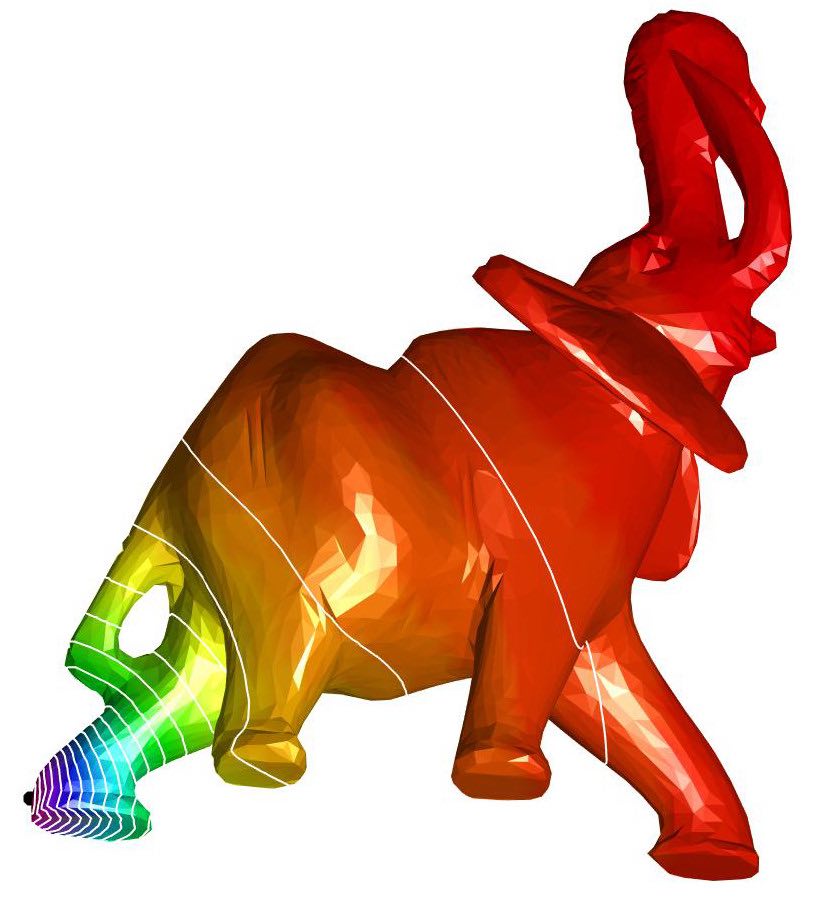}
&\includegraphics[height=80pt]{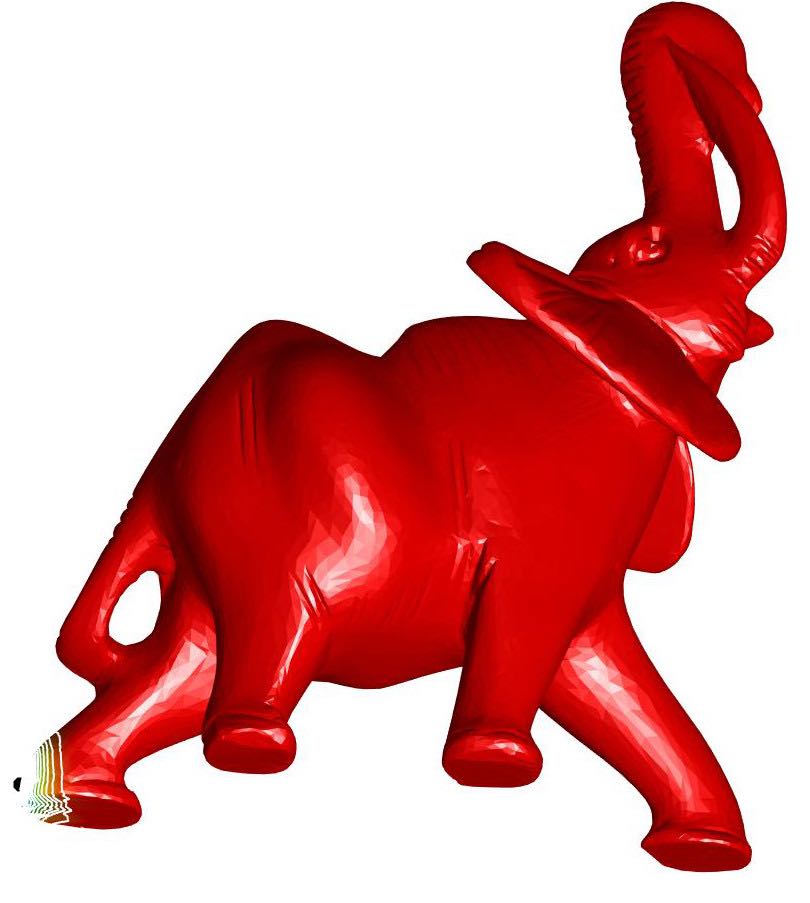}
&\includegraphics[height=80pt]{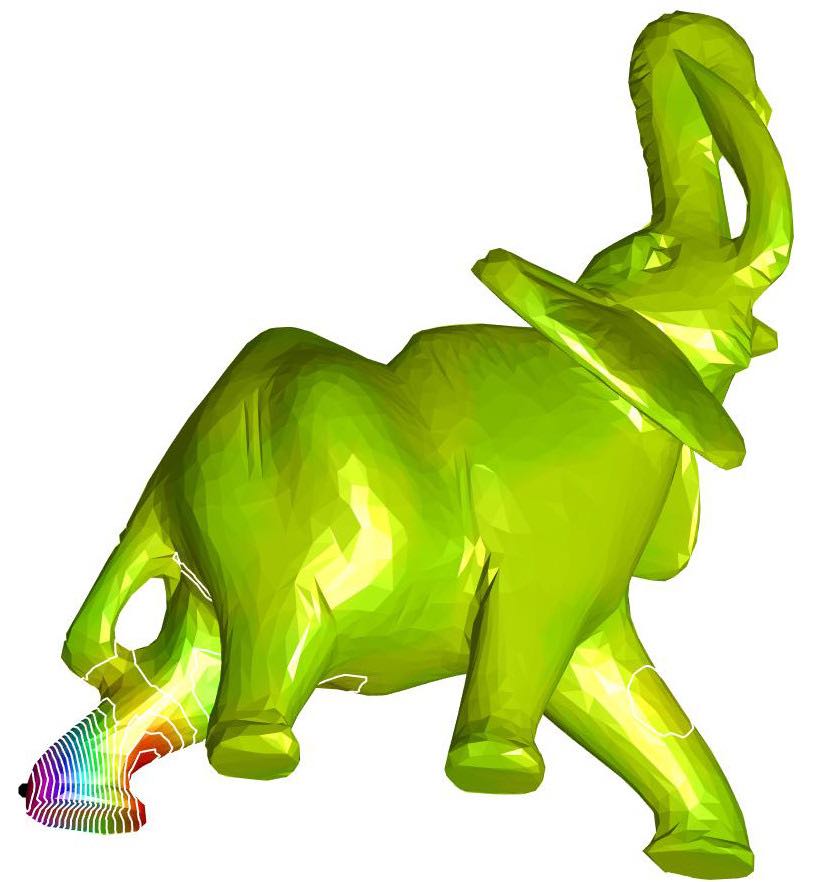}
&\includegraphics[height=80pt]{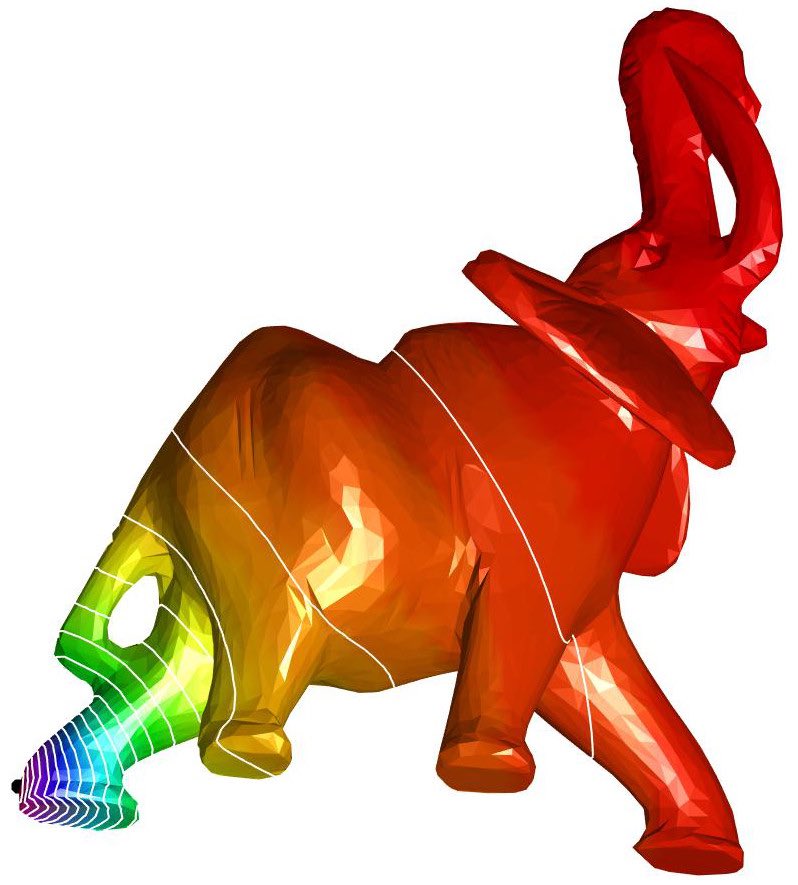}
&\includegraphics[height=80pt]{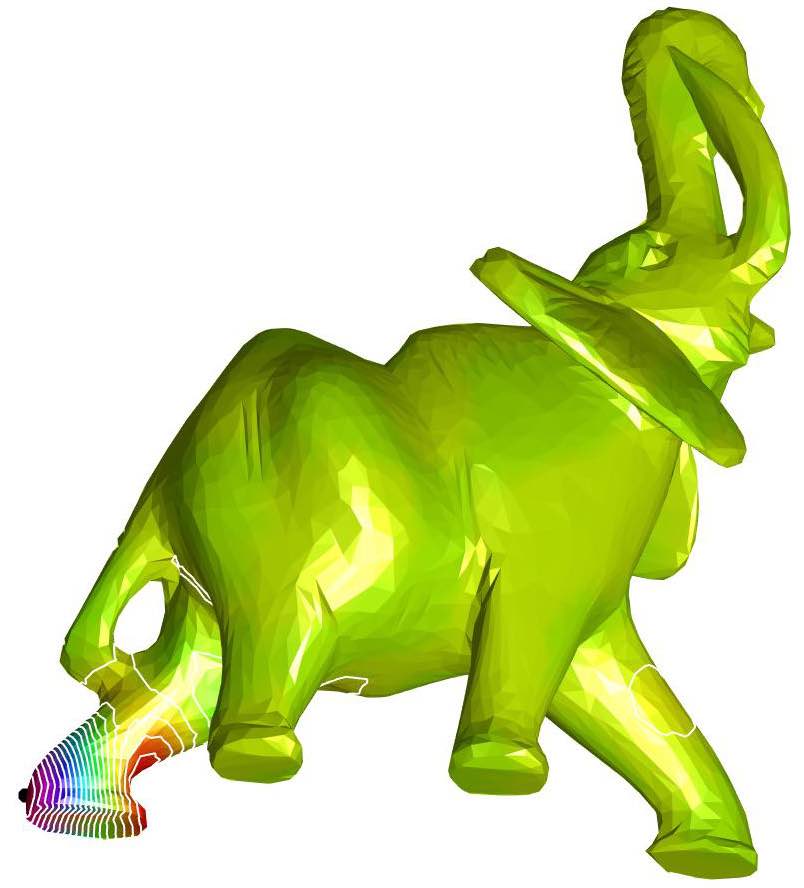}
&\includegraphics[height=80pt]{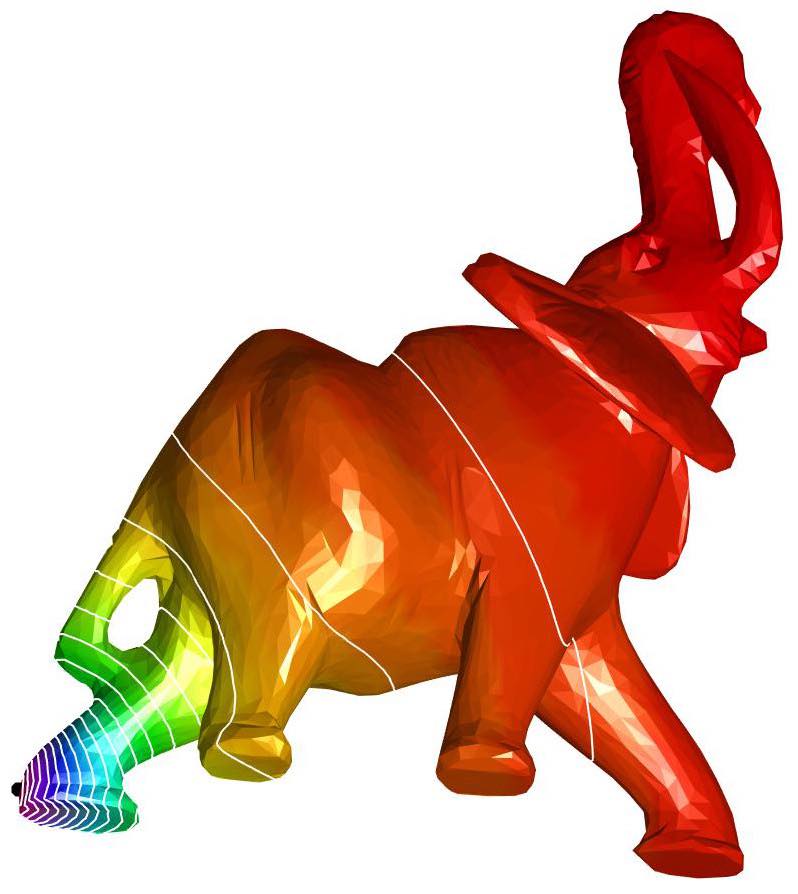}\\
\end{tabular}
\caption{Level-sets of the continuous (FT) and inverse (inverse FT) F-transform of the Dirac~$\delta$-function at a seed point (localised on the elephant feet) induced by the Gaussian, hyperbolic tangent, and polynomial membership functions. The F-transform induced by the hyperbolic tangent and polynomial membership functions is more localised around the seed point than the F-transform induced by the Gaussian membership function, while the inverse F-transform shows an opposite behaviour.\label{fig:3D-FT-IFT}}
\end{figure*}
\subsection{Pseudo-inverse of the F-transform in Hilbert Spaces\label{sec:PSEUDO-INVERSE}}
Given \mbox{$g\in\mathcal{C}^{0}(\Omega)$}, we compute the function \mbox{$f\in\mathcal{L}^{2}(\Omega)$} such that \mbox{$\mathcal{L}_{K}f=g$}. Expressing the functions \mbox{$f=\sum_{n=0}^{+\infty}a_{n}\phi_{n}$} and \mbox{$g=\sum_{n=0}^{+\infty}\langle g,\phi_{n}\rangle_{2}\phi_{n}$} in terms of the eigensystem of~$\mathcal{L}_{K}$ and imposing the previous condition, we get the relation
\begin{equation*}
\mathcal{L}_{K}f=g,\quad
\sum_{n=0}^{+\infty}a_{n}\lambda_{n}\phi_{n}
=\sum_{n=0}^{+\infty}\langle g,\phi_{n}\rangle_{2}\phi_{n};
\end{equation*}
i.e., \mbox{$a_{n}=\frac{\langle g,\phi_{n}\rangle_{2}}{\lambda_{n}}$}. The representation of~$f$ in terms of the spectrum of the integral operator belongs to \mbox{$\mathcal{L}^{2}(\Omega)$} (i.e.,~$f$ is well-defined) if and only if
\begin{equation*}
\sum_{n=0}^{+\infty}\frac{\vert\langle g,\phi_{n}\rangle_{2}\vert^{2}}{\lambda_{n}^{2}}<+\infty;\textrm{ i.e., }
\left(\frac{\vert\langle g,\phi_{n}\rangle_{2}\vert}{\lambda_{n}}\right)_{n=0}^{+\infty}\in\ell_{2}.
\end{equation*}
Recalling that \mbox{$\lim_{n\rightarrow+\infty}\lambda_{n}=0$}, this last condition is satisfied if \mbox{$(\lambda_{n}^{-1})_{n=0}^{+\infty}\in\ell_{2}$}. Then, the \emph{spectral representation of the inverse continuous F-transform} is
\begin{equation}\label{eq:SPECTRAL-PSEUDO-INVERSE}
\mathcal{L}_{K}^{\dag}g
=\sum_{n=0}^{+\infty}\frac{\langle g,\phi_{n}\rangle_{2}}{\lambda_{n}}\phi_{n},
\end{equation}
where \mbox{$\mathcal{L}_{K}^{\dag}:\mathcal{C}^{0}(\Omega)\rightarrow\mathcal{C}^{0}(\Omega)$} is the \emph{pseudo-inverse} of~$\mathcal{L}_{K}$. Given a function~$g$, the \emph{best approximation} \mbox{$\mathcal{L}_{K}f$} of~$g$ with respect to the~$\mathcal{L}^{2}$-norm is achieved for \mbox{$f:=\mathcal{L}_{K}^{\dag}g$}, i.e., \mbox{$\mathcal{L}_{K}^{\dag}g
=\arg\min_{f}\{\|\mathcal{L}_{K}f-g\|_{2}\}$}. If~$\mathcal{L}_{K}$ has a finite number of non-null eigenvalues, then the continuous inverse F-transform~$\mathcal{L}_{K}^{\dag}$ is equal to the integral operator~$\mathcal{L}_{K^{\dag}}$ induced by the \emph{pseudo-inverse kernel} \mbox{$K^{\dag}(\mathbf{p},\mathbf{q})=\sum_{n=0}^{+\infty}\frac{1}{\lambda_{n}}\phi_{n}(\mathbf{p})\phi_{n}(\mathbf{q})$}, and Eq. (\ref{eq:SPECTRAL-PSEUDO-INVERSE}) is rewritten as \mbox{$\mathcal{L}_{K}^{\dag}g=\mathcal{L}_{K^{\dag}}g$}.

\subsection{Continuous F-transform in RHKS\label{sec:CONT-FT-RKHS}}
We derive the conditions for the existence of the inverse continuous F-transform and its spectral representation. Given a compact set~$\Omega$ and a Mercel kernel \mbox{$K:\Omega\times\Omega\rightarrow\mathbb{R}$}, there exists a unique Hilbert Space~$\mathcal{H}_{K}$ of functions on~$\Omega$, endowed with the scalar product \mbox{$\langle f,g\rangle_{K}:=\sum_{i,j}\alpha_{i}\beta_{j}K(\mathbf{p}_{i},\mathbf{p}_{j})$}, \mbox{$f=\sum_{i}\alpha_{i}K(\mathbf{p}_{i},\cdot)$}, \mbox{$g=\sum_{i}\beta_{i}K(\mathbf{p}_{i},\cdot)$}, such that (i) \mbox{$K_{\mathbf{p}}:=K(\mathbf{p},\cdot)$} belongs to~$\mathcal{H}_{K}$, for any \mbox{$\mathbf{p}\in\Omega$}, (ii) the span of the set \mbox{$\{K(\mathbf{p},\cdot)\}_{\mathbf{p}\in\Omega}$} is dense in~$\mathcal{H}_{K}$, and (iii) \mbox{$f(\mathbf{p})=\langle f,K(\mathbf{p},\cdot)\rangle_{K}$} for any~$\mathbf{p}$ in~$\Omega$.

The elements of~$\mathcal{H}_{K}$ are continuous functions on~$\Omega$, and any~$f$ in~$\mathcal{H}_{K}$ is represented as \mbox{$f=\sum_{n=0}^{+\infty}a_{n}\phi_{n}$}, where the series converges uniformly and absolutely. In this case, we define the continuous F-transform as \mbox{$\mathcal{L}_{K}f:=\langle K(\mathbf{p},\cdot),f\rangle_{K}=f(\mathbf{p})$}, i.e., the value \mbox{$f(\mathbf{p})$} is recovered through the continuous F-transform~$\mathcal{L}_{K}$.

The \emph{Reproducing Kernel Hilbert Space}~$\mathcal{H}_{K}$ (RHKS) is defined in terms of the spectrum of the integral operator~$\mathcal{L}_{K}$ as
\begin{equation*}
\mathcal{H}_{K}:=\left\{f\in\mathcal{L}^{2}(\Omega):\quad f=\sum_{n=0}^{+\infty}a_{n}\phi_{n},\quad \left(\frac{a_{n}}{\lambda_{n}^{1/2}}\right)_{n}\in\ell_{2}\right\}.
\end{equation*}
Considering \mbox{$f=\sum_{n=0}^{+\infty}a_{n}\phi_{n}$} and \mbox{$g=\sum_{n=0}^{+\infty}b_{n}\phi_{n}$} in~$\mathcal{H}_{K}$, \mbox{$\langle f,g\rangle_{K}=\sum_{n=0}^{+\infty}\frac{a_{n}b_{n}}{\lambda_{n}}$} is the spectral representation of the scalar product that makes~$\mathcal{H}_{K}$ a Hilbert Space. In particular, \mbox{$\|f\|_{K}^{2}=\sum_{n=0}^{+\infty}\frac{\vert a_{n}\vert^{2}}{\lambda_{n}}$} and \mbox{$\langle\phi_{m},\phi_{n}\rangle_{K}=\frac{1}{\lambda_{n}}\delta_{mn}$}. 

\textbf{Spectral representation of the pseudo-inverse in RKHS}
Recalling that \mbox{$\mathcal{B}:=(\lambda_{n}^{1/2}\phi_{n})_{n=0}^{+\infty}$} is orthonormal with respect to the \mbox{$\langle\cdot,\cdot\rangle_{K}$} scalar product, we express the functions
\begin{equation}\label{eq:RKHS-SPECTRAL}
\left\{
\begin{array}{l}
f=\sum_{n=0}^{+\infty}a_{n}\phi_{n},\quad a_{n}:=\lambda_{n}\langle f,\phi_{n}\rangle_{K},\\
g=\sum_{n=0}^{+\infty}b_{n}\phi_{n},\quad b_{n}:=\lambda_{n}\langle g,\phi_{n}\rangle_{K},
\end{array}
\right.
\end{equation}
in terms of~$\mathcal{B}$. Imposing that \mbox{$\mathcal{L}_{K}f=g$}, we get that
\begin{equation*}
f(\mathbf{p})
=\sum_{n=0}^{+\infty}\langle g,\phi_{n}\rangle_{K}\phi_{n}(\mathbf{p})
=\langle H(\mathbf{p},\cdot),g\rangle_{K}=:\mathcal{L}_{H}^{\dag}g,
\end{equation*}
with \mbox{$H(\mathbf{p},\mathbf{q}):=\sum_{n=0}^{+\infty}\phi_{n}(\mathbf{p})\phi_{n}(\mathbf{q})$} \emph{spectral kernel} in~$\mathcal{H}_{K}$. Indeed, we introduce the continuous (pseudo) inverse F-transform as the integral operator \mbox{$\mathcal{L}_{H}^{\dag}g(\mathbf{p}):=\langle H(\mathbf{p},\cdot),g)\rangle_{K}$}. Noting that \mbox{$\langle f,\phi_{n}\rangle_{K}=\lambda_{n}^{1/2}\langle f,\phi_{n}\rangle_{2}$}, Eq. (\ref{eq:RKHS-SPECTRAL}) is rewritten as
\begin{equation*}
\mathcal{L}_{H}^{\dag}g
=\sum_{n=0}^{+\infty}\frac{\langle g,\phi_{n}\rangle_{2}}{\lambda_{n}^{1/2}}
=\langle \tilde{H}(\mathbf{p},\cdot),g\rangle_{2},\,
\tilde{H}(\mathbf{p},\cdot)\in\mathcal{L}^{2}(\Omega\times\Omega),
\end{equation*}
with \mbox{$\tilde{H}(\mathbf{p},\mathbf{q}):=\sum_{n=0}^{+\infty}\lambda_{n}^{-1/2}\phi_{n}(\mathbf{p})\phi_{n}(\mathbf{q})$} \emph{spectral kernel}. 

\section{Discussion\label{sec:DISCUSSION}}
We now introduce the discretisation of the continuous F-transform and its inverse (Sect.~\ref{sec:DISCRETISATION}); then, we discuss experimental results on 2D and 3D data (Sect.~\ref{sec:EXPERIMENTS}).
\begin{figure*}[t]
\centering
\begin{tabular}{c|ccccc}
(a)~$\alpha:=0$ &(b)~$\alpha:=1\%$ &(c)~$\alpha:=5\%$ &(d)~$\alpha:=10\%$ &(e)~$\alpha:=50\%$ &(f)~$\alpha:=100\%$\\
\hline
$f$\includegraphics[height=90pt]{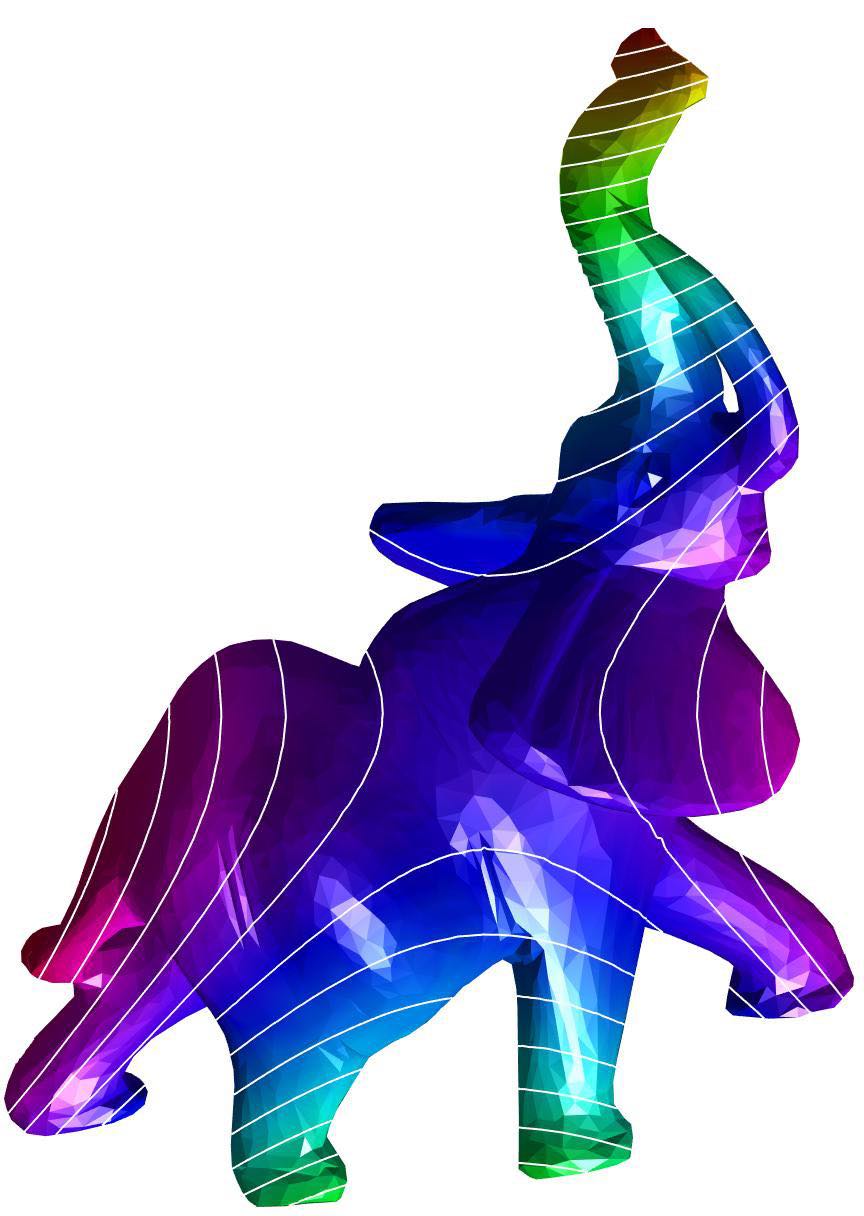}
&\includegraphics[height=90pt]{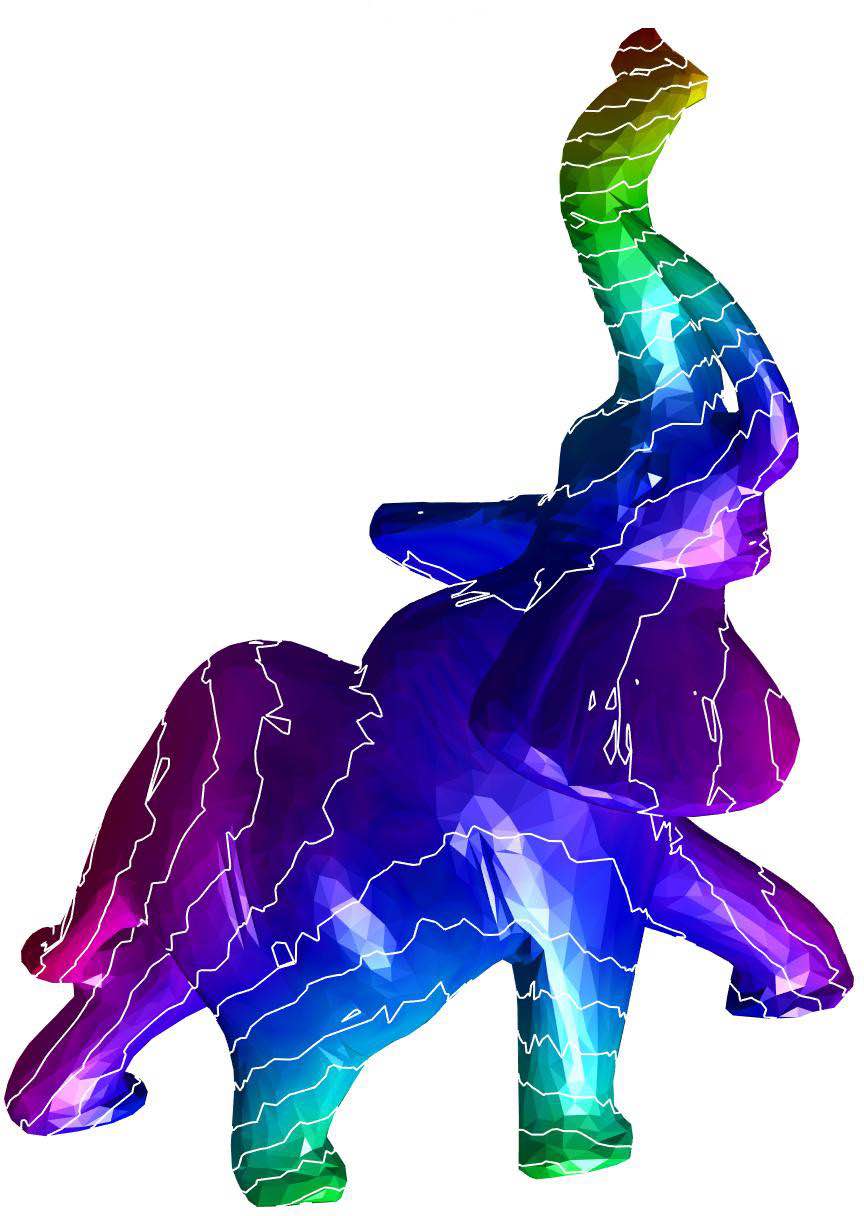}
&\includegraphics[height=90pt]{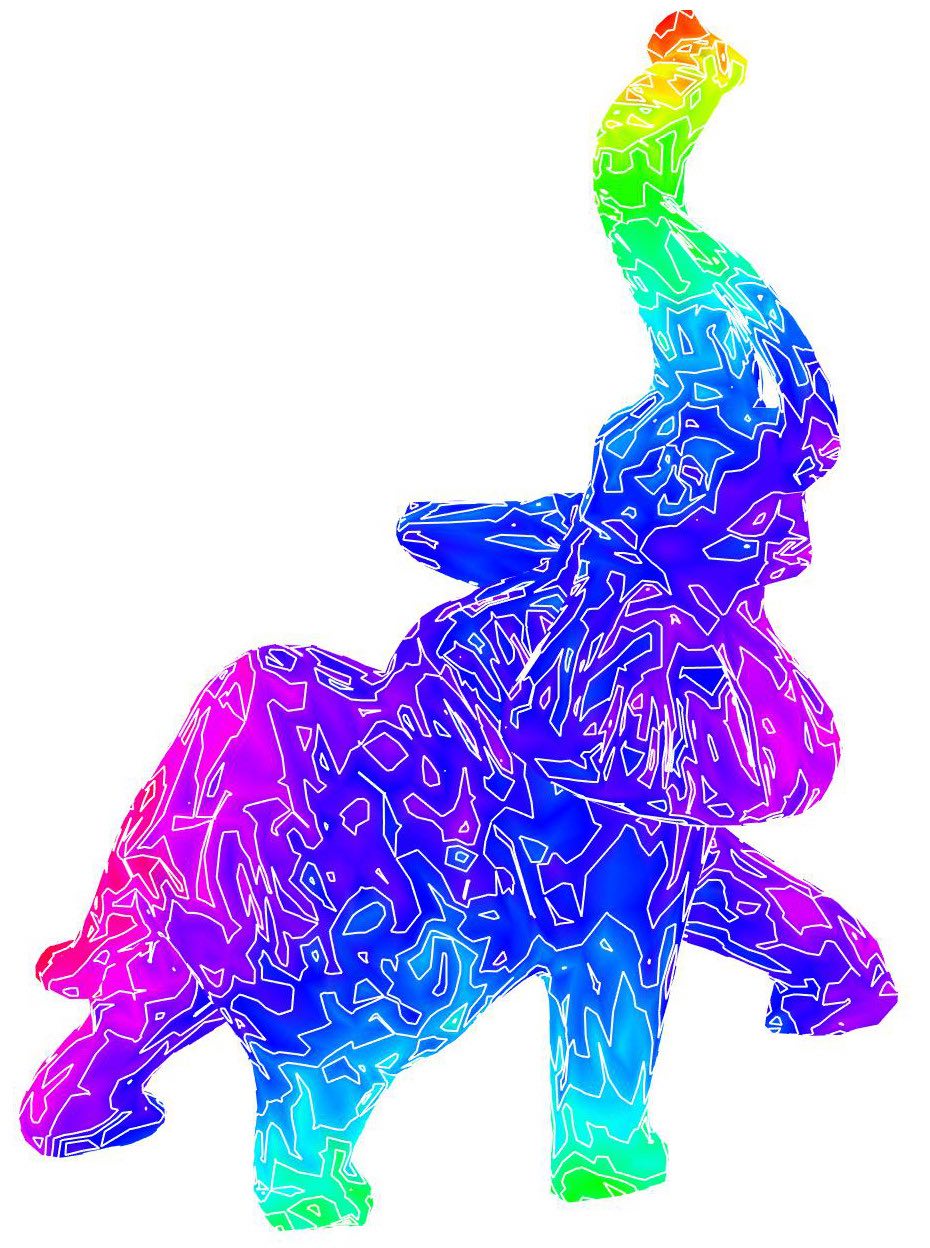}
&\includegraphics[height=90pt]{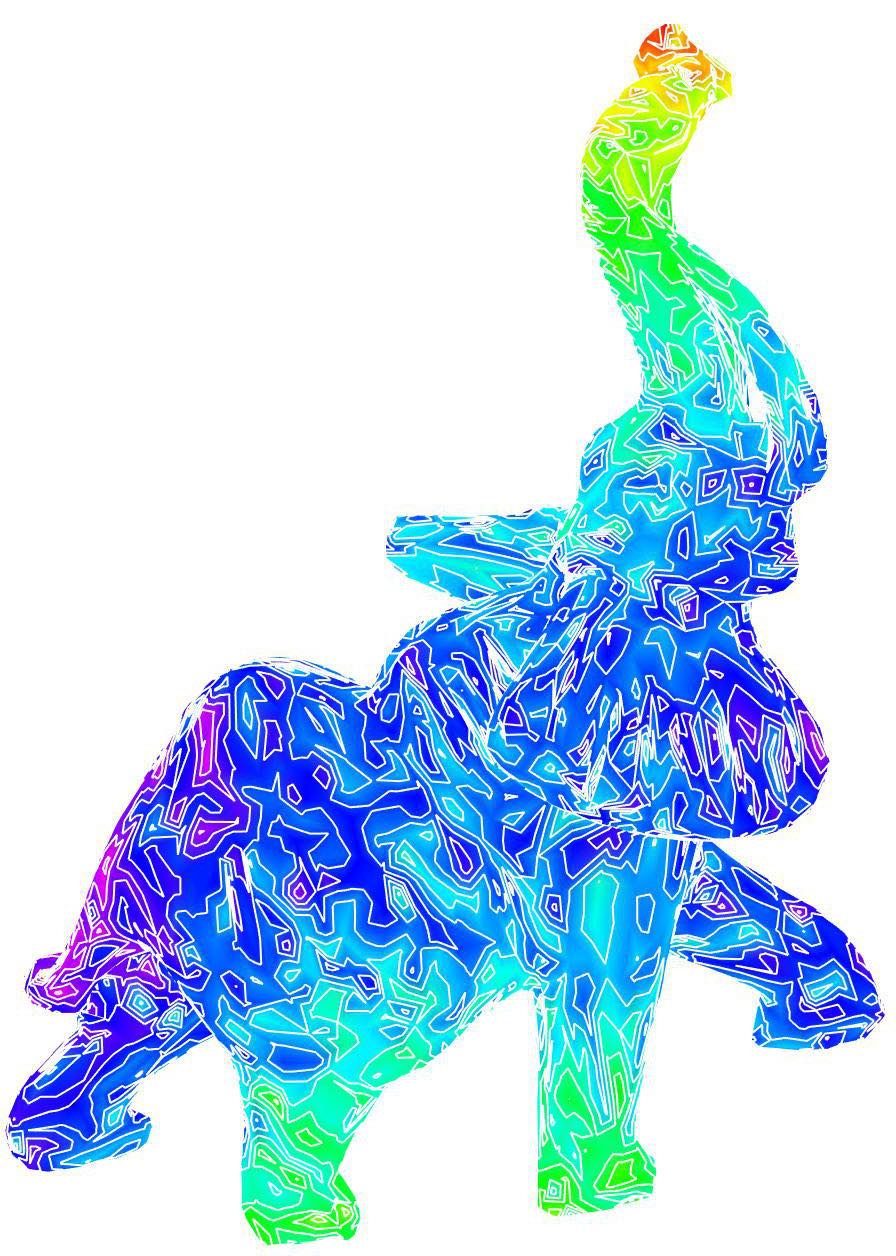}
&\includegraphics[height=90pt]{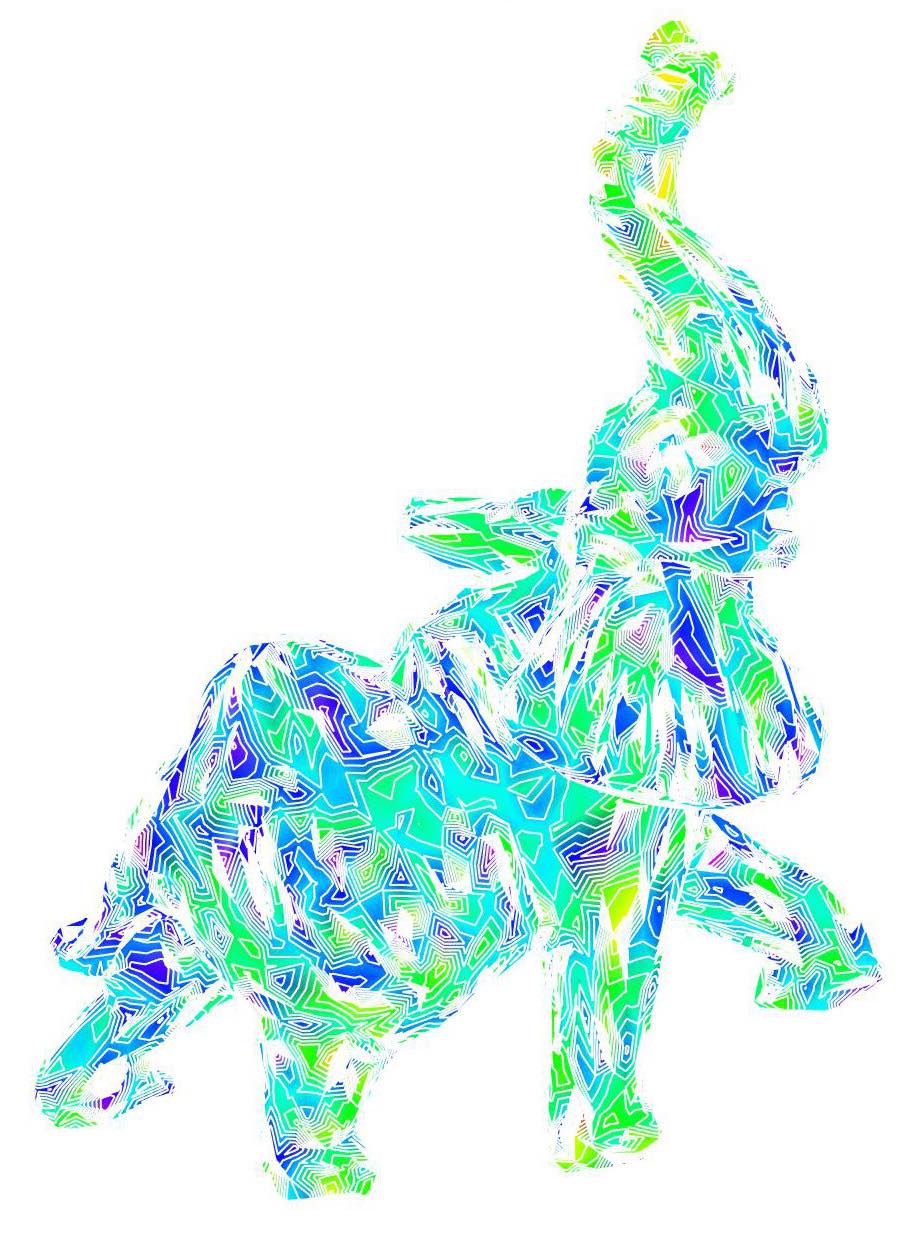}
&\includegraphics[height=90pt]{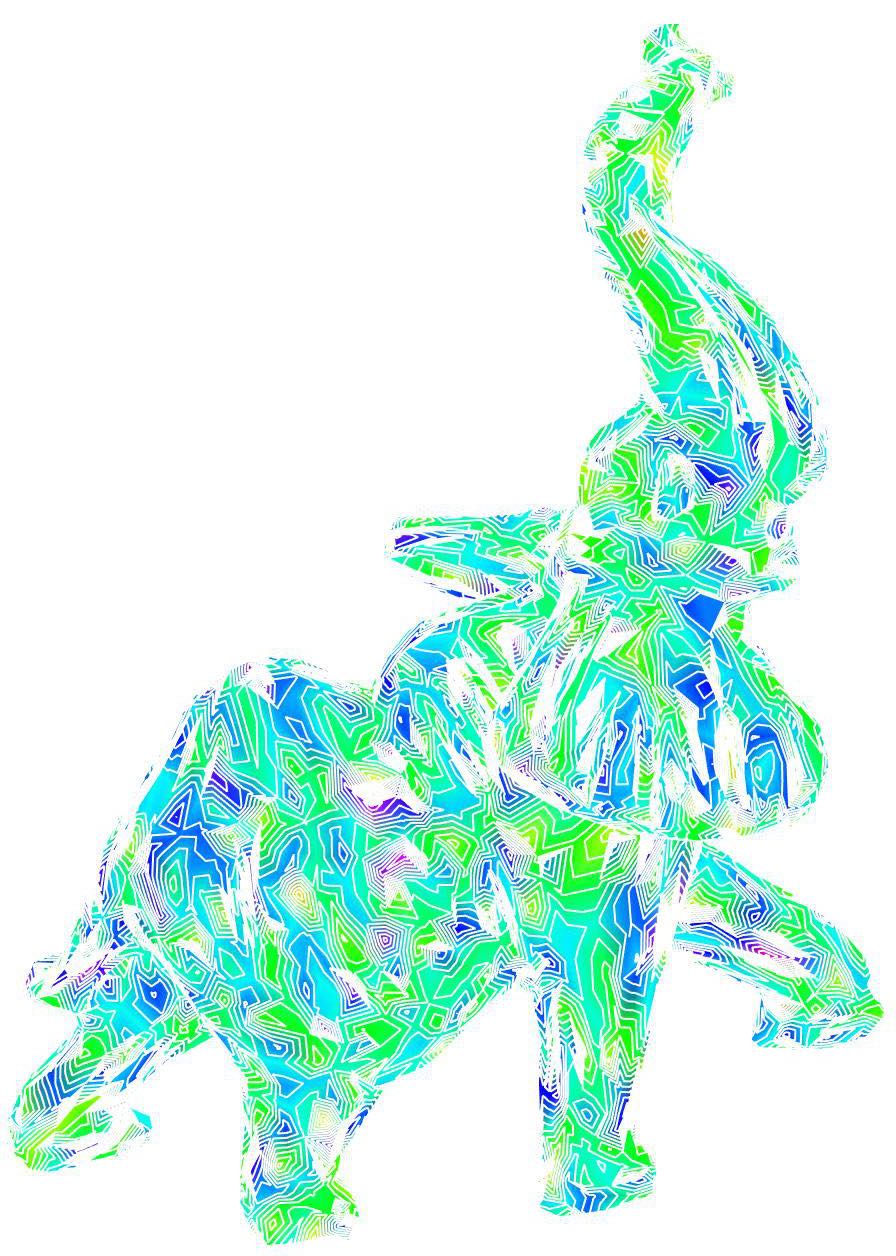}\\
\hline
$\mathcal{F}f$&\includegraphics[height=90pt]{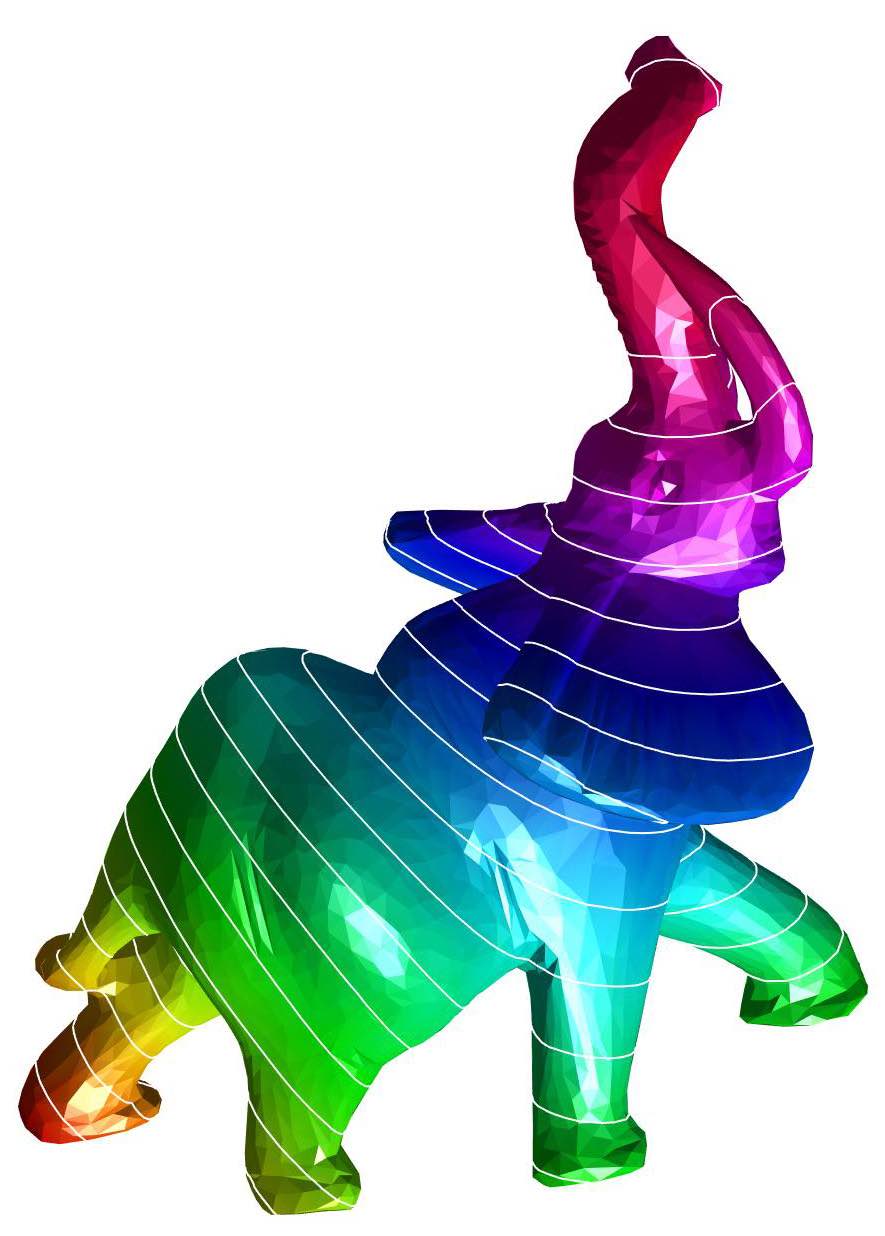}
&\includegraphics[height=90pt]{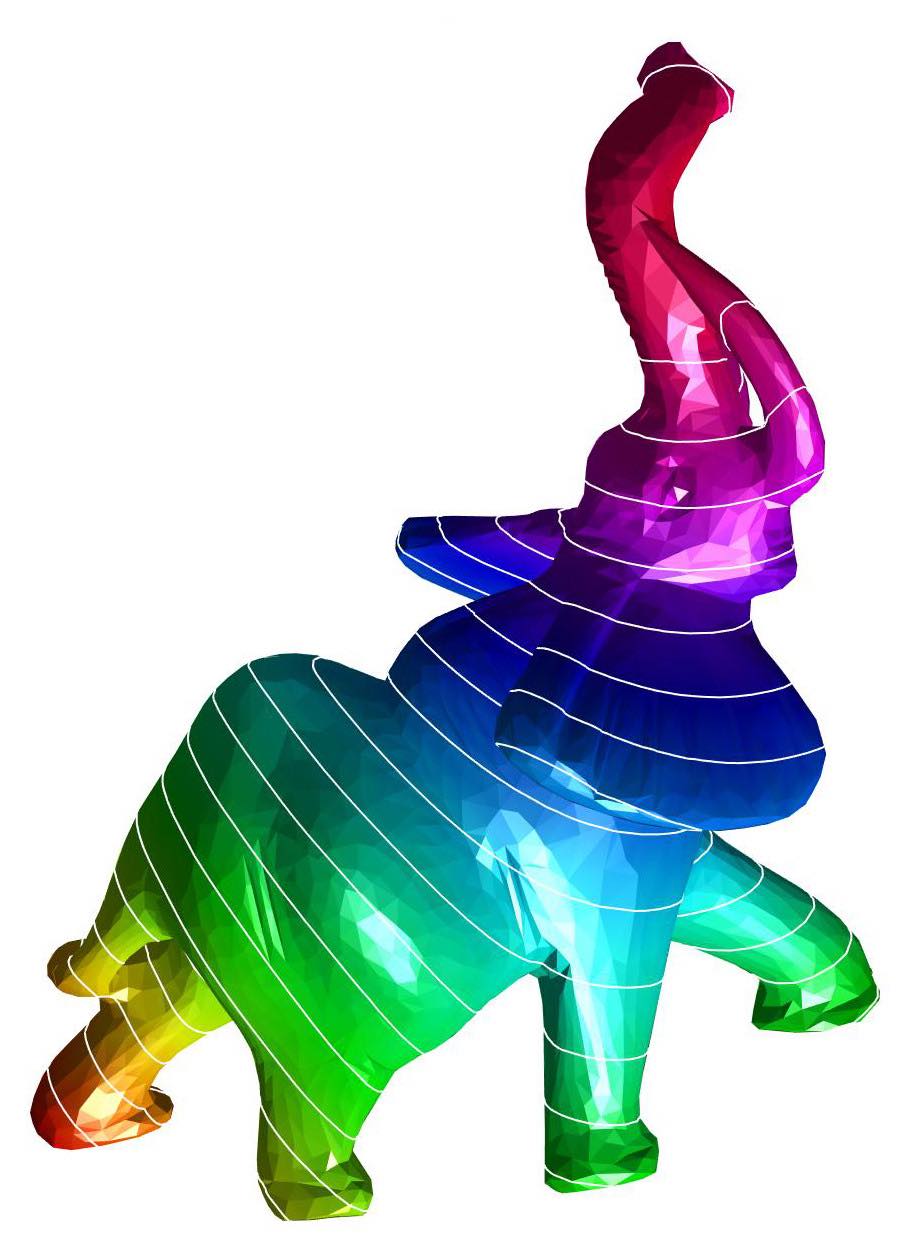}
&\includegraphics[height=90pt]{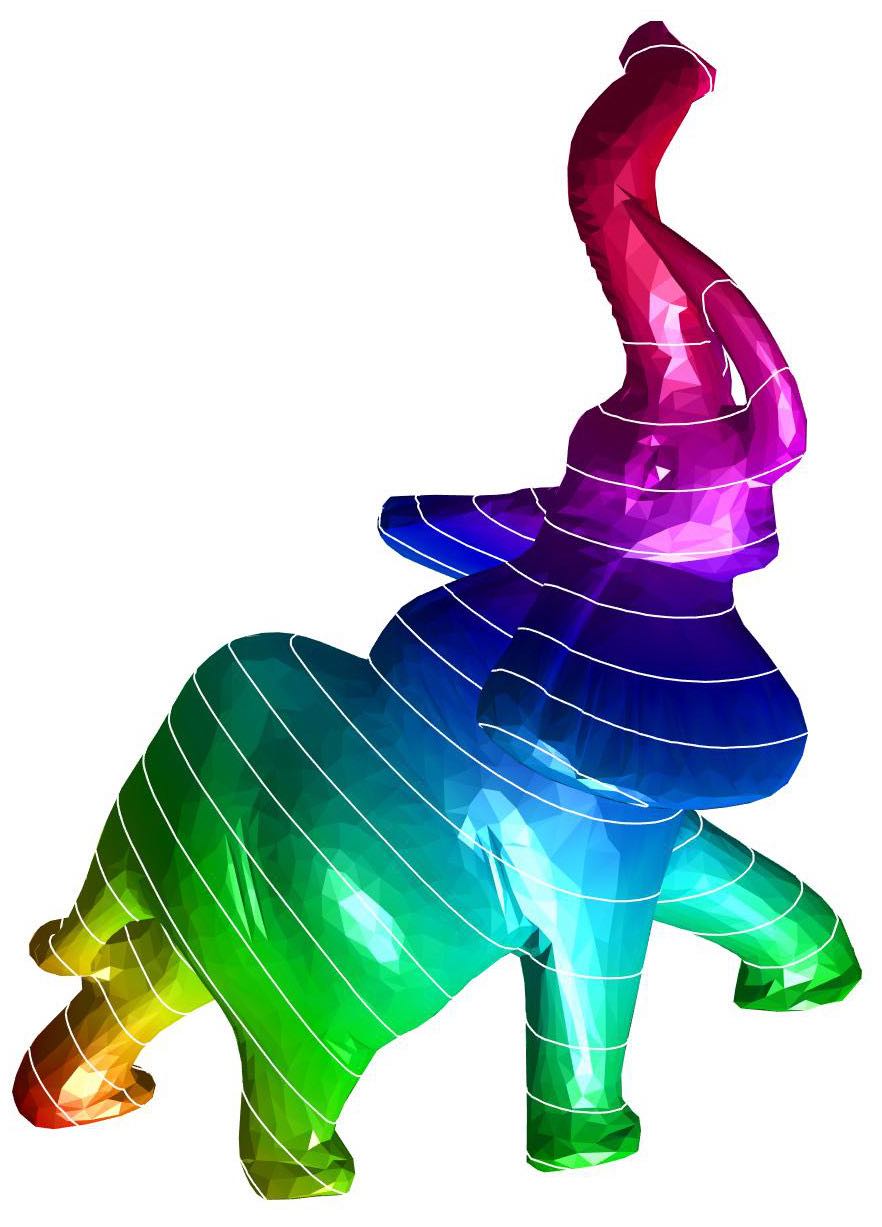}
&\includegraphics[height=90pt]{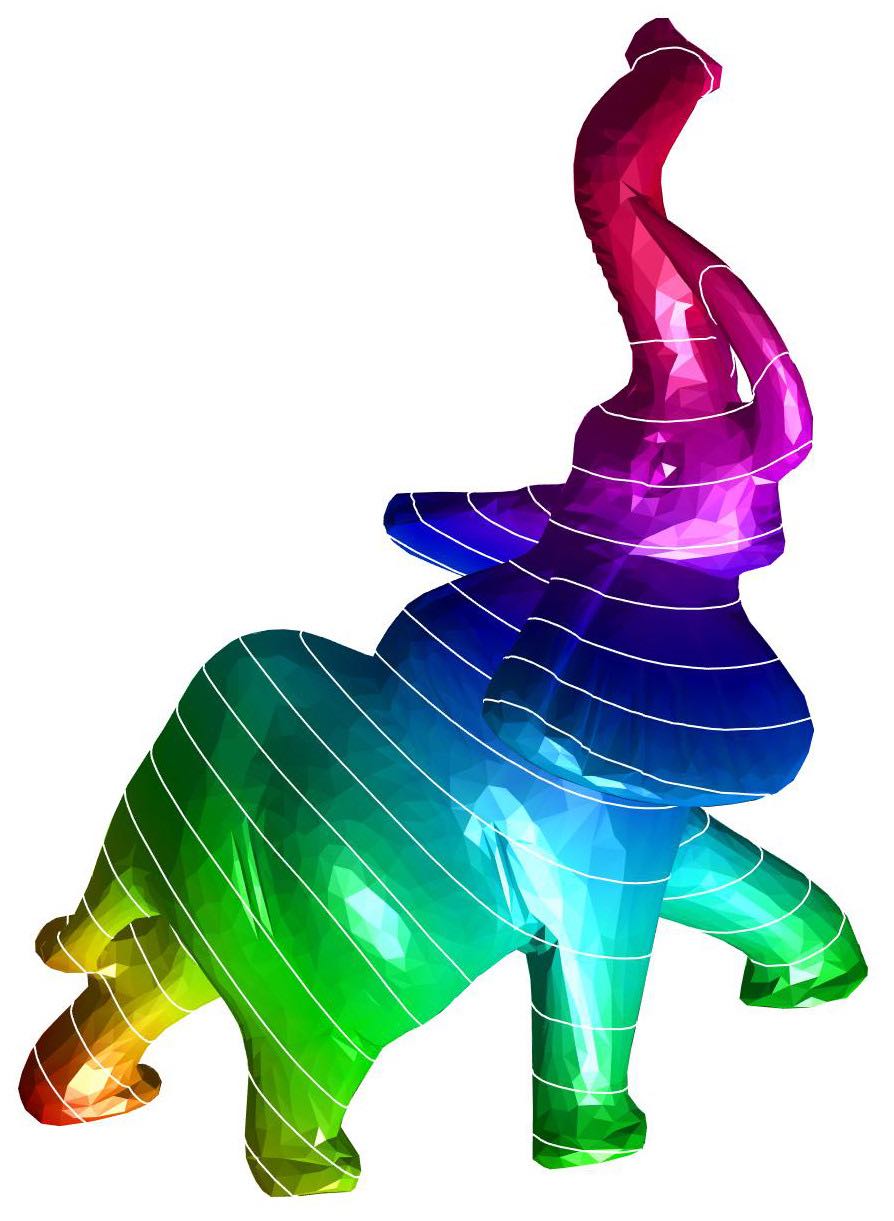}
&\includegraphics[height=90pt]{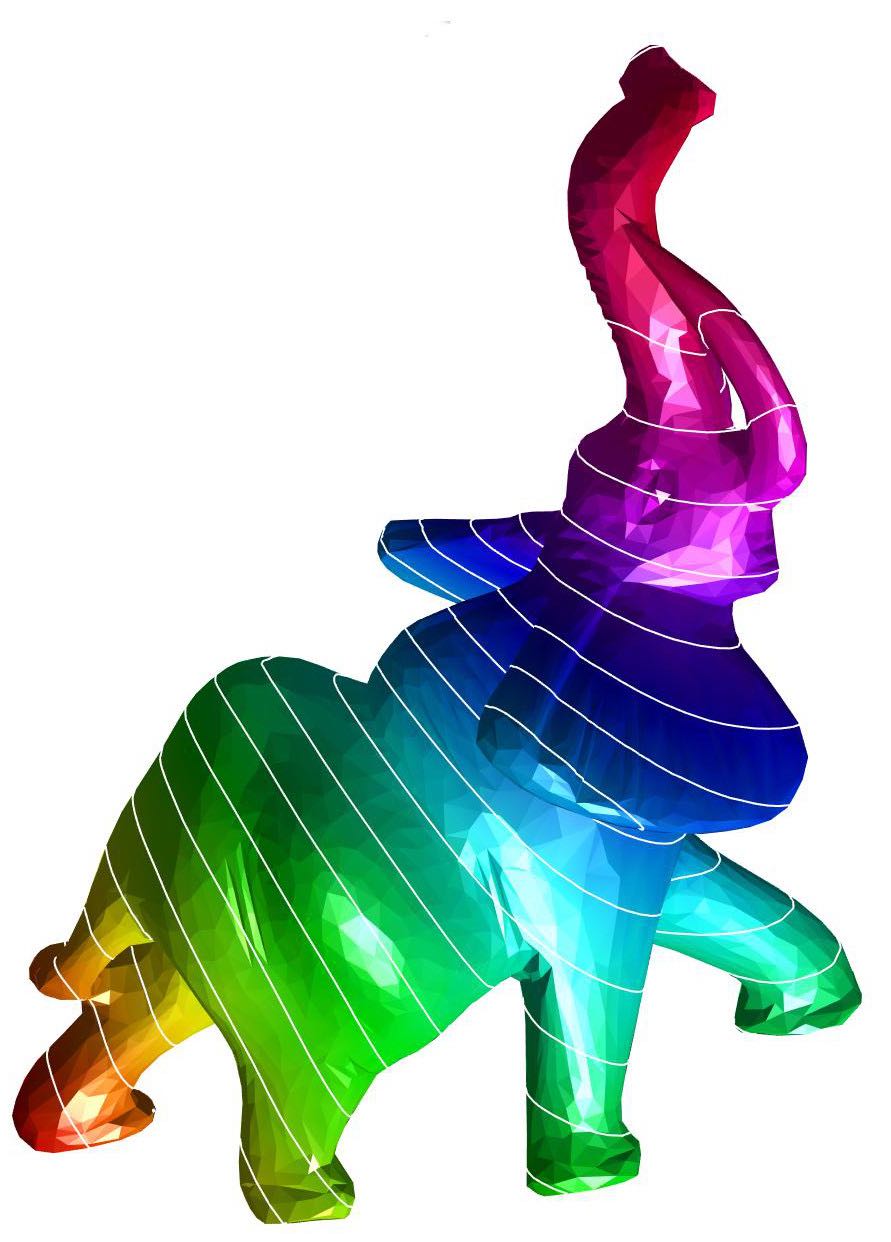}\\
\hline
$\mathcal{F}^{-1}(\mathcal{F}f)$&\includegraphics[height=90pt]{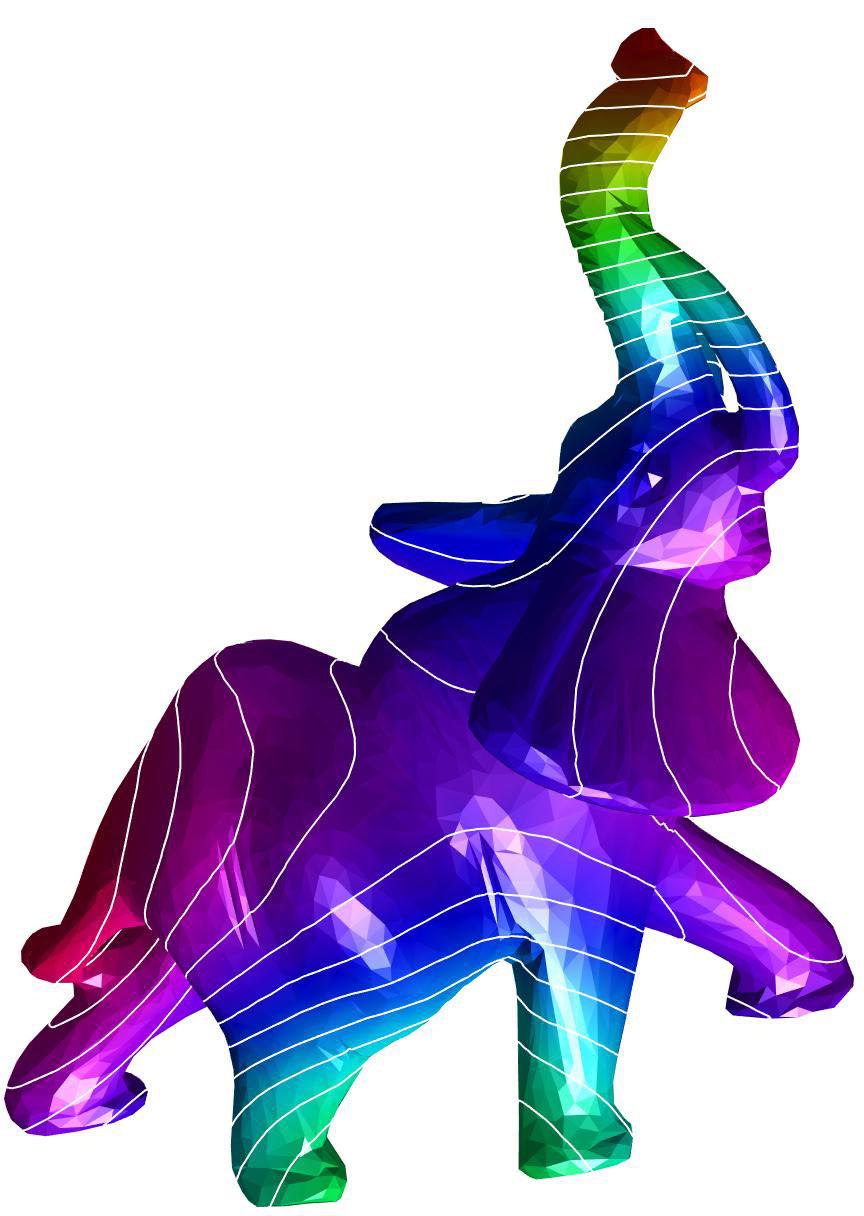}
&\includegraphics[height=90pt]{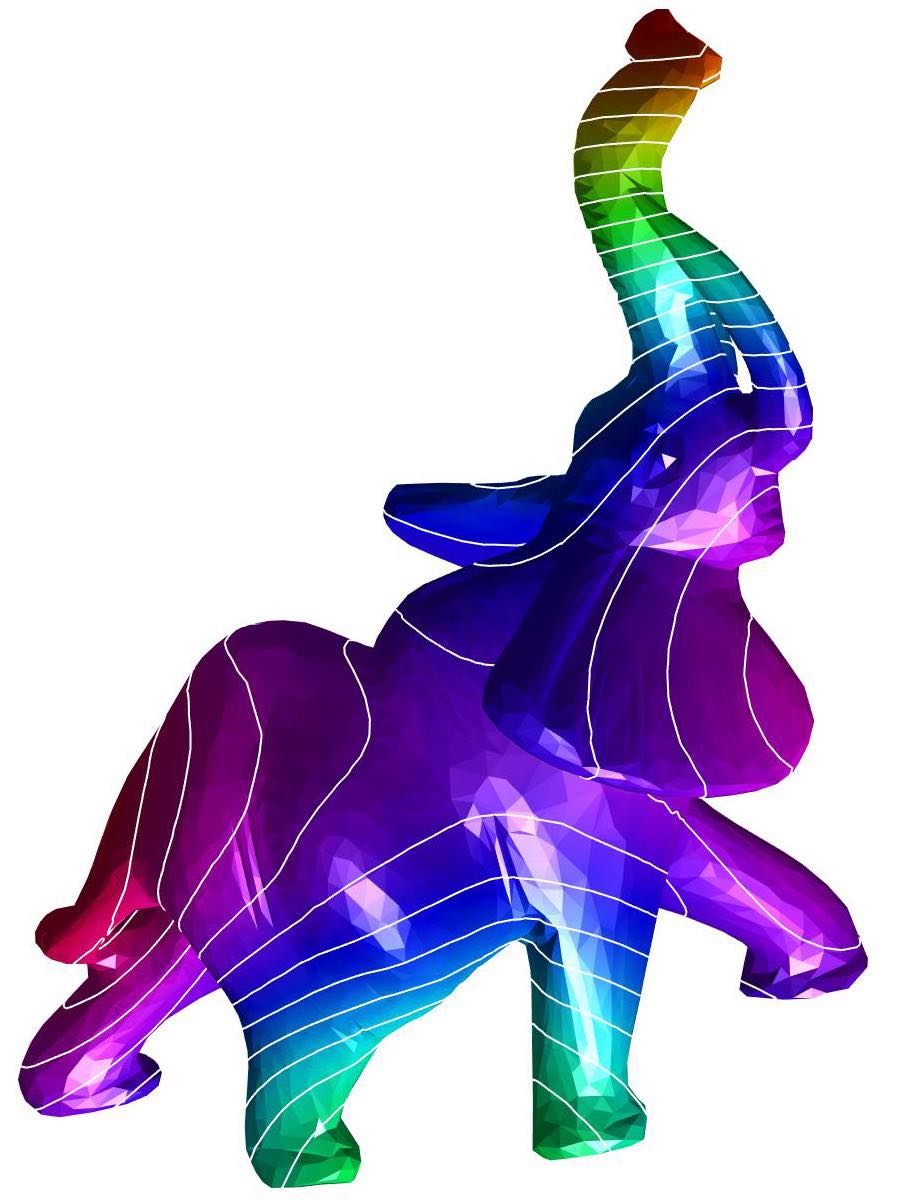}
&\includegraphics[height=90pt]{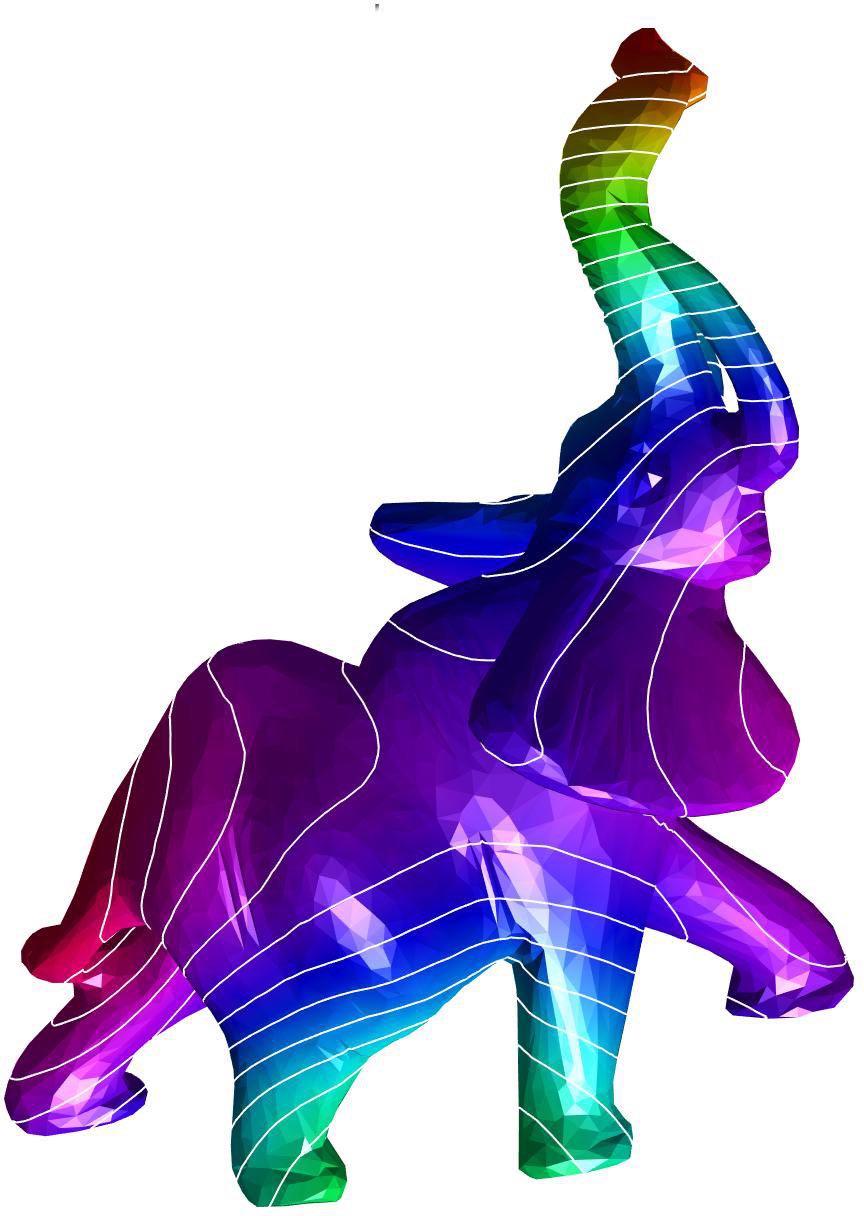}
&\includegraphics[height=90pt]{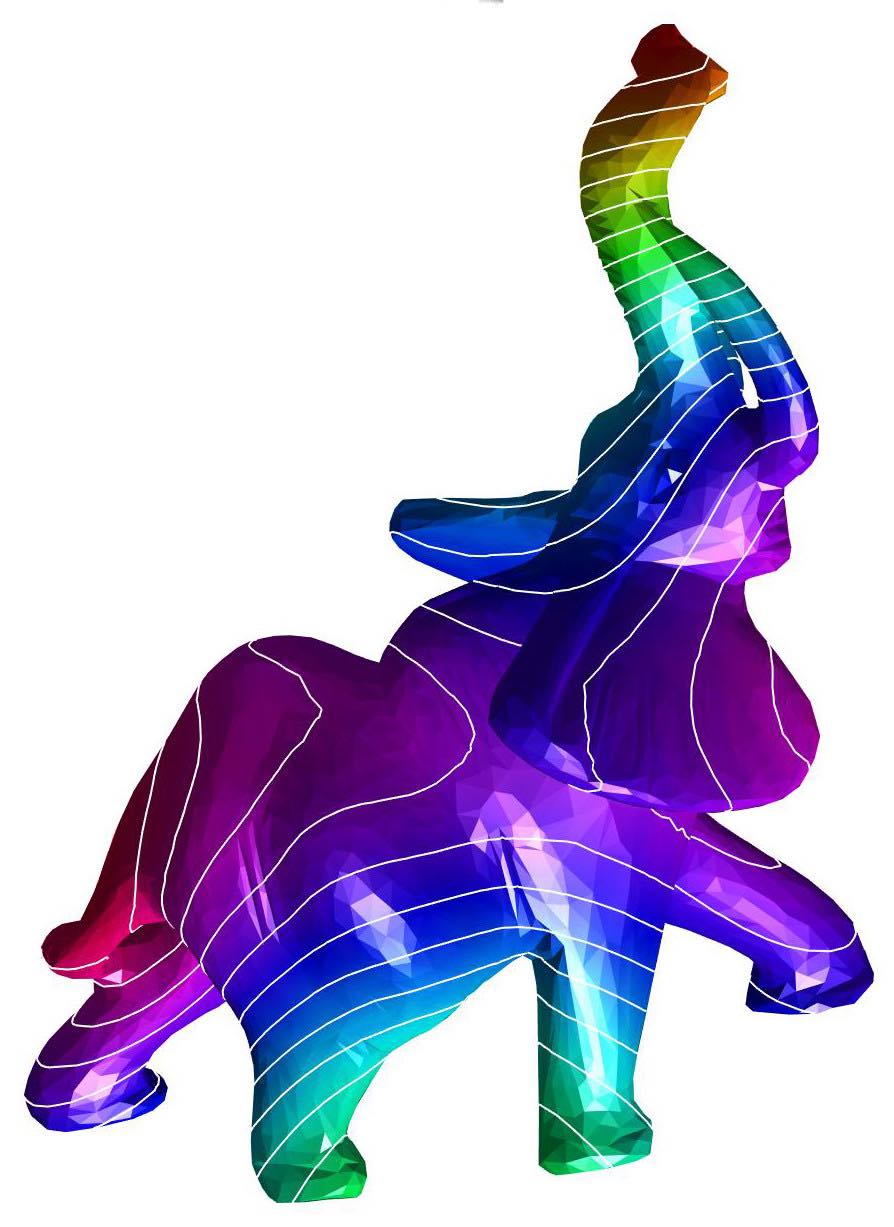}
&\includegraphics[height=90pt]{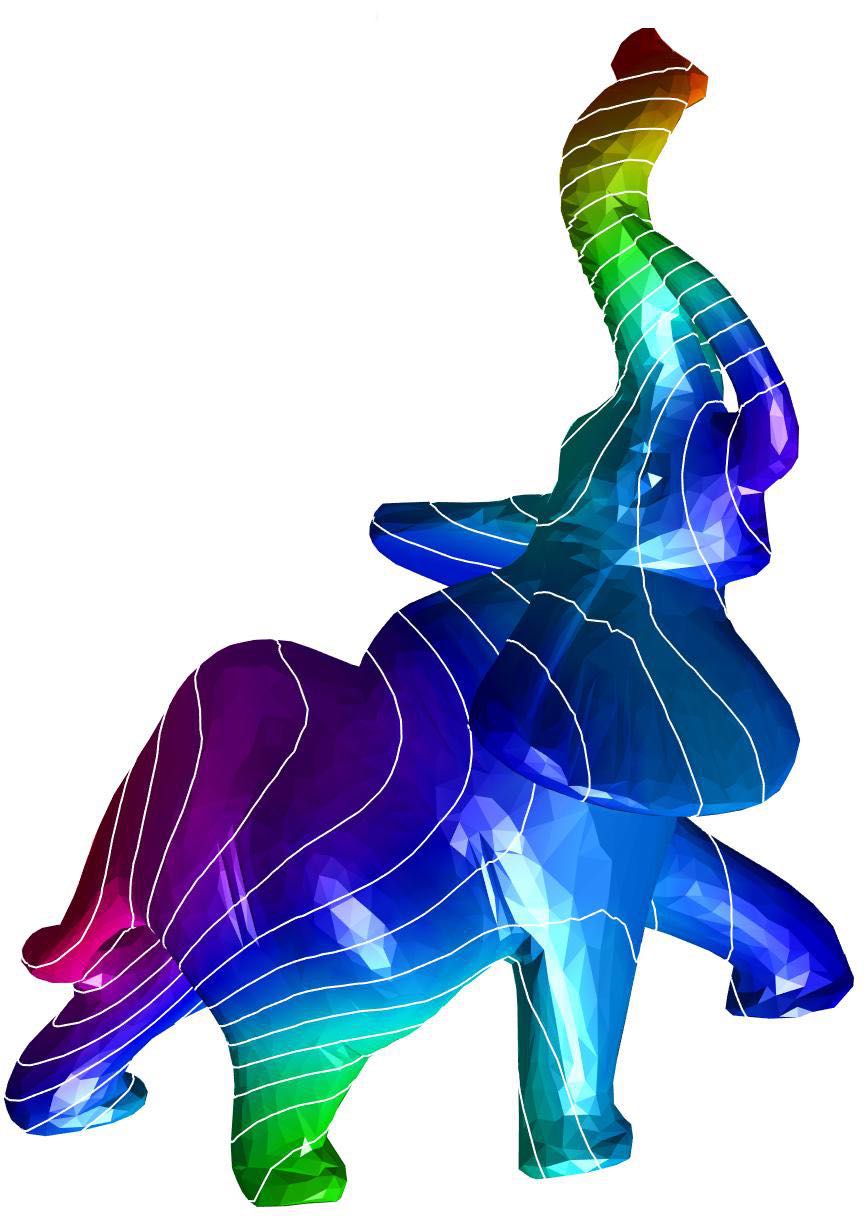}\\
&$\epsilon_{\infty}=1.2\%$	&$\epsilon_{\infty}=2.3\%$	&$\epsilon_{\infty}=4.1\%$	&$\epsilon_{\infty}=5.6\%$ &$\epsilon_{\infty}=7.2\%$
\end{tabular}
\caption{Level-sets (first row) of noisy signals on a 3D shape with an increasing error magnitude~$\alpha$ (from (b) to (f)) achieved by adding a Gaussian noise to (a) an input signal~$f$. Level-sets (second row) of the F-transform \mbox{$\mathcal{F}f$} and (third row) of the reconstructed signal \mbox{$\mathcal{F}^{-1}(\mathcal{F}f)$}. Here, the F-transform and its inverse are induced by the multi-quadratic kernel.\label{fig:3D-FT-IFT-KERNELS-NOISE}}
\end{figure*}
\subsection{Discretisation of the continuous F-transform\label{sec:DISCRETISATION}}
Let us assume that the input domain~$\Omega$ is discretised as triangle mesh or a point set with~$n$ nodes, and that the Laplace-Beltrami operator is discretised as the \mbox{$n\times n$} matrix~$\mbox{$\mathbf{L}:=\mathbf{B}^{-1}\tilde{\mathbf{L}}$}$, where the \emph{mass matrix}~$\mathbf{B}$ is sparse, symmetric and the \emph{stiffness matrix}~$\tilde{\mathbf{L}}$ is sparse, symmetric, and positive semi-definite. Then, the \emph{Laplacian eigensystem} \mbox{$(\lambda_{i},\mathbf{x}_{i})_{i=1}^{n}$}, with \mbox{$\lambda_{i}\leq\lambda_{i+1}$}, satisfies the identity \mbox{$\mathbf{L}\mathbf{x}_{i}=\lambda_{i}\mathbf{x}_{i}$} and the eigenvectors \mbox{$(\mathbf{x}_{i})_{i=1}^{n}$} are orthonormal. On a triangle mesh,~\cite{PATANE2014} the stiffness matrix encodes the variation of the cotangent of the angles of the input triangles and the mass matrix encodes their area. On a point set~\cite{SCHOELKOPF02}, the stiffness matrix is
\begin{equation*}\label{eq:EXP-WEIGHTS}
L(i,j):=
\frac{1}{nt(4\pi t)^{3/2}}
\left\{
\begin{array}{ll}
\exp\left(-\frac{\|\mathbf{p}_{i}-\mathbf{p}_{j}\|_{2}}{4t}\right)				&i\neq j,\\
-\sum_{k\neq i}\exp\left(-\frac{\|\mathbf{p}_{i}-\mathbf{p}_{k}\|_{2}}{4t}\right) &i=j.
\end{array}
\right.
\end{equation*}
To guarantee the sparsity of the Laplacian matrix, for each point~$\mathbf{p}_{i}$ we consider only the entries \mbox{$L(i,j)$} related to the points \mbox{$(\mathbf{p}_{j})_{j\in\mathcal{N}_{\mathbf{p}_{i}}}$} that are closest to~$\mathbf{p}_{i}$ with respect to the Euclidean distance. Finally,~$\mathbf{B}$ is the identity matrix or the diagonal matrix whose non-null entries are the areas of the approximated Voronoi regions associated with the input points.

\textbf{Continuous F-transform and its inverse}
Evaluating the continuous F-transform at~$\mathbf{p}_{i}$, the discrete F-transform (c.f., Eq. (\ref{eq:DISCRETE-FT})) is \mbox{$\mathcal{L}_{K}f(\mathbf{p}_{i})\approx\sum_{j=1}^{s}K(\mathbf{p}_{i},\mathbf{q}_{j})f(\mathbf{q}_{j})$}, i.e., \mbox{$\mathbf{F}:=(\mathcal{L}_{K}f(\mathbf{p}_{i}))_{i=1}^{n}=\mathbf{K}\mathbf{f}$}, where \mbox{$\mathbf{K}:=(K(\mathbf{p}_{i},\mathbf{q}_{j}))_{i=1,\ldots,n}^{j=1,\ldots,s}$} is the Gram matrix associated with the input kernel, \mbox{$\mathbf{F}:=(F_{i})_{i=1}^{n}$}, and \mbox{$\mathbf{f}:=(f(\mathbf{q}_{i}))_{i=1}^{s}$}. In particular, the problem \mbox{$\mathcal{L}_{K}f=g$} is equivalent to the linear system \mbox{$\mathbf{K}\mathbf{f}=\mathbf{g}$}, and the \emph{discrete pseudo-inverse inverse F-transform} is induced by the pseudo-inverse matrix \mbox{$\mathbf{K}^{\dag}=\mathbf{V}\Sigma^{\dag}\mathbf{U}^{\top}$}, where~$\Sigma^{\dag}$ is the diagonal matrix whose entries are the reciprocal of the non-null singular values. For a large number of samples, we cannot compute the entire spectrum of~$\mathbf{K}$; indeed, we evaluate the first~$k$ singular values and vectors with \mbox{$k<<n$}.

\textbf{Computational cost\label{sec:COMPUTATIONAL-COST}}
The truncated approximation of the spectral membership functions with~$k$ number of selected eigenpairs takes from \mbox{$\mathcal{O}(kn\log n)$} to \mbox{$\mathcal{O}(kn^{2})$} time, according to the sparsity of the Laplacian matrix.

\subsection{Experimental results\label{sec:EXPERIMENTS}}
The analytic membership functions induced by the Gaussian kernel (Fig.~\ref{fig:2D-ANALYTIC-GAUSSIAN}) and the diffusion membership functions \mbox{$K_{t}(\mathbf{p},\cdot)$} (Fig.~\ref{fig:2D-DIFFUSION}) are well localised around their seed point and their support reduces as the scale tends to zero, thus showing their multi-scale behaviour with respect to the parameter~$t$. The membership functions induced by the Gaussian and hyperbolic tangent are well localised  around the seed point, similarly to the diffusion kernel, as a matter of the decay of the filter to zero. On the contrary, the membership functions induced by the multi-quadratic and inverse multi-quadratic kernels (Fig.~\ref{fig:2D-DIFFUSION-MF}), as well as the harmonic and biharmonic membership functions (Fig.~\ref{fig:2D-HAR-BIHAR}), are generally not localised around the seed point and are globally-supported. 

To discuss the accuracy and robustness to noise of the continuous F-transform and its inverse, we consider a noisy signal (Fig.~\ref{fig:2D-FT-IFT-KERNELS-NOISE}, first row) on a 2D domain with an increasing error (from (b) to (e)) achieved by adding a Gaussian noise to an input signal~$f$ (Fig.~\ref{fig:2D-FT-IFT-KERNELS-NOISE}(a)). Evaluating its inverse F-transform (Fig.~\ref{fig:2D-FT-IFT-KERNELS-NOISE}, second row) and its reconstruction (Fig.~\ref{fig:2D-FT-IFT-KERNELS-NOISE}, third row; Fig.~\ref{fig:2D-FT-IFT-KERNELS}), we notice the consistency of the behaviour of the continuous F-transform and its inverse with respect to a different level of noise. Here, the inverse F-transform is induced by the multi-quadratic kernel and the \emph{normalised reconstruction error} \mbox{$\epsilon_{\infty}:=\|f-\mathcal{F}^{-1}(\mathcal{F}f)\|_{\infty}/\|f\|_{\infty}$} is measured as~$\mathcal{L}_{\infty}$ error between~$f$ and its reconstruction \mbox{$\mathcal{F}^{-1}(\mathcal{F}f)$}. 

Finally, the generality of the proposed approach allows us to define data-driven membership functions (Fig.~\ref{fig:3D-DIFFUSION-MF}), the F-transform and its inverse (Fig.~\ref{fig:3D-FT-IFT}) on arbitrary data, in terms of connectivity and dimensionality and with guarantee on the approximation accuracy and stability to noise (Fig.~\ref{fig:3D-FT-IFT-KERNELS-NOISE}).

\section{Conclusions and future work\label{sec:CONCLUSION}}
This work has introduced the continuous F-transform, as a generalisation of the discrete F-transform, which is ubiquitous in several fields, such as fuzzy logic, fuzzy modelling, and artificial intelligence. This generalisation is based on integral operators induced by symmetric kernels, whose properties naturally extend to the F-transform to data with an arbitrary dimension and structure. As future work, we will further study the definition and properties of the data-driven F-transform, with a focus on the construction of a larger class of data-driven membership functions, which adapt to the input data in terms of sampling density and encode the underlying geometric and topological properties. Finally, we will address the efficient computation of the data-driven F-transform in terms of approximation accuracy and numerical robustness.

\subsection*{Acknowledgements}
We thank the Reviewers for their thorough review and constructive comments, which helped us to improve the technical part and presentation of the revised paper. This work is partially supported by the H2020 ERC Advanced Grant CHANGE, grant agreement N. 694515.
\begin{IEEEbiography}
[{\includegraphics[width=1in,height=1.10in,clip,keepaspectratio]{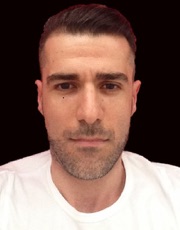}}]{Giuseppe Patan\`e}
is senior researcher at CNR-IMATI. Since 2001, his research is mainly focused on Data Science. He obtained the National Scientific Qualification as Full Professor of Computer Science. He is author of scientific publications on international journals and conference proceedings, and tutor of Ph.D. and Post.Doc students. He is responsible of R$\&$D activities in national and European projects.
\end{IEEEbiography}
\end{document}